\documentclass[preprints,article,accept,moreauthors,pdftex,10pt,a4paper]{mdpi}
\firstpage{1} 
\pubvolume{xx}
\issuenum{1}
\articlenumber{1}
\pubyear{2018}
\copyrightyear{2018}
\externaleditor{Academic Editor: Luigi Brugnano}
\history{Received: date; Accepted: date; Published: date}

\renewcommand{\S}{\mathbf{S}} 
\newcommand{\Q}{\mathbf{Q}} 
\newcommand{\E}{\mathbf{E}} 
\renewcommand{\u}{\mathbf{u}} 
\newcommand{\F}{\mathbf{F}} 
\renewcommand{\G}{\mathbf{G}} 
\newcommand{\w}{\mathbf{w}} 
\newcommand{\q}{\mathbf{q}} 
\newcommand{\x}{\mathbf{x}} 
\newcommand{\dev}{\textnormal{dev}} 
\renewcommand{\d}{\, d} 
\newcommand{\halb}{\frac{1}{2}} 
\newcommand{\diff}[2]{\frac{\partial {#1} }{\partial {#2} } }

\Title{Efficient implementation of ADER discontinuous Galerkin schemes for a scalable hyperbolic PDE engine}

\Author{Michael Dumbser $^{1,\dagger,*}$\orcidA{}, Francesco Fambri $^{1}$\orcidB{}, Maurizio Tavelli $^{1}$, 
	\\ Michael Bader $^{2}$  and Tobias Weinzierl $^{3}$}

\AuthorNames{Michael Dumbser, Francesco Fambri, Maurizio Tavelli, Michael Bader and Tobias Weinzierl}

\address{%
$^{1}$ \quad University of Trento, Italy; michael.dumbser@unitn.it, francesco.fambri@unitn.it   \\
$^{2}$ \quad Technical University Munich, Germany; bader@in.tum.de   \\
$^{3}$ \quad University of Durham, United Kingdom; tobias.weinzierl@durham.ac.uk   \\
}

\corres{Correspondence: michael.dumbser@unitn.it}

\firstnote{Current address: University of Trento, Via Mesiano 77, I-38123 Trento, Italy} 

\abstract{In this paper we discuss a new and very efficient implementation of high order accurate ADER discontinuous Galerkin (ADER-DG) finite element schemes on modern massively parallel supercomputers. The numerical methods apply to a very broad  
	class of nonlinear systems of hyperbolic partial differential equations. 
ADER-DG schemes are by construction communication avoiding and cache blocking and are furthermore very well-suited for  vectorization, so that they 
appear to be a good candidate for the future generation of exascale supercomputers. We introduce the 
numerical algorithm and show some applications to a set of hyperbolic equations with increasing
level of complexity, ranging from the compressible Euler equations over the equations 
of linear elasticity and the unified Godunov-Peshkov-Romenski (GPR) model of continuum mechanics 
to general relativistic magnetohydrodynamics (GRMHD) and the Einstein field 
equations of general relativity.  
We present strong scaling results of the new ADER-DG schemes up to 180,000 CPU cores. To our
knowledge, these are the largest runs ever carried out with high order ADER-DG schemes
for nonlinear hyperbolic PDE systems. We also provide
a detailed performance comparison with traditional Runge-Kutta DG schemes. }


\begin{document}


\section{Introduction}

Hyperbolic partial differential equations are omnipresent in the mathematical description of 
time-dependent processes in fluid and solid mechanics, in engineering and geophysics, as well 
as in plasma physics and even in general relativity. Among the most widespread applications 
nowadays are i) computational fluid mechanics in mechanical and aerospace engineering, in particular 
compressible gas dynamics at high Mach numbers; ii) geophysical and environmental free surface flows 
in oceans, rivers and lakes, describing natural hazards such as tsunami wave propagation,  
landslides, storm surges and floods; iii) seismic, acoustic and electromagnetic wave propagation processes
in the time domain are described by systems of hyperbolic partial differential equations, namely the 
equations of linear elasticity, the acoustic wave equation and the well-known Maxwell equations; 
iv) high energy density plasma flows in nuclear fusion reactors as well as astrophysical plasma flows 
in the solar system and the universe, using either the Newtonian limit or the complete equations
in full general relativity; v) the Einstein field equations of general relativity, which govern the 
evolution of dynamic spacetimes, can be written under the form of a nonlinear system of hyperbolic 
partial differential equations.  

The main challenge of nonlinear hyperbolic PDE arises from the fact that they can contain at the same
time smooth solutions (like sound waves) as well as small scale structures (e.g. turbulent vortices), 
but they can also develop discontinuous solutions (shock waves) after finite times, even when starting 
from perfectly smooth initial data. 
These discontinuities were first discovered by Bernhard Riemann in his ground breaking work on the 
propagation of waves of finite amplitude in air \cite{riemann1,riemann2}, where the term 
\textit{finite} should actually be understood in the sense of \textit{large}, rather than simple sound 
waves of \textit{infinitesimal} strength that have been considered in the times before Riemann. In the 
abstract of his work, Riemann stated that his discovery of the shock waves might probably not be of 
practical use for applied and experimental science, but should be mainly understood as a contribution 
to the theory of nonlinear partial differential equations. Several decades later, shock waves were also 
observed experimentally, thus confirming the new and groundbreaking mathematical concept of Riemann.    

The connection between symmetries and conservation laws were established in the work of Emmy Noether
\cite{noether} at the beginning of the 20th century, while the first methods for the 
\textit{numerical solution} of hyperbolic conservation laws go back to famous mathematicians such as 
Courant and Friedrichs and co-workers \cite{CFL,CIR,Courant62,Courant62b}. The connection between
hyperbolic conservation laws, symmetric hyperbolic systems in the sense of Friedrichs \cite{FriedSymm}  
and thermodynamics was established for the first time by Godunov in 1961 \cite{God1961} and was 
rediscovered again by Friedrichs and Lax in 1971 \cite{FriedLax1971}. Within this theoretical framework of 
symmetric hyperbolic and thermodynamically compatible (SHTC) systems, established by Godunov and Romenski
\cite{Godunov:1995a,Godunov:2003a} it is possible to write down the Euler equations of compressible gas dynamics,  
the magnetohydrodynamics (MHD) equations \cite{God1972MHD}, the equations of nonlinear elasticity  \cite{GodunovRomenski72}, as well as a rather wide class of nonlinear hyperbolic conservation laws 
\cite{Rom1998} with very interesting mathematical properties and structure.  
Very recently, even a novel and unified formulation of continuum physics, including solid and fluid 
mechanics only as two particular cases of a more general model, have been 
cast into the form of a single SHTC system \cite{PeshRom2014,HPRmodel,HPRmodelMHD,LagrangeHPR}.         
In the 1940ies and 1950ies major steps forward in numerical methods for hyperbolic PDE have been made in 
the ground-breaking contributions  of von Neumann and Richtmyer 
\cite{vonneumann50} and Godunov \cite{godunov:linear}. While the former introduce an artificial viscosity 
to stabilize the numerical scheme in the presence of discontinuities, the latter constructs his scheme
starting from the most elementary problem in hyperbolic conservation laws for which an exact solution 
is still available, the so-called Riemann problem. The Riemann problem consists in a particular Cauchy 
problem where the initial data consist of two piecewise constant states, separated by a discontinuity.  
In the absence of source terms, its solution is self--similar. While provably robust, these schemes are 
only first order accurate in space and time and thus only applicable to flows with shock waves, but
not to those also involving smooth sound waves and turbulent small scale flow structures. In his paper 
\cite{godunov:linear}, Godunov has also proven that any linear numerical scheme that is required to  
be monotone can be at most of order one, which is the well-known Godunov barrier theorem. The 
main goal in the past decades was to find ways how to circumvent it, since it only applies to linear schemes. 
The first successful nonlinear monotone and higher order accurate schemes were the method of Kolgan 
\cite{Kolgan:1972} and the schemes of van Leer \cite{leer2,leer5}. Subsequently, many other higher
order nonlinear schemes have been proposed, such as the ENO \cite{HartenENO} and WENO schemes 
\cite{shu_efficient_weno} and there is a rapidly growing literature on the subject. In this paper we
mainly focus on a rather recent family of schemes, which is of the discontinuous finite element type,
namely the so-called discontinuous Galerkin (DG) finite element methods, which were systematically 
introduced for hyperbolic conservation laws in a well-known series of papers by Cockburn and Shu and 
collaborators \cite{cbs1,cbs2,cbs3,cbs4,CBS-convection-diffusion}. For a review on high order DG methods
and WENO schemes, the reader is referred 
to \cite{CBS-convection-dominated} and \cite{ShuDGWENOReview}. 
In this paper we use a particular variant of the DG scheme that is called
ADER discontinuous Galerkin scheme \cite{Dumbser2008,ADERNC,DGLimiter1,DGLimiter2,DGLimiter3}, 
where ADER stands for arbitrary high order schemes using derivatives, first developed by Toro et al. 
in the context of high order finite volume schemes \cite{toro3,toro4,titarevtoro,Toro:2006a}.  
In comparison to traditional \textit{semi-discrete} DG schemes, which mainly use Runge-Kutta time integration, 
ADER-DG methods are \textit{fully-discrete} and are based on a predictor-corrector approach that allows to achieve 
a naturally cache-blocking and communication-avoiding scheme, which reduces the amount of necessary MPI communications 
to a minimum.  
These properties make the method well suitable for high performance computing (HPC).

\section{High order ADER discontinuous Galerkin finite element schemes}

In this paper we consider hyperbolic PDE with non-conservative products and algebraic source
terms of the form (see also \cite{Dumbser2008,ADERNC}) 
\begin{equation}
\label{eqn.pde}
\frac{\partial \mathbf{Q} }{\partial t} + \nabla \cdot \mathbf{F}\left( \mathbf{Q} \right) 
+ \mathbf{B}(\Q) \cdot \nabla \mathbf{Q}  = \mathbf{S}(\mathbf{Q}),
\end{equation}
where $t \in \mathbb{R}_0^+$ is the time, $\mathbf{x} \in \Omega \subset \mathbb{R}^d$ is the spatial position vector in $d$ space dimensions, 
$\Q \in \Omega_{Q} \subset \mathbb{R}^m$ is the state vector, $\F(\Q)$ is the nonlinear flux 
tensor, $\mathbf{B}(\Q) \cdot \nabla \mathbf{Q} $ is a non-conservative product and $\S(\Q)$ is 
a purely algebraic source term. Introducing the system matrix  $\mathbf{A}(\Q) = \partial \F / \partial \Q + \mathbf{B}(\Q)$ 
the above system can also be written in quasi-linear form as 
\begin{equation}
\frac{\partial \mathbf{Q} }{\partial t}   
+ \mathbf{A}(\Q) \cdot \nabla \mathbf{Q}  = \mathbf{S}(\mathbf{Q}).
\end{equation}
The system is hyperbolic if for all $\mathbf{n} \neq 0$ and for all $\Q \in \Omega_Q$ the 
matrix $\mathbf{A}(\Q) \cdot \mathbf{n}$ has $m$ real eigenvalues and a full set of 
$m$ linearly independent right eigenvectors. The system \eqref{eqn.pde} is provided with an  
initial condition $\Q(\x,0) = \Q_0(\x)$ and appropriate boundary conditions on $\partial \Omega$. 
\textcolor{black}{In some parts of the paper we will also make use of the vector of primitive (physical) 
variables denoted by $\mathbf{V}=\mathbf{V}(\mathbf{Q})$. For very complex PDE systems, like the general 
relativistic MHD equations, it may be much easier to express the flux tensor $\mathbf{F}$ in terms of $\mathbf{V}$ 
rather than in terms of $\mathbf{Q}$, but the evaluation of $\mathbf{V}=\mathbf{V}(\mathbf{Q})$ can become 
very complicated. } 

\subsection{Unlimited ADER-DG scheme and Riemann solvers}

We cover the computational domain $\Omega$ with a set of non-overlapping
Cartesian control volumes in space $\Omega_{i} = [x_i - \halb \Delta x_i,
x_i + \halb \Delta x_i] \times [y_i - \halb \Delta y_i, y_i + \halb
\Delta y_i] \times [z_i - \halb \Delta z_i, z_i + \halb \Delta z_i] $.
Here, $\mathbf{x}_i = (x_i, y_i, z_i)$ denotes the barycenter of cell
$\Omega_i$ and $\Delta \mathbf{x}_i = (\Delta x_i,\Delta y_i,\Delta z_i)$
is the mesh spacing associated with $\Omega_i$ in each space dimension.  
The domain $\Omega = \bigcup \Omega_i$ is the union of all spatial
control volumes. A key ingredient of the ExaHyPE engine \textcolor{blue}{http://exahype.eu} 
is a cell-by-cell adaptive mesh refinement (AMR), which is built upon the space-tree 
implementation Peano \cite{Peano1,Peano2}. For further details about cell-by-cell AMR, see 
\cite{Khokhlov1998}. For AMR in combination with better than second order accurate 
finite volume and DG schemes with time-accurate local time stepping (LTS) \textcolor{black}{and
for a literature overview of state of the art AMR methods}, the 
reader is referred to \cite{AMR3DCL,AMR3DNC,Zanotti2015c,Zanotti2015d,ADERDGVisc,ADERGRMHD}
\textcolor{black}{and references therein. Since the main focus of this paper is not on AMR,  
	at this point we can only give a very brief summary of  
existing AMR methods and codes for hyperbolic PDE, without pretending to be complete. 
Starting point of adaptive mesh refinement for hyperbolic conservation laws was of course the pioneering work of Berger et al.  \cite{Berger-Oliger1984,Berger-Jameson1985,Berger-Colella1989}, 
who were the first to introduce a patched-based block-structured AMR method. 
Further developments are reported in \cite{LeVequeCLAWPACK,BergerLeveque1998,Bell1994} based on the second order 
accurate  wave-propagation algorithm of LeVeque. 
For computational astrophysics, relevant AMR techniques have been documented, e.g., in \cite{Dezeeuw1993,BalsaraAMR,Teyssier2002,Keppens2003,Ziegler2008,Mignone2012,Cunningham2009,Keppens2012,BHAC}, including the 
RAMSES, PLUTO, NIRVANA, AMRVAC and BHAC codes. For a recent and more complete survey of high level AMR codes, 
the reader is referred to the review paper \cite{AMRSurvey}. 
Better than second order accurate finite volume and finite difference schemes on adaptive meshes can be found, e.g., 
in  \cite{Mulet1,Colella2009,Burger2012,Ivan2014,Buchmuller2015,SCR:CWENOquadtree,FDWENOAMR}. }

In the following, the discrete solution of the PDE system \eqref{eqn.pde} is 
denoted by $\mathbf{u}_h$ and is defined in terms of tensor
products of piecewise polynomials of degree $N$ in each spatial
direction. The discrete solution space is denoted by $\mathcal{U}_h$ in the following. 
Since we adopt a discontinuous Galerkin (DG) finite element method, the numerical 
solution $\u_h$ is allowed to \textit{jump} across element interfaces, as in the
context of finite volume schemes. 
\textcolor{black}{Within each spatial control volume $\Omega_i$ the discrete solution $\u_h$ restricted to that
control volume} is written at time $t^n$  
in terms of some nodal spatial basis functions $\Phi_l(\mathbf{x})$ and some unknown degrees of 
freedom $\hat{\mathbf{u}}_{i,l}^n$:  
\begin{equation}
\left. \mathbf{u}_h(\mathbf{x},t^n) \right|_{\Omega_i} = \sum \limits_l \hat{\mathbf{u}}_{i,l}
\Phi_l(\mathbf{x}) := \hat{\mathbf{u}}_{i,l}^n \Phi_l(\mathbf{x})\,,
\label{eqn.ansatz.uh}
\end{equation}
where $l=(l_1,l_2,l_3)$ is a multi-index and the spatial basis functions 
$\Phi_l(\mathbf{x}) = \varphi_{l_1}(\xi) \varphi_{l_2}(\eta)
\varphi_{l_3}(\zeta)$ are generated via tensor products of 
one-dimensional nodal basis functions $\varphi_{k}(\xi)$ on the reference
interval $[0,1]$. The transformation from physical coordinates $\mathbf{x} \in
\Omega_i$ to reference coordinates $\boldsymbol{\xi}=\left(
\xi,\eta,\zeta \right) \in [0,1]^d$ is given by the linear mapping 
$\mathbf{x} =
\mathbf{x}_i - \halb \Delta \mathbf{x}_i + (\xi \Delta x_i, \eta \Delta
y_i, \zeta \Delta z_i)^T$. For the one-dimensional basis functions
$\varphi_k(\xi)$ we use the Lagrange interpolation polynomials passing
through the Gauss-Legendre quadrature nodes $\xi_j$ of an $N+1$ point
Gauss quadrature formula. Therefore, the nodal basis functions satisfy the 
interpolation property $\varphi_k(\xi_j) = \delta_{kj}$, where
$\delta_{kj}$ is the usual Kronecker symbol, and the resulting basis
is \textit{orthogonal}. Furthermore, due to this particular 
choice of a \textit{nodal} tensor-product basis, the entire scheme can be
written in a dimension-by-dimension fashion, where all integral operators
can be decomposed into a sequence of one-dimensional operators acting
only on the $N+1$ degrees of freedom in the respective dimension.
For details on multi-dimensional quadrature, see the well-known book
of Stroud \cite{stroud}. 

In order to derive the ADER-DG method, we first multiply the governing PDE
system \eqref{eqn.pde} with a test function $\Phi_k \in \mathcal{U}_h$
and integrate over the space-time control volume $\Omega_i \times
[t^n;t^{n+1}]$. This leads to  
\begin{equation}
\label{eqn.pde.nc.gw1}
\int \limits_{t^n}^{t^{n+1}} \int \limits_{\Omega_i}
\Phi_k \frac{\partial \Q}{\partial t} \d\x \d t
+\int \limits_{t^n}^{t^{n+1}} \int \limits_{\Omega_i}
\Phi_k \left( \nabla \cdot \F(\Q) + \mathbf{B}(\Q) \cdot \nabla \Q  \right) \d\x \d t
= \int \limits_{t^n}^{t^{n+1}} \int \limits_{\Omega_i}   \Phi_k
\mathbf{S}(\Q) \d\x \d t\,,
\end{equation}
with $\d \mathbf{x} = \d x\,\d y\,\d z$.  As already mentioned before, the
discrete solution is allowed to jump across element interfaces, which means
that the resulting jump terms have to be taken properly into account. In 
our scheme this is achieved via numerical flux functions (approximate Riemann solvers) 
and via the path-conservative approach that was developed by Castro and Par\'es in the 
finite volume context \cite{Castro2006,Pares2006}. It has later been also extended  
to the discontinuous Galerkin finite element framework in 
\cite{Rhebergen2008,ADERNC,USFORCE2}. In classical Runge-Kutta DG schemes, only a weak
form in space of the PDE is obtained, while time is still kept continuous, thus reducing 
the problem to a nonlinear system of ODE, which is subsequently integrated with standard
ODE solvers in time. However, this requires MPI communication in each Runge-Kutta
stage. Furthermore, each Runge-Kutta stage requires accesses to the entire discrete 
solution in memory. In the ADER-DG framework, a completely different paradigm is used. 
Here, higher order in time is achieved with the use of an element-local 
space-time predictor, denoted by $\mathbf{q}_h(\mathbf{x},t)$ in the following, and which 
will be discussed in more detail later. Using \eqref{eqn.ansatz.uh}, integrating the 
first term by parts in time and integrating the flux divergence term by parts in space,
 taking into account the jumps between elements and making use of this local space-time 
 predictor solution $\mathbf{q}_h$ instead of $\Q$, the weak formulation 
 \eqref{eqn.pde.nc.gw1} can be rewritten as
\begin{eqnarray}
\label{eqn.pde.nc.gw2}
\left( \int \limits_{\Omega_i}  \Phi_k \Phi_l \d\x \right)
\left( \hat{\mathbf{u}}^{n+1}_{i,l} - \hat{\mathbf{u}}^{n}_{i,l}  \right)
+ \int \limits_{t^n}^{t^{n+1}} \! \! \int \limits_{\partial \Omega_i}
\Phi_k \mathcal{D}^- \left( \q_h^-, \q_h^+ \right) \cdot \mathbf{n} \,\d S \d t 
- \int \limits_{t^n}^{t^{n+1}} \! \! \int \limits_{\Omega_i^\circ}
\left( \nabla \Phi_k \cdot \F(\q_h)  \right)  \d\x \d t  
+  && \nonumber \\  
 + \int \limits_{t^n}^{t^{n+1}} \! \! \int \limits_{\Omega_i^\circ}
\Phi_k \left( \mathbf{B}(\q_h) \cdot \nabla \q_h  \right)  \d\x \d t
= \int \limits_{t^n}^{t^{n+1}} \! \! \int \limits_{\Omega_i}
\Phi_k \mathbf{S}(\q_h)  \d\x \d t\,,
\end{eqnarray}
where the first integral leads to the element mass matrix,
which is diagonal since our basis is orthogonal. The boundary integral 
contains the approximate Riemann solver and accounts for the jumps across
element interfaces, also in the presence of non-conservative products. 
The third and fourth integral account for the smooth part of the flux and the
non-conservative product, while the right hand side takes into account
the presence of the algebraic source term. According to the framework
of path-conservative schemes \cite{Pares2006,Castro2006,ADERNC,USFORCE2}, the 
jump terms are defined via a path-integral in phase space between the boundary 
extrapolated states at the left $\mathbf{q}_h^-$ and at the right $\q_h^+$ of the 
interface as follows:
\begin{equation}
\label{eqn.pc.scheme}
\mathcal{D}^-\left( \q_h^-, \q_h^+ \right) \cdot \mathbf{n} =
\frac{1}{2} \left( \F(\q_h^+) + \F(\q_h^-) \right) \cdot \mathbf{n} + 
\frac{1}{2} \left( \, \int \limits_{0}^1 \mathbf{B}(\boldsymbol{\psi})
\cdot \mathbf{n} \, \d s - \boldsymbol{\Theta} \right) \left( \q_h^+ - \q_h^- \right),
\end{equation}
with $\mathbf{B} \cdot \mathbf{n} = \mathbf{B}_1 n_1 + \mathbf{B}_2 n_2 +
\mathbf{B}_3 n_3$. Throughout this paper, we use the simple straight--line segment path 
\begin{equation}
\boldsymbol{\psi} = \boldsymbol{\psi}(\q_h^-, \q_h^+, s) =
\q_h^- + s \left( \q_h^+ - \q_h^- \right),
\qquad 0 \leq s \leq 1\,.
\end{equation}
In order to achieve exactly well-balanced schemes for certain classes of hyperbolic 
equations with non-conservative products and source terms, the segment
path is not sufficient and a more elaborate choice of the path becomes necessary, 
see e.g. \cite{MuellerToro1,MuellerToro2,GaburroDumbserSWE,GaburroDumbserEuler}. 
In relation  \eqref{eqn.pc.scheme} above the symbol $\boldsymbol{\Theta} > 0$ 
denotes an appropriate numerical \textcolor{black}{dissipation} matrix. Following 
\cite{ADERNC,USFORCE,OsherNC}, the path integral that appears in 
\eqref{eqn.pc.scheme} can be simply evaluated via some sufficiently accurate
numerical quadrature formulae. We typically use a three-point Gauss-Legendre
rule in order to approximate the path-integral. For a simple path-conservative 
Rusanov-type method \cite{ADERNC,CastroPardoPares}, the \textcolor{black}{numerical dissipation} 
matrix reads 
\begin{equation}
\label{eqn.rusanov}
\boldsymbol{\Theta}_{\textnormal{Rus}} = s_{\max} \mathbf{I},
\qquad \textnormal{with} \qquad s_{\max} = \max \left( \left| \lambda(\q_h^-)
\right|, \left| \lambda(\q_h^+) \right| \right),  
\end{equation}
where $\mathbf{I}$ denotes the identity matrix and $s_{\max}$ is the maximum wave speed 
(eigenvalue $\lambda$ of matrix $\mathbf{A} \cdot \mathbf{n}$) at the element interface. 
In order to reduce numerical dissipation, one can use
better Riemann solvers, such as the Osher-type schemes proposed in \cite{OsherNC,OsherUniversal}, 
or the recent extension of the original HLLEM method of Einfeldt and Munz \cite{Einfeldt1991} 
to general conservative and non-conservative hyperbolic systems recently put forward in  
\cite{HLLEMNC}.   
The choice of the approximate Riemann solver and therefore of the viscosity matrix $\boldsymbol{\Theta}$
completes the numerical scheme \eqref{eqn.pde.nc.gw2}. In the next subsection, we shortly discuss
the computation of the element--local space-time predictor $\q_h$, which is a key ingredient
of our high order accurate and communication-avoiding ADER-DG schemes.  

\subsection{Space-time predictor and suitable initial guess}
As already mentioned previously, the element-local space-time predictor is an important 
\textit{key feature} of ADER-DG schemes and is briefly discussed in this section. 
The computation of the predictor solution $\q_h(\mathbf{x},t)$ is based on a weak formulation 
of the governing PDE system in space-time and was first introduced in \cite{DumbserEnauxToro,Dumbser2008,DumbserZanotti}. Starting from the known solution 
$\u_h(\mathbf{x},t^n)$ at time $t^n$ and following the terminology of Harten et al. \cite{eno}, 
we solve a so-called Cauchy problem \textit{in the small}, i.e. without considering the 
interaction with the neighbor elements. In the ENO scheme 
of Harten et al. \cite{eno} and in the original ADER approach of Toro and
Titarev \cite{titarevtoro,toro4,Toro:2006a} the strong differential form of the 
PDE was used, together with a combination of Taylor series expansions and the 
so-called Cauchy-Kovalewskaya procedure. The latter is very cumbersome, or becomes 
even unfeasible for very complicated nonlinear hyperbolic PDE systems, since it requires a lot of 
analytic manipulations of the governing PDE  
system, in order to replace time derivatives with known space derivatives at time $t^n$. 
This is achieved by successively differentiating the governing PDE system with 
respect to space and time and inserting the resulting terms into the 
Taylor series. For an explicit example of the Cauchy-Kovalewskaya procedure 
applied to the three-dimensional Euler equations of compressible gas
dynamics and the MHD equations, see \cite{DumbserKaeser07} and \cite{taube_jsc}. 
Instead, the local space-time discontinuous Galerkin predictor introduced in 
\cite{DumbserEnauxToro,Dumbser2008,DumbserZanotti}, requires only 
pointwise evaluations of the fluxes, source terms and non-conservative products.  
For element $\Omega_i$ the predictor solution $\q_h$ is now expanded in terms of a local space-time basis
\begin{equation}
\label{eqn.spacetime}
\left. \q_h(\mathbf{x},t) \right|_{\Omega^{st}_i} = \sum \limits_l \theta_l(\mathbf{x},t)
\hat{\mathbf{q}}^i_l := \theta_l(\mathbf{x},t) \hat{\mathbf{q}}^i_l\,,
\end{equation}
with the multi-index $l=(l_0,l_1,l_2,l_3)$ and where the space-time basis
functions $\theta_l(\mathbf{x},t) = \varphi_{l_0}(\tau)
\varphi_{l_1}(\xi) \varphi_{l_2}(\eta) \varphi_{l_3}(\zeta) $ are again
generated from the same one-dimensional nodal basis functions
$\varphi_{k}(\xi)$ as before, i.e. the Lagrange interpolation polynomials
of degree $N$ passing through $N+1$ Gauss-Legendre quadrature nodes. The
spatial mapping $\mathbf{x} = \mathbf{x}(\boldsymbol{\xi})$ is also the
same as before and the coordinate time is mapped to the reference time
$\tau \in [0,1]$ via $t = t^n + \tau \Delta t$. Multiplication of the PDE
system \eqref{eqn.pde} with a test function $\theta_k$ and integration 
over the space-time control volume $\Omega^{st}_i = \Omega_i \times [t^n,t^{n+1}]$ yields 
the following weak form of the governing PDE, which is 
\textit{different} from \eqref{eqn.pde.nc.gw1}, because now the test and 
basis functions are both time dependent:
\begin{equation}
\label{eqn.pde.st1}
\int \limits_{t^n}^{t^{n+1}} \! \! \int \limits_{\Omega_i}
\theta_k(\x,t) \frac{\partial \q_h}{\partial t} \d\x \d t
+\int \limits_{t^n}^{t^{n+1}} \! \! \int \limits_{\Omega_i}
\theta_k(\x,t) \left( \nabla \cdot \F(\Q) + \mathbf{B}(\q_h) \cdot \nabla \q_h  \right) \d\x \d t
= \int \limits_{t^n}^{t^{n+1}} \! \! \int \limits_{\Omega_i}
\theta_k(\x,t) \mathbf{S}(\q_h) \d\x \d t\,.
\end{equation}
Since we are only interested in an element local predictor solution,
i.e. without considering interactions with the neighbor elements 
we do not yet take into account the jumps in $\q_h$ across the element 
interfaces, because this will be done in the final corrector step of the ADER-DG 
scheme \eqref{eqn.pde.nc.gw2}. Instead, we introduce the known discrete 
solution $u_h(\x,t^n)$ at time $t^n$. For this purpose, the first term is
integrated by parts in time. This leads to
\begin{eqnarray}
\label{eqn.pde.st2}
\int \limits_{\Omega_i}   \theta_k(\x,t^{n+1}) \q_h(\x,t^{n+1}) \d\x 
- \int \limits_{t^n}^{t^{n+1}} \! \! \int \limits_{\Omega_i}
\frac{\partial }{\partial t} \theta_k(\x,t) \q_h(\x,t)
\d\x \d t
 - \int \limits_{\Omega_i}   \theta_k(\x,t^{n}) \u_h(\x,t^{n}) \d\x
 = && \nonumber \\
\int \limits_{t^n}^{t^{n+1}} \! \! \int \limits_{\Omega_i^\circ}  \theta_k(\x,t)
 \nabla \cdot \F(\q_h)  \d\x \d t + 
\int \limits_{t^n}^{t^{n+1}} \! \! \int \limits_{\Omega_i^\circ}  \theta_k(\x,t)
\left( \mathbf{S}(\q_h) - \mathbf{B}(\q_h) \cdot \nabla \q_h \right) \d\x \d t. &&
\end{eqnarray}
Using the local space-time ansatz \eqref{eqn.spacetime} Eq.
\eqref{eqn.pde.st2} becomes an element-local nonlinear system for the
unknown degrees of freedom $\hat{\mathbf{q}}_{i,l}$ of the space-time
polynomials $\q_h$. The solution of \eqref{eqn.pde.st2} can be 
found via a simple and fast converging fixed point iteration (a discrete Picard
iteration) as detailed 
e.g. in \cite{Dumbser2008,HidalgoDumbser}. For linear homogeneous 
systems, the discrete Picard iteration converges in a finite number of at most $N+1$ steps, 
since the involved iteration matrix is nilpotent, see \cite{Jackson}. 

However, we emphasize that the choice of an appropriate \textit{initial guess}  
$\q_h^0(\x,t)$ for $\q_h(\x,t)$ is of fundamental importance to obtain a  
faster convergence and thus a computationally more efficient scheme. 
For this purpose, one can either use an extrapolation of $\q_h$ from the previous 
time interval $[t^{n-1},t^n]$, as suggested e.g. in \cite{ADERPrim}, or one can 
employ a second-order accurate MUSCL-Hancock-type approach, as proposed in 
\cite{HidalgoDumbser}, which is based on discrete derivatives computed
at time $t^n$. The initial guess is most conveniently written in terms of a 
Taylor series expansion of the solution in time, where then suitable
approximations of the time derivatives are computed. In the following 
we introduce the operator 
\begin{equation}
\boldsymbol{\mathcal{L}}(\u_h(\x,t^n)) = \S(\u_h(\x,t^n)) - \nabla \cdot \F(\u_h(\x,t^n)) - 
\mathbf{B}\left(\u_h(\x,t^n)\right) \cdot  \nabla \u_h(\x,t^n), 
\end{equation} 
which is an approximation of the time derivative of the solution. 
The second-order accurate MUSCL-type initial guess \cite{HidalgoDumbser} then reads 
\begin{equation} 
\q^0_h(\x,t) = \u_h(\x,t^n) + \left(t - t^n\right)  \boldsymbol{\mathcal{L}}(\u_h(\x,t^n)), 
\end{equation} 
while a third-order accurate initial guess for $\q_h(\x,t)$ is given by  
\begin{equation}
\q^0_h(\x,t) = \u_h(\x,t^n) + \left(t - t^n\right)  \boldsymbol{k}_1 + 
\frac{1}{2} \left(t - t^n\right)^2 \frac{\left( \boldsymbol{k}_2 - \boldsymbol{k}_1 \right)}{\Delta t}. 
\end{equation} 
Here, we have used the abbreviations $\boldsymbol{k}_1 := \boldsymbol{\mathcal{L}}\left(\u_h(\x,t^n)
\right)$ and $\boldsymbol{k}_2 :=
\boldsymbol{\mathcal{L}}\left(\u_h(\x,t^n) + \Delta t \boldsymbol{k}_1
\right)$. For an initial guess of even higher order of accuracy, it is possible to 
use the so--called continuous extension Runge-Kutta (CERK) schemes of Owren and
Zennaro \cite{OwrenZennaro}; see also \cite{Gassner2011a} for the use of CERK
time integrators in the context of high order discontinuous Galerkin finite element
methods. 
If an initial guess with polynomial degree $N-1$ in time is chosen, it is sufficient to use 
\emph{one single} Picard iteration to solve \eqref{eqn.pde.st2} to the desired accuracy.   

At this point, we make some comments about a suitable data-layout for high order ADER-DG schemes. 
In order to compute the discrete derivative operators needed in the predictor \eqref{eqn.pde.st2}, 
especially for the computation of the discrete gradient $\nabla \q_h$, 
it is very convenient to use an array-of-struct (AoS) data structure. In this way, the first or 
fastest-running unit-stride index is the one associated with the $m$ quantities contained in the 
vector $\Q$, while the other indices are associated with the space--time degrees of freedom, 
i.e. we arrange the data contained in the set of degrees of freedom $\hat{\mathbf{q}}^i_l$ as 
$\hat{\mathbf{q}}^i_{v,l_1,l_2,l_3,l_0}$, with $1 \leq v \leq m$ and $1 \leq l_k \leq N+1$. 
The discrete derivatives in any spatial and in time direction can then be simply computed by the 
multiplication of a subset of the degrees of freedom with the transpose of a small 
$(N+1) \times (N+1)$ matrix $D_{kl}$ from the right, which reads
\begin{equation}
  D_{kl} = \frac{1}{h} \left( \int \limits_0^1 \phi_k(\xi) \phi_m(\xi) d \xi \right)^{-1} 
           \left( \int \limits_0^1 \phi_m(\xi) \frac{\partial \phi_l(\xi)}{\partial \xi} d \xi \right),  
\end{equation} 
where $h$ is the respective spatial or temporal step size in the corresponding coordinate direction, i.e.  
either $\Delta x_i$, $\Delta y_i$, $\Delta z_i$ or $\Delta t$.  
For this purpose, the optimized library for small matrix 
multiplications \texttt{libxsmm} can be employed \textcolor{black}{on Intel machines}, 
see \cite{libxsmm} and \cite{SeisSol1,SeisSol2}  for more details.   
However, the AoS data layout is \textit{not} convenient for \textit{vectorization} of the 
PDE evaluation in ADER-DG scheme, since vectorization of the fluxes, source terms and non-conservative products 
should preferably be done over the integration points $l$. For this purpose, we convert the AoS
data layout \textit{on the fly} into a struct-of-array (SoA) data layout via appropriate transposition 
of the data and then call the physical  
flux function $\F(\q_h)$ as well as the combined algebraic source term and non-conservative 
product contained in the expression $\S(\q_h) - \mathbf{B}(\q_h) \cdot \nabla \q_h$ simultaneously
for a subset of \texttt{VECTORLENGTH} space-time degrees of freedom, where \texttt{VECTORLENGTH} is 
the length of the AVX registers of modern Intel Xeon CPUs, i.e. 4 for those with the old 256 bit 
AVX and AVX2 instruction sets (Sandy Bridge, Haswell, Broadwell) and 8 for the latest Intel Xeon 
scalable CPUs with 512 bit AVX instructions (Skylake). The result of the vectorized evaluation 
of the PDE, which is still in SoA format, is then converted back to the AoS data layout using 
appropriate vectorized shuffle commands. 

The element-local space-time predictor is arithmetically very intensive, but at the same time it is also
by construction cache-blocking. 
While in traditional RKDG schemes, each Runge-Kutta stage requires touching
all spatial degrees of freedom of the entire domain once per Runge-Kutta stage, in our ADER-DG approach 
the spatial degrees of freedom $\u_h$ need to be loaded only once per element and time step, and from those 
all space-time degrees of freedom of $\q_h$ are computed. Ideally, this procedure fits entirely into the
L3 cache or even into the L2 cache of the CPU, at least up to a certain critical polynomial degree 
$N_c=N_c(m)$, which is a function 
of the available L3 or L2 cache size, but also of the number of quantities $m$ to be evolved in the PDE system. 

Last but not least, it is important to note that it is possible to hide the entire MPI communication that 
is inevitably needed on distributed memory supercomputers behind the space-time predictor. For this purpose, the predictor
is first invoked on the MPI boundaries of each CPU, which then immediately sends the boundary-extrapolated data
$\q_h^-$ and $\q_h^+$ to the neighbor CPUs. While the messages containing the data of these non-blocking MPI send and receive commands are sent around, each CPU can compute the space-time predictor of purely interior elements that do not need any MPI communication.   

For an efficient task-based formalism used within ExaHyPE in the context of shared memory
parallelism, see \cite{StopTalking}. 
This completes the description of the efficient implementation of the unlimited ADER-DG schemes
used within the ExaHyPE engine. 

\subsection{A posteriori subcell finite volume limiter}

In regions where the discrete solution is smooth, there is indeed no need for using
nonlinear limiters. However, in the presence of shock waves, discontinuities or 
strong gradients, and taking into account the fact that even a \textit{smooth signal} 
may become \textit{non-smooth} on the discrete level if 
it is \textit{underresolved} on the grid, we have to supplement our high order 
unlimited ADER-DG scheme described above with a nonlinear limiter. 

In order to build a simple, robust and accurate limiter, we follow the  
ideas outlined in \cite{DGLimiter1,DGLimiter2,DGLimiter3,ALEDG}, where a novel 
\textit{a posteriori} limiting strategy for ADER-DG schemes was developed, 
based on the ideas of the MOOD paradigm introduced in \cite{CDL1,CDL2,CDL3,ADER_MOOD_14}
in the finite volume context.   
In a first run, the unlimited ADER-DG scheme is used and produces a so-called
\textit{candidate solution}, denoted by $\u_h^*(\x,t^{n+1})$ in the following. 
This candidate solution is then checked \textit{a posteriori} against several 
physical and numerical detection criteria. For example, we require some relevant
physical quantities of the solution 
to be positive (e.g. pressure and density), we require the absence of floating point 
errors (\texttt{NaN}) and we impose a relaxed discrete maximum principle (DMP) in 
the sense of polynomials, see \cite{DGLimiter1}.  As soon as one of these
detection criteria is not satisfied, a cell is marked as troubled zone and
is scheduled for limiting. 
  
A cell  
$\Omega_i$ that has been marked for limiting is now split into $(2N+1)^d$
finite volume subcells, which are denoted by $\Omega_{i,s}$. They 
satisfy $\Omega_i = \bigcup_s \Omega_{i,s}$. Note that this very fine 
division of a DG element into finite volume subcells does \textit{not}
reduce the time step of the overall ADER-DG scheme, since the CFL number
of explicit DG schemes scales with $1/(2N+1)$, while the CFL number of
finite volume schemes (used on the subgrid) is of the order of unity.
The discrete solution in the subcells $\Omega_{i,s}$ is represented at
time $t^n$ in terms of \textit{piecewise constant} subcell averages
$\bar{\u}^n_{i,s}$, i.e.
\begin{equation}
\label{eqn.subcellaverage}
\bar{\u}^n_{i,s} = \frac{1}{|\Omega_{i,s}|} \int \limits_{\Omega_{i,s}}
\Q(\x,t^n) \d\x\,.
\end{equation}
These subcell averages are now evolved in time with a second or third order
accurate finite volume scheme, which actually looks very similar to the 
previous ADER-DG scheme \eqref{eqn.pde.nc.gw2}, with the difference that now
the test function is unity and the spatial control volumes $\Omega_i$ are replaced
by the sub-volumes $\Omega_{i,s}$: 
\begin{equation}
\label{eqn.pde.nc.fv}
\left( \bar{\mathbf{u}}^{n+1}_{i,s} - \bar{\mathbf{u}}^{n}_{i,s}  \right)
+ \int \limits_{t^n}^{t^{n+1}} \! \! \int \limits_{\partial \Omega_{i,s}}
\mathcal{D}^-\left( \q_h^-, \q_h^+ \right) \cdot \mathbf{n} \, \d S \d t 
+ \int \limits_{t^n}^{t^{n+1}} \! \! \int \limits_{\Omega_{i,s}^\circ}
\left( \mathbf{B}(\q_h) \cdot \nabla \q_h  \right)  \d\x \d t 
= \int \limits_{t^n}^{t^{n+1}} \! \! \int \limits_{\Omega_{i,s}}
\mathbf{S}(\q_h)  \d\x \d t\,.
\end{equation}
Here we use again a space-time predictor solution $\q_h$, but which is
now computed from an initial condition given by a second order TVD 
reconstruction polynomial, or from a WENO \cite{shu_efficient_weno}
or CWENO reconstruction \cite{LPR:99,LPR:2001,SCR:CWENOquadtree,ADER_CWENO} 
polynomial $\w_h(\x,t^n)$ computed from the cell averages 
$\bar{\u}^n_{i,s}$ via an appropriate reconstruction operator. The predictor
is either computed via a standard second order MUSCL-Hancock-type 
strategy, or via the space-time DG approach of \eqref{eqn.pde.st2}, but where 
the initial data $\u_h(\x,t^n)$ are now replaced by $\w_h(\x,t^n)$ and the 
spatial control volumes $\Omega_i$ are replaced by the subcells $\Omega_{i,s}$. 

Once all subcell averages $\bar{\mathbf{u}}^{n+1}_{i,s}$ inside a cell
$\Omega_i$ have been computed according to \eqref{eqn.pde.nc.fv}, the
limited DG polynomial $\u'_h(\x,t^{n+1})$ at the next time level is
obtained again via a classical constrained least squares reconstruction
procedure requiring
\begin{equation}
\frac{1}{|\Omega_{i,s}|} \int \limits_{\Omega_{i,s}} \u'_h(\x,t^{n+1}) \d\x
= \bar{\u}^{n+1}_{i,s} \qquad \forall \Omega_{i,s} \in \Omega_i,
\qquad \textnormal{ and } \qquad \int \limits_{\Omega_{i}} \u'_h(\x,t^{n+1})
\d\x = \sum \limits_{\Omega_{i,s} \in \Omega_i} |\Omega_{i,s}| \bar{\u}^{n+1}
{i,s}\,.
\end{equation}
Here, the second relation is a constraint and means conservation at the
level of the control volume $\Omega_i$. This completes the brief
description of the subcell finite volume limiter used here.

\section{Some examples of typical PDE systems solved with the ExaHyPE engine}

The great advantage of ExaHyPE over other existing PDE solvers is its great flexibility and versatility
for the solution of a very wide class of hyperbolic PDE systems \eqref{eqn.pde}.  
The implementation of the numerical method and the definition of the PDE system to be solved are
completely independent of each other. The compute kernels are provided either as generic   or as an optimized implementation 
for the general PDE system given by \eqref{eqn.pde}, while the user only needs to provide  
particular implementations of the functions $\F(\Q)$, $\mathbf{B}(\Q)$ and $\S(\Q)$. It is 
obviously also possible to drop terms that are not needed. This allows to solve all the 
PDE systems listed below in one single software package. \textcolor{black}{ In all numerical examples shown below, 
we have used a CFL condition of the type 
\begin{equation} 
 \Delta t \leq \frac{\alpha}{ \frac{|\lambda^x_{\max}|}{\Delta x} + 
 	                          \frac{|\lambda^y_{\max}|}{\Delta y} + 
 	                          \frac{|\lambda^z_{\max}|}{\Delta z} },  
\end{equation}
where $\Delta x$, $\Delta y$ and $\Delta z$ are the mesh spacings and $|\lambda^x_{\max}|$, $|\lambda^x_{\max}|$ 
and $|\lambda^x_{\max}|$ are the maximal absolute values of the eigenvalues (wave speeds) of the matrix 
$\mathbf{A} \cdot \mathbf{n}$ in $x$, $y$ and $z$  
direction, respectively. The coefficient $\alpha < 1/(2N+1)$ can be obtained via a numerical von Neumann stability 
analysis and is reported for some relevant $N$ in \cite{Dumbser2008}. 
} 

\subsection{The Euler equations of compressible gas dynamics}
The Euler equations of compressible gas dynamics are among the simplest nonlinear systems of 
hyperbolic conservation laws. They only involve a conservative flux $\F(\Q)$ and read 
\begin{equation}
\label{eqn.euler} 
\frac{\partial}{\partial t}  \left( \begin{array}{c} \rho \\ \rho \mathbf{v} \\ \rho E \end{array} \right) + \nabla \cdot \left( \begin{array}{c} \rho \mathbf{v} \\ \rho \mathbf{v} \otimes \mathbf{v} + p \mathbf{I} \\ \mathbf{v} \left(\rho E + p \right) \end{array} \right) = 0. 
\end{equation}
Here, $\rho$ is the mass density, $\mathbf{v}$ is the fluid velocity, $\rho E$ is the total energy density and $p$ is the fluid pressure, which is related to $\rho$, $\rho E$ and $\mathbf{v}$ via the so--called equation of state (EOS).   
In the following we show the computational results for two test problems. The first one is the smooth isentropic vortex test case first proposed in \cite{HuShuTri} and also used in \cite{DGLimiter1}, which has an exact solution and is therefore suitable for a numerical convergence study. Some results of \cite{DGLimiter1} are summarized in Table \ref{tab.vortex}
below, where $N_x$ denotes the number of cells per space dimensions. From the results one can conclude that the high order ADER-DG schemes converge with the designed order of accuracy  
in both space and time. In order to give a quantitative assessment for the cost of the scheme, we define and provide the TDU metric, which is the cost per degree of freedom update per CPU 
core, see also \cite{Dumbser2008}.  
The TDU metric is easily computed by dividing the measured wall clock time (WCT) of a simulation by the number of elements per CPU core and time steps carried out, and by the number of spatial degrees of freedom per element, i.e. $(N+1)^d$. 
With the appropriate initial guess and AVX 512 vectorization of the code discussed in the previous section, the cost for updating one single degree of freedom for a fourth order 
ADER-DG scheme ($N=3$) for the 3D compressible Euler equations is as low as  
TDU=0.25 $\mu$s when using one single CPU core of a new Intel i9-7900X Skylake test workstation 
with 3.3 GHz nominal clock speed, 32 GB of RAM and a total number of 10 CPU cores. This cost  
metric can be directly compared with the cost to update one single point or control volume 
of existing finite difference and finite volume schemes.    

\begin{figure}
	\centering
	\begin{tabular}{cc}
		\includegraphics[width=0.45\textwidth]{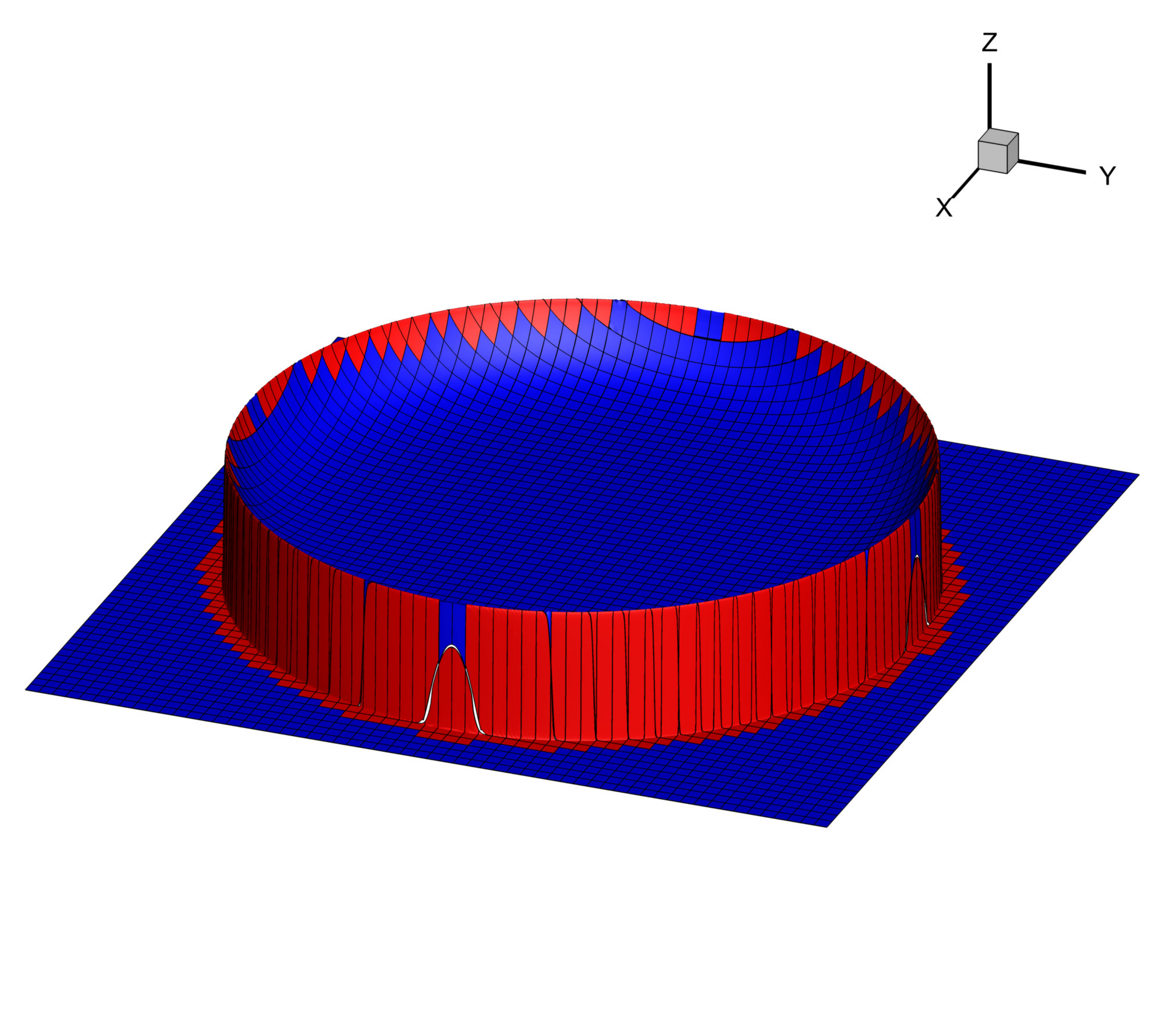} &  
		\includegraphics[width=0.45\textwidth]{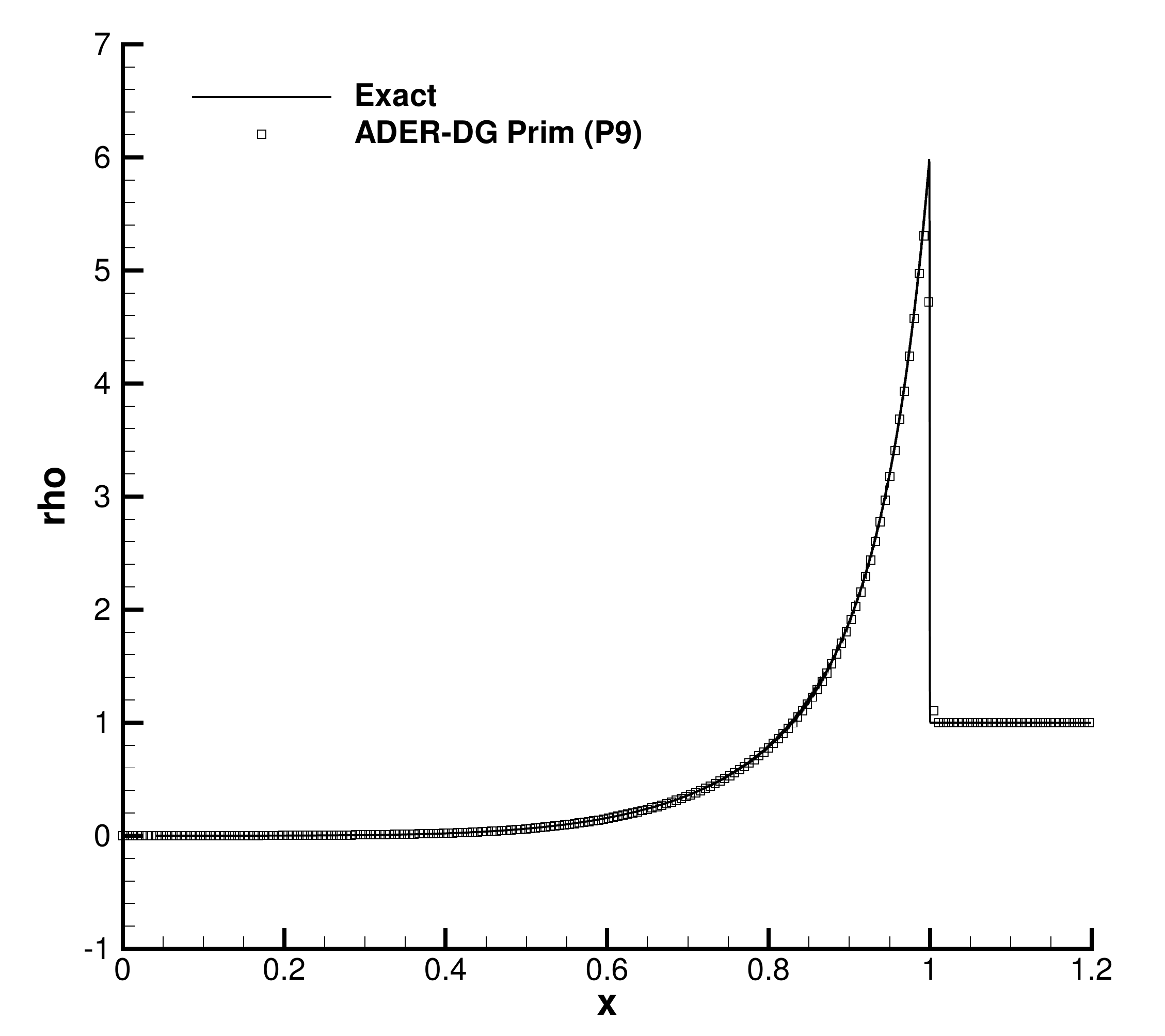}
	\end{tabular} 
	\caption{Sedov blast wave problem using an ADER-DG $P9$ scheme with \textit{a posteriori}
	subcell finite volume limiter using predictor and limiter in primitive variables, see \cite{ADERPrim}. Unlimited cells are depicted in blue, while limited cells are highlighted in red (left). 1D cut through the numerical solution  and comparison with the exact solution (right). }
	\label{fig.sedov}
\end{figure}

 \begin{table}[!htbp] 
	\centering
	\begin{tabular}{ccccccccc}
		\hline
		& $N_x$ & $L^1$ error & $L^2$ error & $L^\infty$ error & $L^1$ order & $L^2$ order & $L^\infty$ order &
		Theor. \\
		\hline
		\multirow{4}{*}{{{$N=3$}}}
		& 25	& 5.77E-04	& 9.42E-05	& 7.84E-05	& ---	& ---	& --- & \multirow{4}{*}{4}\\
		& 50	& 2.75E-05	& 4.52E-06	& 4.09E-06	& 4.39	& 4.38	& 4.26 &\\
		& 75	& 4.36E-06	& 7.89E-07	& 7.55E-07	& 4.55	& 4.30	& 4.17 &\\
		& 100	& 1.21E-06	& 2.37E-07	& 2.38E-07	& 4.46	& 4.17	& 4.01 &\\
		\cline{2-8}
		\hline
		\multirow{4}{*}{{{$N=4$}}}
		& 20	& 1.54E-04	& 2.18E-05	& 2.20E-05	& ---	& ---	& --- & \multirow{4}{*}{5}\\		
		& 30	& 1.79E-05	& 2.46E-06	& 2.13E-06	& 5.32	& 5.37	& 5.75 &\\
		& 40	& 3.79E-06	& 5.35E-07	& 5.18E-07	& 5.39	& 5.31	& 4.92 &\\
		& 50	& 1.11E-06	& 1.61E-07	& 1.46E-07	& 5.50	& 5.39	& 5.69 &\\
		\cline{2-8}
		\hline
		\multirow{4}{*}{{{$N=5$}}}
		& 10	& 9.72E-04	& 1.59E-04	& 2.00E-04	& ---	& ---	& --- & \multirow{4}{*}{6}\\		
		& 20	& 1.56E-05	& 2.13E-06	& 2.14E-06	& 5.96	& 6.22	& 6.55 &\\
		& 30	& 1.14E-06	& 1.64E-07	& 1.91E-07	& 6.45	& 6.33	& 5.96 &\\
		& 40	& 2.17E-07	& 2.97E-08	& 3.59E-08	& 5.77	& 5.93	& 5.82 &\\
		\cline{2-8}
		\hline
	\end{tabular}
	\caption{ \label{tab.vortex} $L^1, L^2$ and $L^\infty$ errors and numerical convergence rates 
		obtained for the two--dimensional isentropic vortex test problem using different unlimited ADER-DG schemes, see \cite{DGLimiter1}. }
\end{table}
In the following we show the results obtained with an ADER-DG scheme using piecewise polynomials
of degree $N=9$ for a very stringent test case, which is the so-called Sedov blast wave problem
detailed in \cite{Sedov,SedovExact,ALEDG,ADERPrim}. It consists in an explosion propagating in a zero pressure gas, leading to an infinitely strong shock wave. In our setup, the outer pressure is set to $10^{-14}$, i.e. close to machine zero. In order to get a robust numerical
scheme, it is useful to perform the reconstruction step in the subcell finite volume limiter as well as the space-time predictor of the ADER-DG scheme in primitive variables, see \cite{ADERPrim}. The computational results obtained are shown in Fig. \ref{fig.sedov}, where we can observe a very good agreement with the reference solution. One furthermore can see that the discrete solution respects the circular
symmetry of the problem and the \textit{a posteriori} subcell limiter is only acting in 
the vicinity of the shock wave.  
 
\subsection{A novel diffuse interface approach for linear seismic wave propagation in complex geometries}

\begin{figure}
	\centering
	\begin{tabular}{cc}
		\includegraphics[width=0.45\textwidth]{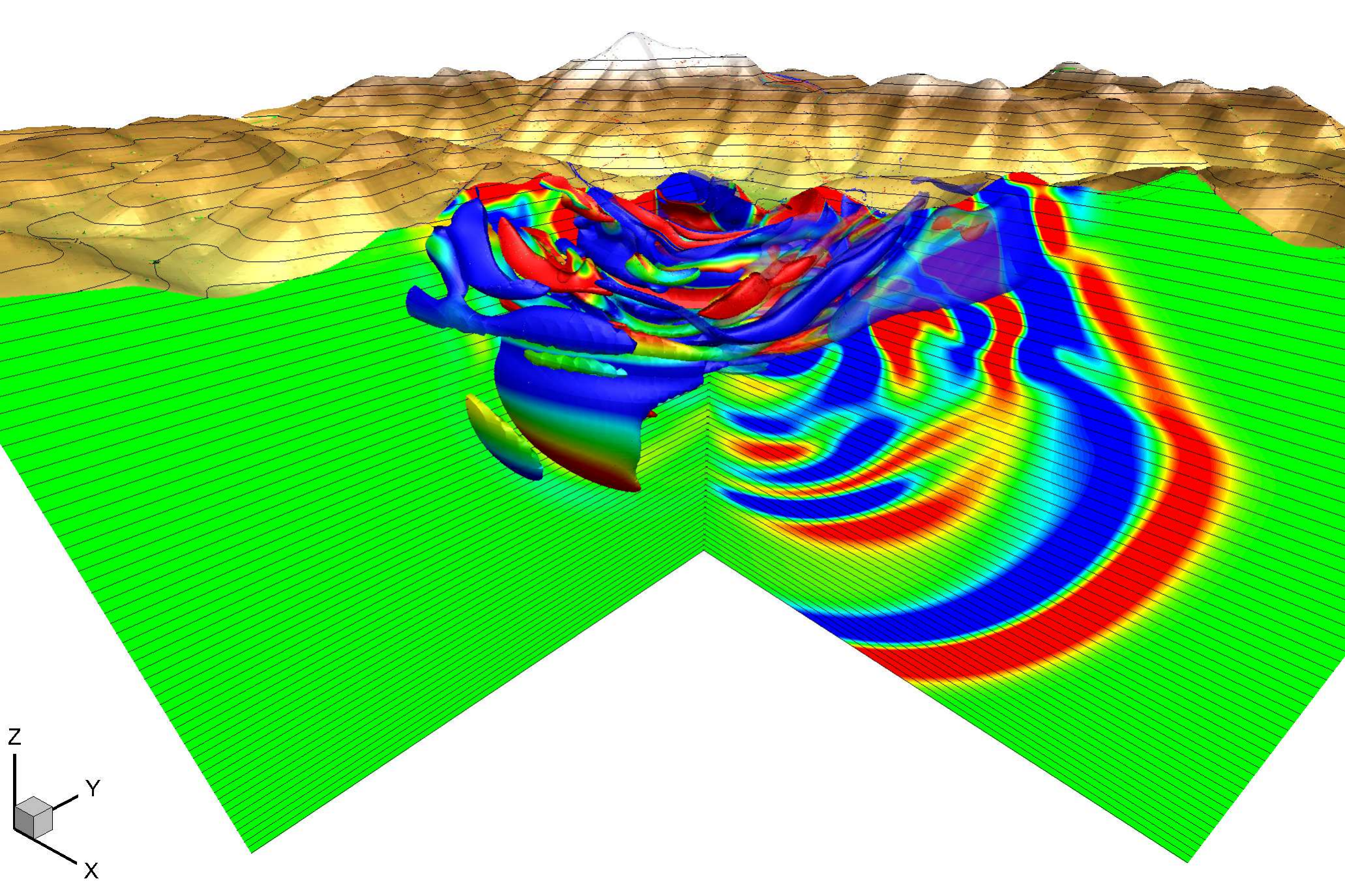} &  
		\includegraphics[width=0.45\textwidth]{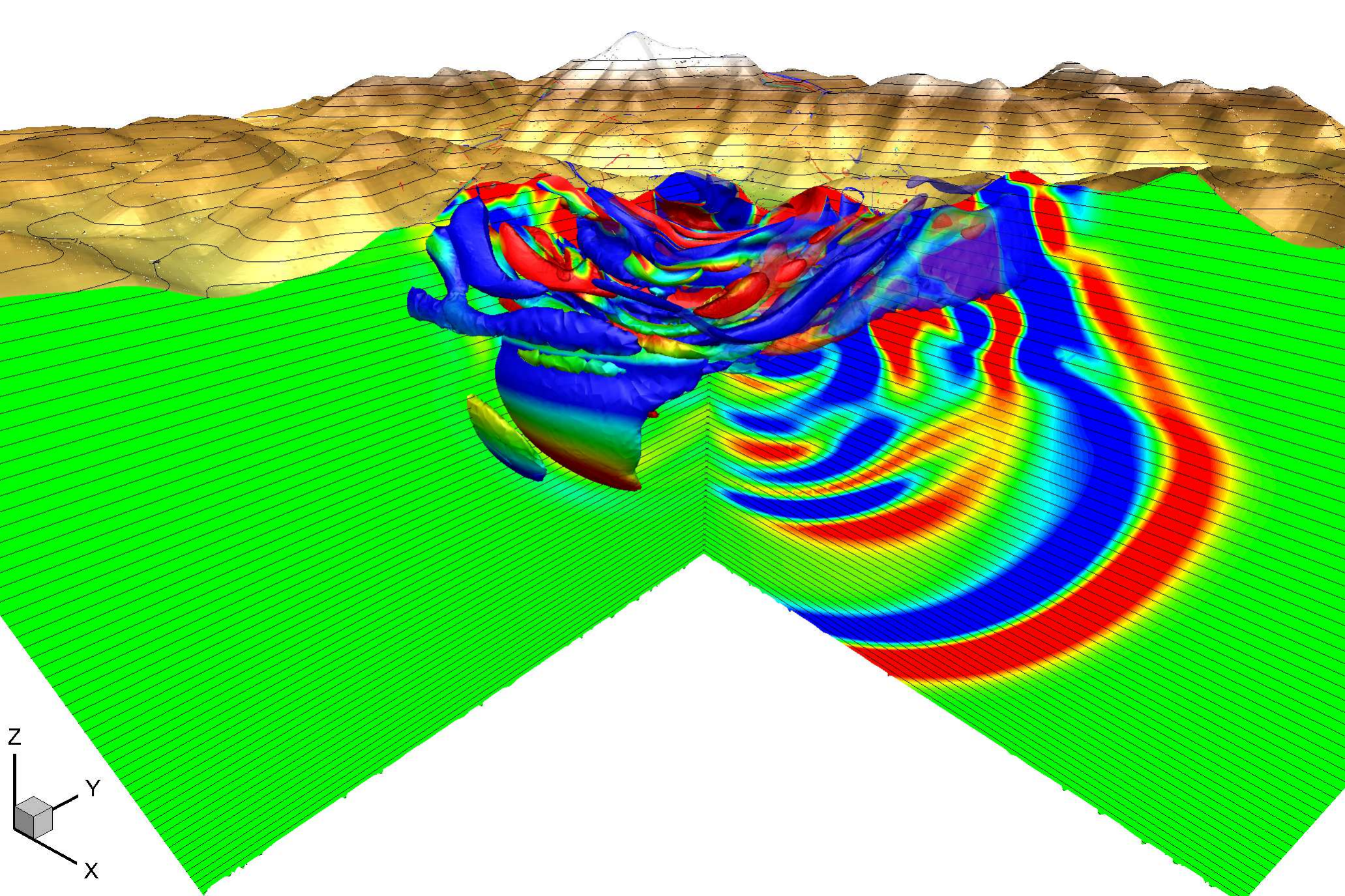}
	\end{tabular} 
	\caption{Wave field of a seismic wave propagation problem with the novel diffuse interface
		approach on adaptive Cartesian grids developed in \cite{AMRDIM} (left) compared with the reference solution obtained on a classical boundary-fitted unstructured tetrahedral mesh \cite{gij2} (right).}
	\label{fig.dim}
\end{figure}

\begin{figure}
	\centering
		\includegraphics[width=0.95\textwidth]{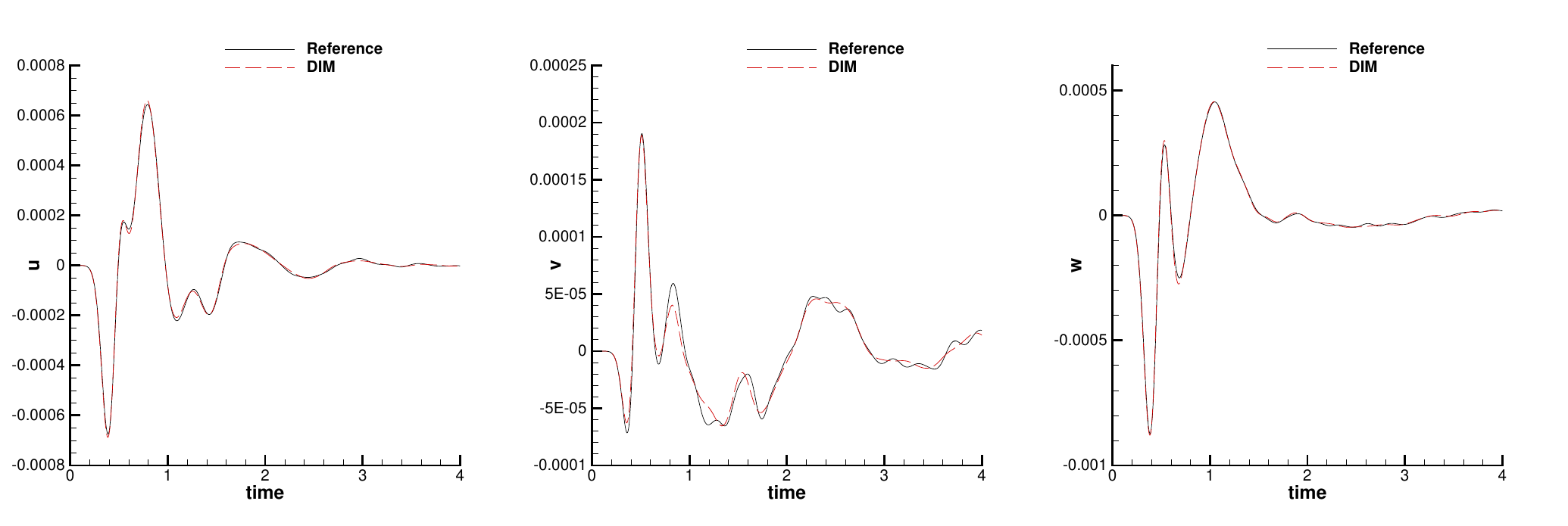} \\   
		\includegraphics[width=0.95\textwidth]{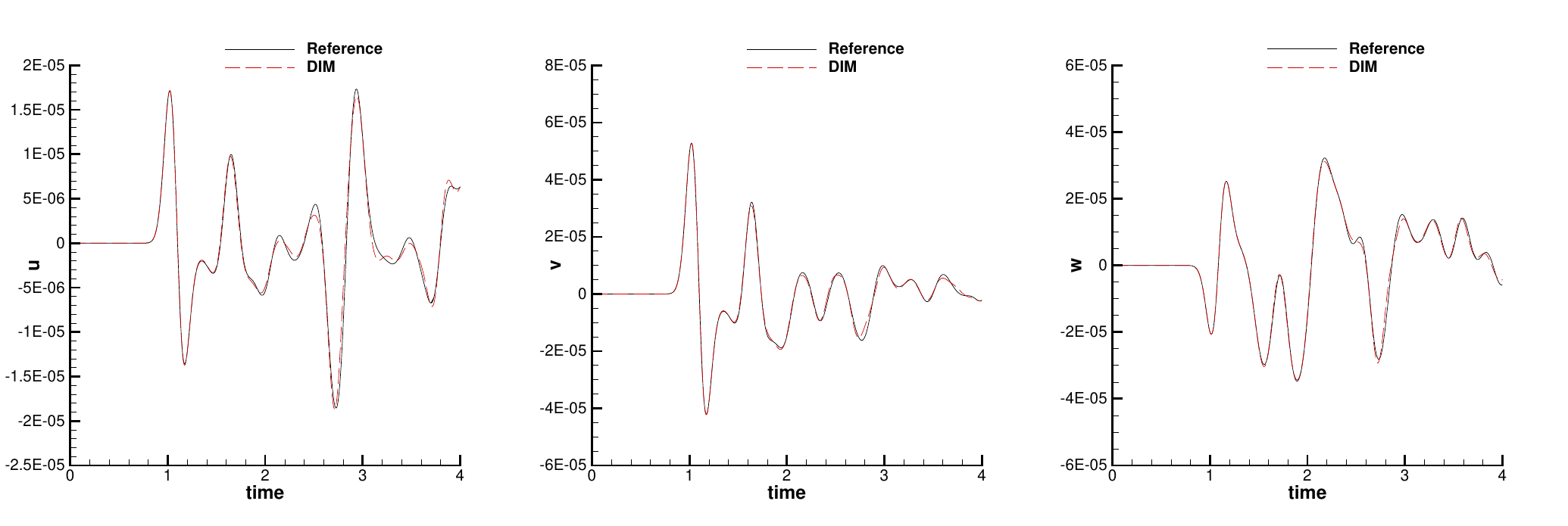}
	\caption{Seismogram recordings in two observation points obtained with the diffuse interface approach on adaptive Cartesian meshes \cite{AMRDIM} and with a reference solution 
	obtained with high order ADER-DG schemes on boundary-fitted unstructured meshes \cite{gij2}. }
	\label{fig.dim.seis}
\end{figure}

Seismic wave propagation problems in complex 3D geometries are often very challenging due to 
the geometric complexity. Standard approaches either use regular curvilinear boundary-fitted meshes, or unstructured tetrahedral or hexahedral meshes. In all cases, a certain amount of user interaction for grid generation is required. Furthermore, the geometric complexity can 
have a negative impact on the admissible time step size due to the CFL condition, since the
mesh generator may create elements with very bad aspect ratio, so--called sliver elements. 
In the case of regular curvilinear grids, the Jacobian of the mapping may become 
ill-conditioned and thus reduce the admissible time step size. In \cite{AMRDIM} a novel
diffuse interface approach has been forwarded, where only the definition of a scalar volume
fraction function $\alpha$ is required, where $\alpha=1$ is set inside the solid medium, 
and $\alpha=0$ in the surrounding gas or vacuum. The governing PDE system proposed in \cite{AMRDIM} reads
\begin{eqnarray}
\frac{\partial \boldsymbol{\sigma}}{\partial t } - \E(\lambda, \mu) \cdot \frac{1}{\alpha} \nabla {(\alpha \mathbf{v})} + \frac{1}{\alpha}\E(\lambda, \mu) \cdot {\mathbf{v}} \otimes \nabla \alpha = \mathbf{S}_\sigma, \label{eq:3.00}\\
\diff{\alpha \mathbf{v}}{t}-\frac{\alpha}{\rho}\nabla \cdot \boldsymbol{\sigma} - \frac{1}{\rho} \boldsymbol{\sigma} \nabla \alpha=\mathbf{S}_v, \\ 
\label{eq:3.0}
\diff{\alpha}{t} = 0,  \qquad 
\diff{\lambda}{t}=0,   \qquad 	\diff{\mu}{t}=0,  \qquad \diff{\rho}{t}=0, \label{eq:3.1}
\end{eqnarray}
and clearly falls into the class of PDE systems described by \eqref{eqn.pde}. Here, $\boldsymbol{\sigma}$ denotes the symmetric stress tensor, $\mathbf{v}$ is the velocity
vector, $\alpha \in [0,1]$ is the volume fraction, $\lambda$ and $\mu$ are the Lam\'e constants
and $\rho$ is the density of the solid medium. The elasticity tensor $\E$ is a function of $\lambda$ and $\mu$ and relates stress and strain via the Hooke law. 
The last four quantities obey trivial 
evolution equations, which state that these parameters remain constant in time. \textcolor{black}{However, they still need to be properly included in the evolution system, since they have an influence on the solution of the Riemann problem.}
An analysis 
of the eigenstructure of \eqref{eq:3.00} - \eqref{eq:3.1} shows that the eigenvalues are
all real and are \textit{independent} of the volume fraction function $\alpha$. Furthermore,
the exact solution of a generic Riemann problem with $\alpha=1$ on the left and $\alpha=0$ on
the right yields the free surface boundary condition $\boldsymbol{\sigma} \cdot \mathbf{n}=0$
at the interface, see \cite{AMRDIM} for details. In this new approach, the mesh generation problem can be fully avoided, since all that is needed is the specification of the scalar
volume fraction function $\alpha$, which is set to unity inside the solid and to zero outside. 
A realistic 3D wave propagation example based on real DTM data of the Mont Blanc region is shown in Figures \ref{fig.dim} and \ref{fig.dim.seis}, where the 3D contour colors of the 
wave field  as well as a set of seismogram recordings in two receiver points are reported. 
For this simulation, a uniform Cartesian base-grid of $80^3$ elements was used, together with
one level of AMR refinement close to the free surface boundary determined by the DTM model. 
A fourth order ADER-DG scheme ($N=3$) has been used in this simulation. We stress that the
entire setup of the computational model in the diffuse interface approach is completely 
automatic, and no manual user interaction was required. 
The reference solution was obtained with a high order ADER-DG scheme \textcolor{black}{of the same polynomial 
degree $N=3$} using an unstructured boundary-fitted tetrahedral mesh \cite{gij2} of \textcolor{black}{similar 
spatial resolution, containing a total of $1,267,717$ elements}. We observe an excellent agreement between the  
two simulations, which were obtained with two completely different PDE systems on two 
different grid topologies.  

\subsection{The unified Godunov-Peshkov-Romenski model of continuum mechanics (GPR)} 

A major achievement of ExaHyPE was the first successful numerical solution of the 
unified first order symmetric hyperbolic and thermodynamically compatible 
Godunov--Peshkov--Romenski (GPR) model of continuum mechanics, see \cite{HPRmodel,HPRmodelMHD}. 
The GPR model is based on the seminal papers by Godunov and Romenski \cite{Godunov:1961a, GodunovRomenski72, Rom1998} 
on inviscid symmetric hyperbolic systems. The dissipative mechanisms, which allow to model
both plastic solids as well as viscous fluids within one single set of equations were added 
later in the groundbreaking work of Peshkov and Romenski in \cite{PeshRom2014}. The GPR model
is briefly outlined below, while for all details the interested reader is referred to 
\cite{PeshRom2014,HPRmodel,HPRmodelMHD}. The governing equations read     
\begin{equation}
\label{eqn.conti} 
\frac{\partial \rho}{\partial t} + \frac{\partial}{\partial x_k} \left(  \rho u_k \right)  = 0, 
\end{equation} 

\begin{equation}
\label{eqn.momentum} 
\frac{\partial \rho u_i}{\partial t} + \frac{\partial}{\partial x_k} \left( \partial \rho u_i u_k + p \delta_{ik} - \sigma_{ik} \right) = 0, 
\end{equation} 

\begin{equation}
\label{eqn.deformation} 
\frac{\partial A_{ik}}{\partial t} + \frac{\partial \left(A_{im} u_m \right)}{\partial x_k}  + u_j \left( \frac{\partial A_{ik}}{\partial x_j} - \frac{\partial A_{ij}}{\partial x_k}  \right)  = 
- \frac{\psi_{ik}}{\theta_1(\tau_1)}, 
\end{equation} 

\begin{equation}
\label{eqn.heatflux} 
\frac{\partial \rho J_i}{\partial t} + \frac{\partial }{\partial x_k} \left( \rho J_i u_k + T \delta_{ik}  \right) = -
\frac{1}{\theta_2(\tau_2)} \rho H_i, 
\end{equation} 

\begin{equation}
\label{eqn.energy} 
\frac{\partial \rho E}{\partial t} + \frac{\partial  }{\partial x_k} \left( u_k \rho E + u_i \left( p \delta_{ik} - \sigma_{ik} \right) + q_k \right) = 0. 
\end{equation} 
Furthermore, the system is also endowed with an entropy inequality, see \cite{HPRmodel}. 
Here, $\rho$ is the mass density, $[u_i] = \mathbf{v} = (u,v,w)$ is the velocity vector, $p$ is a non-equilibrium pressure, $[A_{ik}] = \mathbf{A}$ is the distorsion field, 
$[J_i]=\mathbf{J}$ is the thermal impulse vector, $T$ is the temperature and $\rho E$ is the total energy density that is defined according to \cite{HPRmodel} as 
\begin{equation}
\label{eqn.totenergy} 
\rho E = \rho e + \halb \rho \mathbf{v}^2 +  \frac{1}{4} \rho c_s^2 \, \textnormal{tr} \left( (\textnormal{dev} \G )^T (\textnormal{dev} \G) \right) + \halb \rho \alpha^2  \mathbf{J}^2  
\end{equation} 
in terms of the specific internal energy $e=e(p,\rho)$ given by the usual equation of state (EOS), 
the kinetic energy, the energy stored in the medium due to deformations and in the thermal impulse.  
Furthermore, $\G = \mathbf{A}^T \mathbf{A}$ is a metric tensor induced by the distortion field $\mathbf{A}$, which allows to measure distances and thus deformations in the medium, $c_s$ is 
the shear sound speed and $\alpha$ is a heat wave propagation speed; the symbol 
$\textnormal{dev} \G = \G - \frac{1}{3} \textnormal{tr} \, \G$ indicates the trace-free part of the metric tensor $\G$.   From the definition 
of the total energy 
\eqref{eqn.totenergy} and the relations $H_i = E_{J_i}$, $\psi_{ik} = E_{A_{ij}}$, 
$\sigma_{ik} = -\rho A_{mi}  E_{A_{mk}}$, $T=E_S$ and $q_k = E_S E_{J_k}$ the shear
stress tensor and the heat flux read $\boldsymbol{\sigma} = -\rho c_s^2 \G \dev {\G}$ 
and $\mathbf{q} = \alpha^2 T \mathbf{J}$. It can furthermore be shown via formal asymptotic 
expansion \cite{HPRmodel} that via an appropriate choice of $\theta_1$ and $\theta_2$ 
in the stiff relaxation limit $\tau_1 \to 0$ and $\tau_2 \to 0$, the stress 
tensor and the heat flux tend to those of the compressible Navier-Stokes equations 
\begin{equation}
\boldsymbol{\sigma} \to \mu \left( \nabla \mathbf{v} + \nabla \mathbf{v}^T - \frac{2}{3} 
\left( \nabla \cdot \mathbf{v} \right)  
\mathbf{I} \right) \qquad \textnormal{and} \qquad 
\mathbf{q} \to - \lambda \nabla T,   
\label{eqn.aslim}
\end{equation}   
with transport coefficients $\mu=\mu(\tau_1,c_s)$ and $\lambda = \lambda(\tau_2,\alpha)$ 
related to the relaxation times $\tau_1$ and $\tau_2$ and
to the propagation speeds $c_s$ and $\alpha$, respectively. For a complete derivation, see 
\cite{HPRmodel,HPRmodelMHD}. In the opposite limit $\tau_1 \to \infty$ the model describes
an ideal elastic solid with large deformations. This means that elastic solids as well
as viscous fluids can be described at the aid of the same mathematical model. 
\textcolor{black}{At this point we stress that numerically we always solve the unified \textit{first order
hyperbolic} PDE system \eqref{eqn.conti}-\eqref{eqn.energy}, even in the stiff relaxation limit \eqref{eqn.aslim}, 
when the compressible Navier-Stokes-Fourier system is retrieved asymptotically. We emphasize that we never 
need to discretize any parabolic terms, since the hyperbolic system \eqref{eqn.conti}-\eqref{eqn.energy} 
with algebraic relaxation source terms fits perfectly into the framework of Eqn. \eqref{eqn.pde}. }

\begin{figure}[!t]
	\centering
	\begin{tabular}{cc} 
	\includegraphics[width=0.45\textwidth]{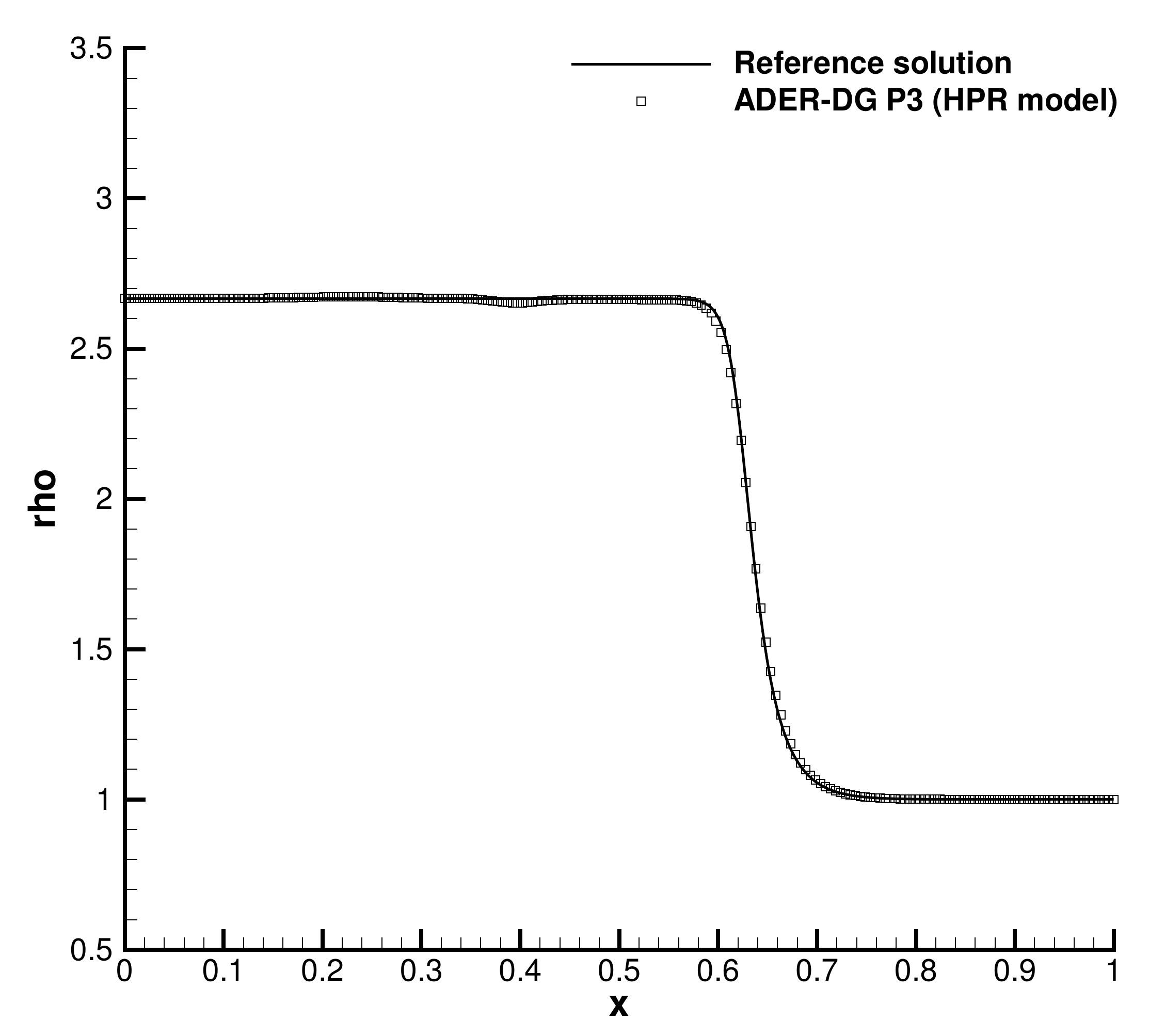} &    
	\includegraphics[width=0.45\textwidth]{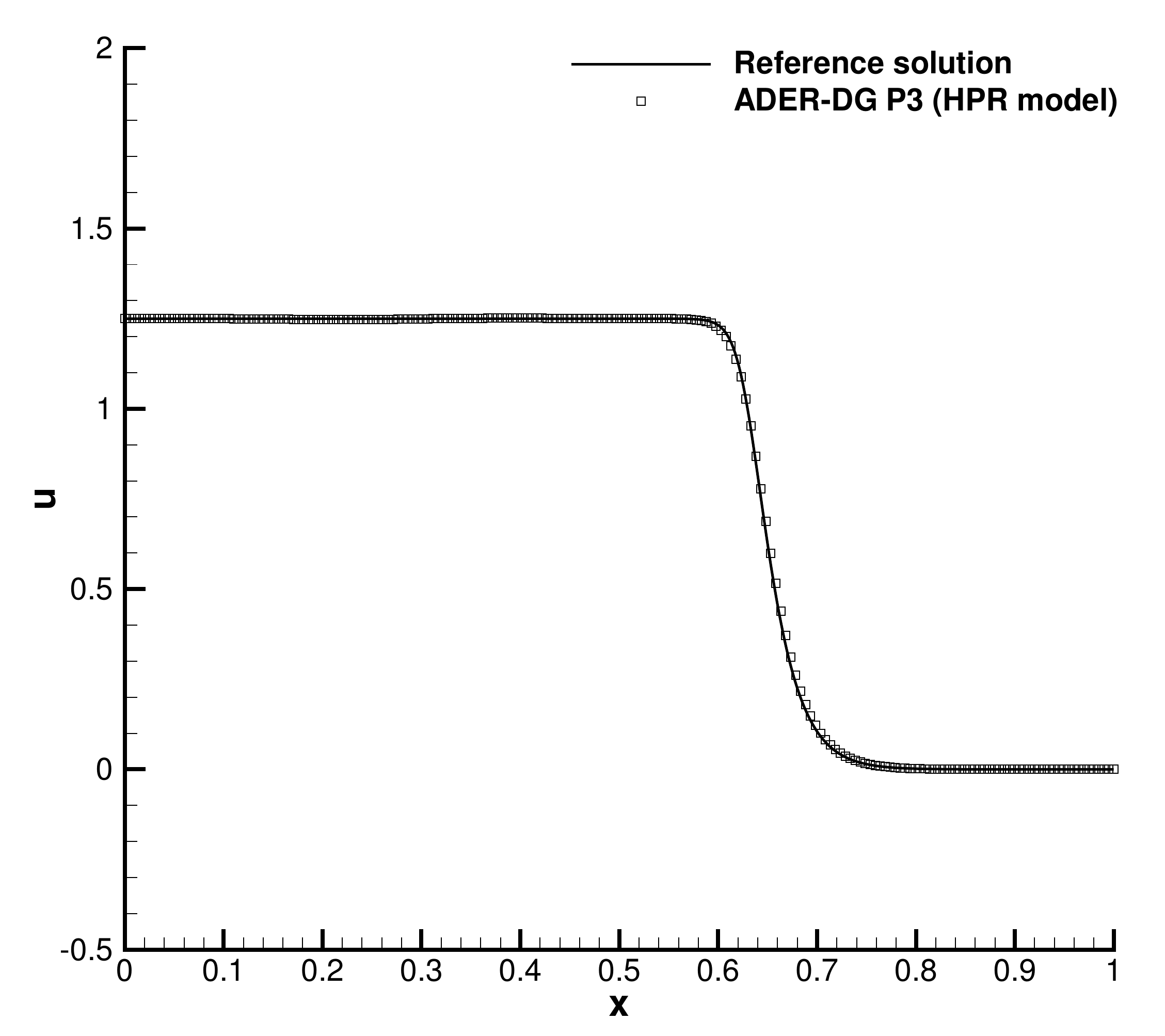} \\ 
	\includegraphics[width=0.45\textwidth]{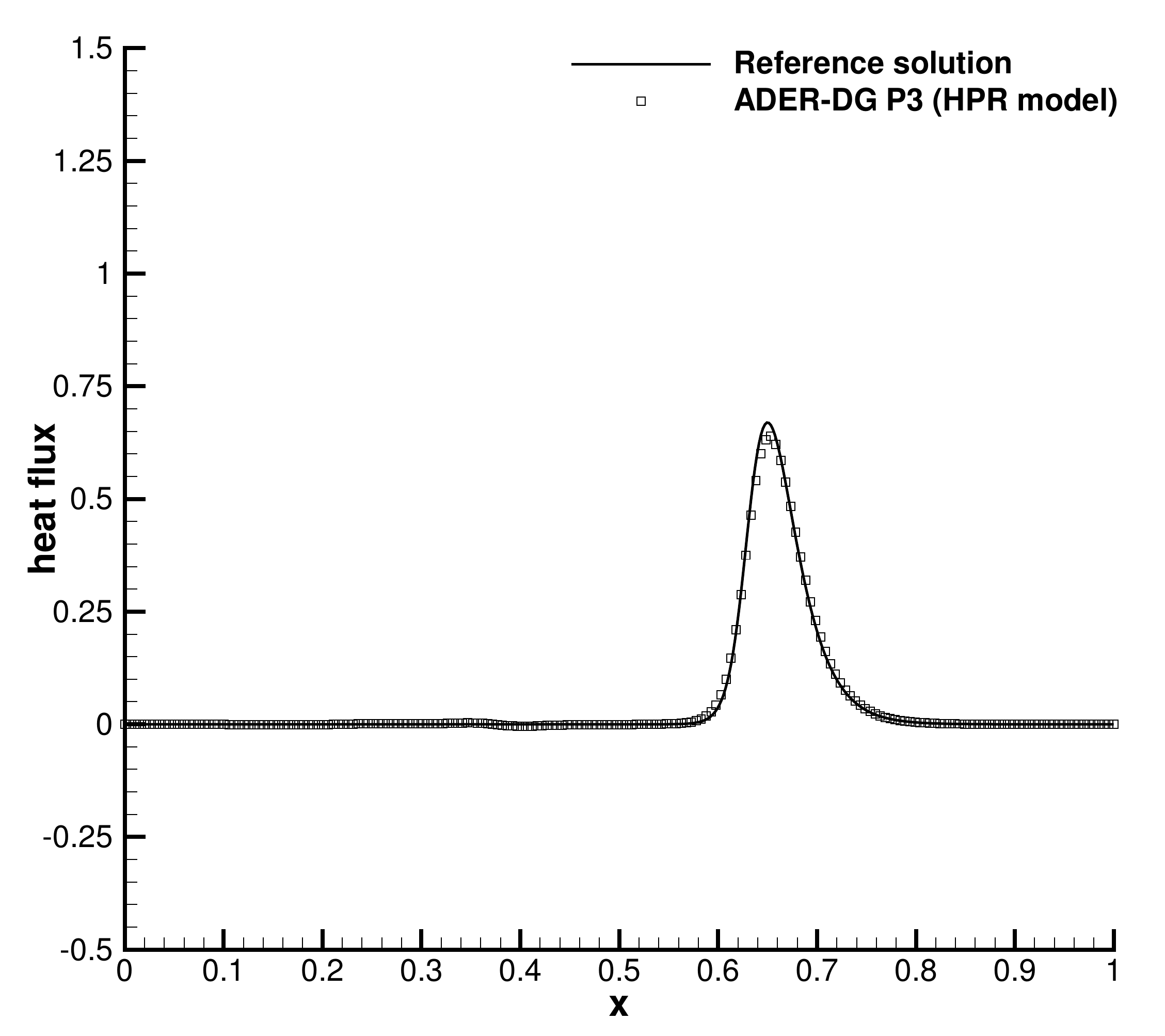} &    
\includegraphics[width=0.45\textwidth]{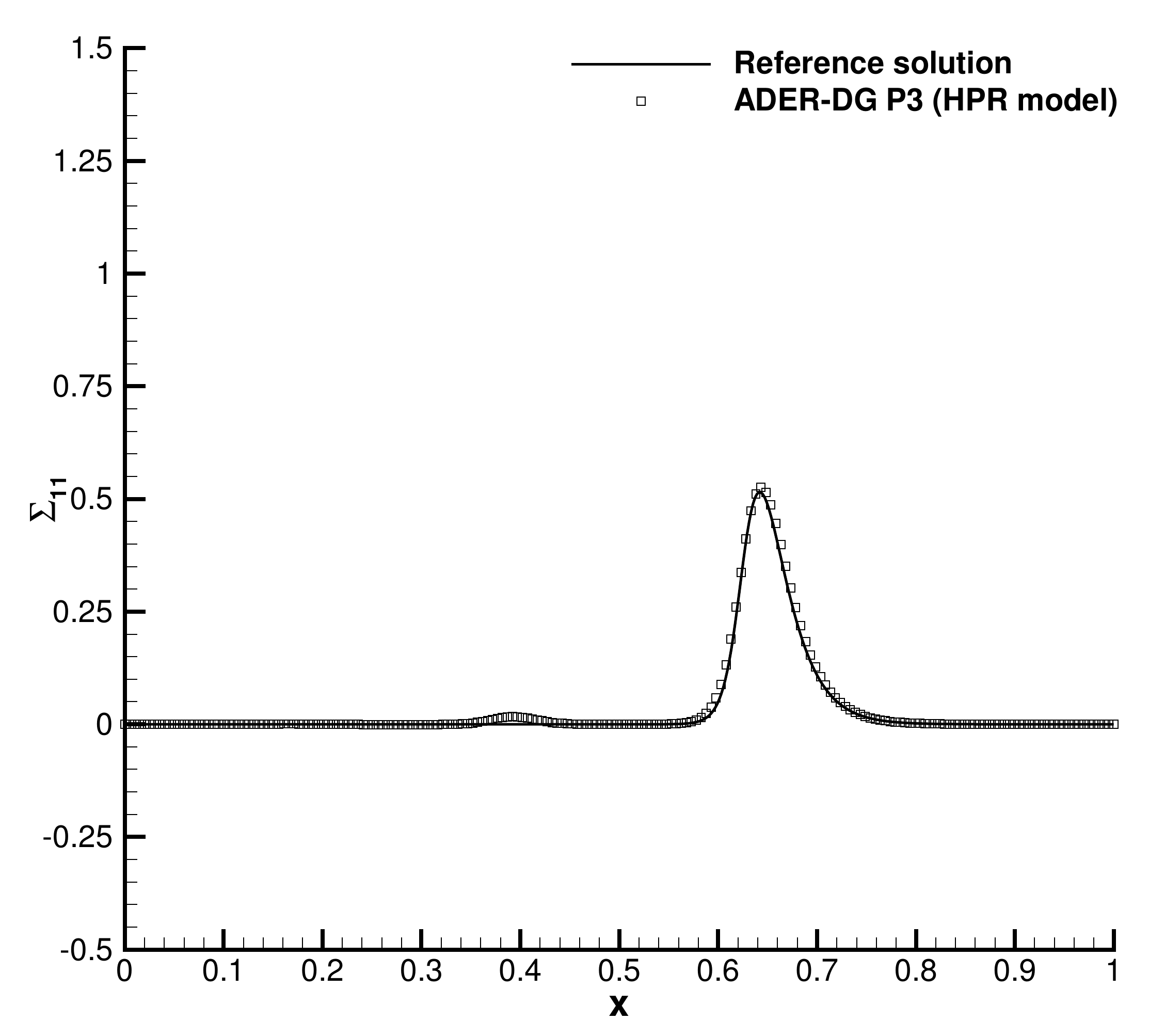}   
	\end{tabular} 
	\caption{Viscous heat conducting shock. Comparison of the exact solution of the compressible Navier-Stokes equations with the numerical solution of the GPR model based on ADER-DG $P3$ schemes. 
	Density profile (top left), velocity profile (top right), heat flux (bottom left) and stress $\sigma_{11}$ (bottom right). }
	\label{fig.gpr} 
\end{figure}

In the Fig. \ref{fig.gpr} we show numerical results obtained in \cite{HPRmodel} for a viscous heat 
conducting shock wave and the comparison with the exact solution of the compressible 
Navier-Stokes equations. 

\subsection{The equations of ideal general relativistic magnetohydrodynamics (GRMHD)}  

A very challenging PDE system is given by the equations of ideal general relativistic magnetohydrodynamics
(GRMHD). The governing PDE are a result of the Einstein field equations and can be written in compact 
covariant notation as follows: 
\begin{equation}
\nabla_{\mu} T^{\mu \nu} =0, \qquad \textnormal{and} \qquad \nabla_{\mu} \,^*F^{\mu \nu} =0 \qquad
\textnormal{and} \qquad \nabla_{\mu} (\rho u^\mu) =0, 
\end{equation} 
where $\nabla_{\mu}$ is the usual covariant derivative operator, $T^{\mu \nu}$ is the energy-momentum
tensor,  $\,^*F^{\mu \nu}$ is the Faraday tensor and $u^\mu$ is the four-velocity. The compact 
equations  above can be expanded into a so-called 3+1 formalism, which can be cast into the form 
of \eqref{eqn.pde}, see \cite{Anton06,ZannaZanotti,ADERGRMHD} for more details.  The final evolution 
system involves 9 field variables plus the 4-metric of the background space-time, which is supposed
to be stationary here. A numerical convergence study for the large amplitude Alfv\'en wave test problem
described in \cite{ZannaZanotti} solved in the domain $\Omega=[0,2\pi]^3$ up to $t=1$ and carried out 
with high order ADER-DG schemes in \cite{ADERGRMHD} is reported in the Table \ref{tab.conv.comp} below, 
where we also show a direct comparison with high order Runge-Kutta discontinuous Galerkin schemes. 
We observe that the ADER-DG schemes are competitive with RKDG methods, even for this very complex 
system of hyperbolic PDE. \textcolor{black}{The results reported in Table \ref{tab.conv.comp} refer to the 
non-vectorized version of the code.} Further significant performance improvements are expected from a 
carefully vectorized implementation of the GRMHD equations, \textcolor{black}{in particular concerning the
vectorization of the cumbersome conversion of the vector of conservative variables to the vector 
of primitive variables, i.e. the function $\mathbf{V}=\mathbf{V}(\mathbf{Q})$. For the GRMHD system $\mathbf{V}$ cannot be computed 
analytically in terms of $\mathbf{Q}$, but requires the iterative solution of one nonlinear scalar algebraic 
equation together with the computation 
of the roots of a third order polynomial, see \cite{ZannaZanotti} for details. In our vectorized implementation
of the PDE, we have therefore in particular vectorized the primitive variable recovery via a direct implementation  
in AVX intrinsics. We have furthermore made use of careful auto-vectorization via the compiler for the 
evaluation of the physical flux function 
and for the non-conservative product. Thanks to this vectorization effort, on one single CPU core of an Intel 
i9-7900X Skylake test workstation with 3.3 GHz nominal clock frequency and using AVX 512 the CPU time necessary for a 
single degree of freedom update (TDU) for a fourth order ADER-DG scheme ($N=3$)  
could be reduced to TDU = 2.3 $\mu$s for the GRMHD equations in three space dimensions. } 

\begin{table*} 
	\centering
	\begin{tabular}{cccc|cccc}
		\hline
		$N_x$ & $L_2$ error &  $L_2$ order & WCT [s] &  $N_x$ & $L_2$ error &  $L_2$ order & WCT [s]   \\ 
		\hline
		\multicolumn{4}{c|}{ADER-DG ($N=3$)} & 	\multicolumn{4}{c}{RKDG ($N=3$)} \\ 
		\hline
		8   & 7.6396E-04	&      	& 0.093 	& 8 	& 8.0909E-04	&     	& 0.107 	 \\ 
		16	& 1.7575E-05	&  5.44	& 1.371 	& 16	& 2.2921E-05	& 5.14	& 1.394 	 \\ 
		24	& 6.7968E-06	&  2.34	& 6.854 	& 24	& 7.3453E-06	& 2.81	& 6.894 	 \\ 
		32	& 1.0537E-06	&  6.48	& 21.642	& 32	& 1.3793E-06	& 5.81	& 21.116	 \\  
		\hline
		\multicolumn{4}{c|}{ADER-DG ($N=4$)} & 	\multicolumn{4}{c}{RKDG ($N=4$)} \\ 
		\hline	
		8   & 6.6955E-05	&      	& 0.363 	& 8 	& 6.8104E-05	&     	& 0.456 	 \\ 
		16	& 2.2712E-06	&  4.88	& 5.696 	& 16	& 2.3475E-06	& 4.86	& 6.666 	 \\ 
		24	& 3.3023E-07	&  4.76	& 28.036	& 24	& 3.3731E-07	& 4.78	& 29.186	 \\ 
		32	& 7.4728E-08	&  5.17	& 89.271	& 32	& 7.7084E-08	& 5.13	& 87.115	 \\  
		\hline
		\multicolumn{4}{c|}{ADER-DG ($N=5$)} & 	\multicolumn{4}{c}{RKDG ($N=5$)} \\ 
		\hline	
		8   & 5.2967E-07	&      	& 1.090  	& 8 	& 5.7398E-07	&     	& 1.219  	 \\ 
		16	& 7.4886E-09	&  6.14	& 16.710 	& 16	& 8.1461E-09	& 6.14	& 17.310 	 \\ 
		24	& 7.1879E-10	&  5.78	& 84.425 	& 24	& 7.7634E-10	& 5.80	& 83.777 	 \\ 
		32	& 1.2738E-10	&  6.01	& 263.021	& 32	& 1.3924E-10	& 5.97	& 260.859	 \\  
		\hline	
	\end{tabular}
	\caption{ \label{tab.conv.comp} Accuracy and cost comparison between ADER-DG and RKDG schemes of different orders for the GRMHD equations in three space dimensions. The errors refer to the variable $B_y$. The table also contains  
		total wall clock times (WCT) measured in seconds using 512 MPI ranks of the SuperMUC phase I system at the LRZ in Garching, Germany. } 
\end{table*}

\begin{figure*}
	\begin{center}
		\includegraphics[width=0.4\textwidth]{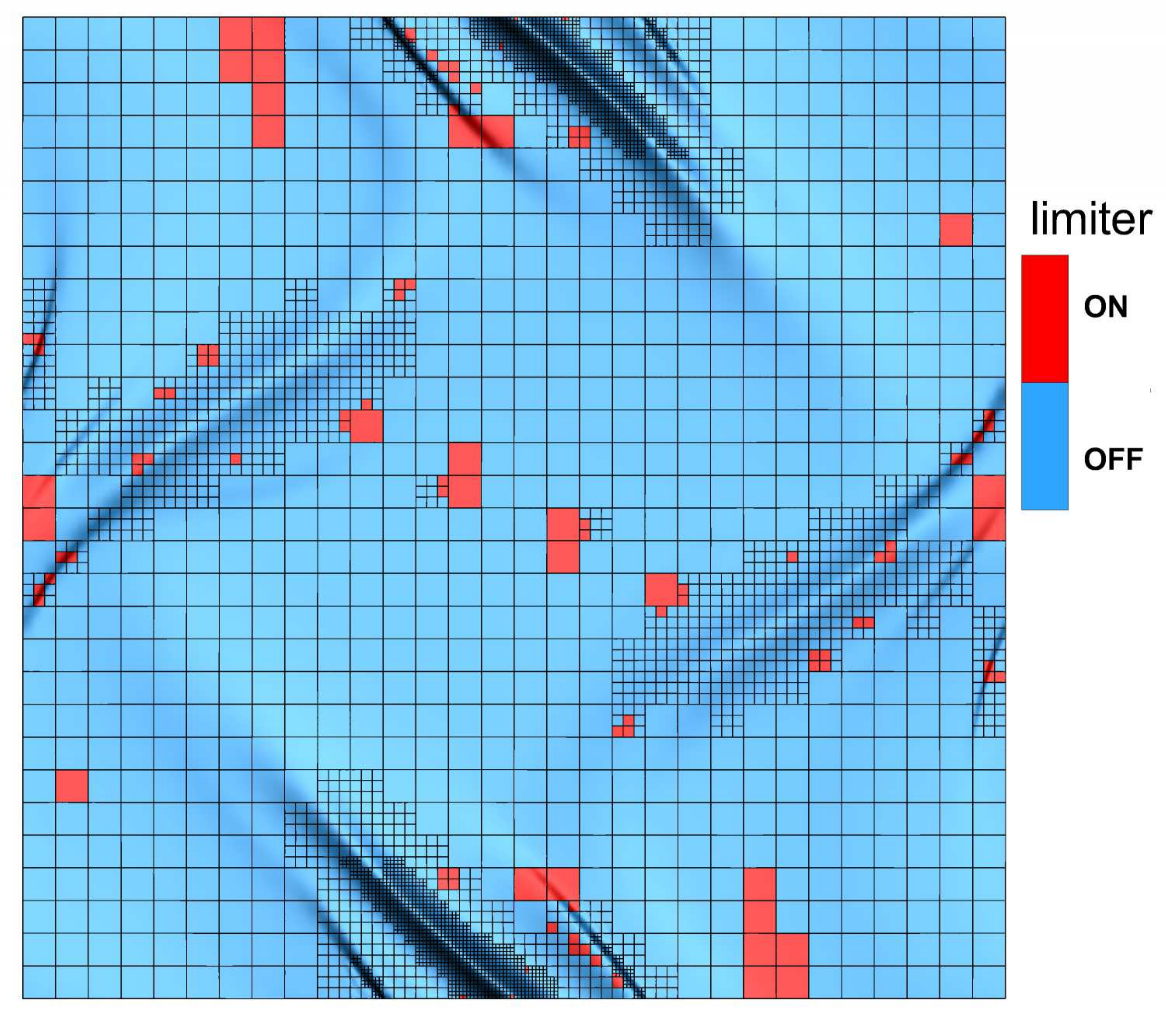} 
		\includegraphics[width=0.4\textwidth]{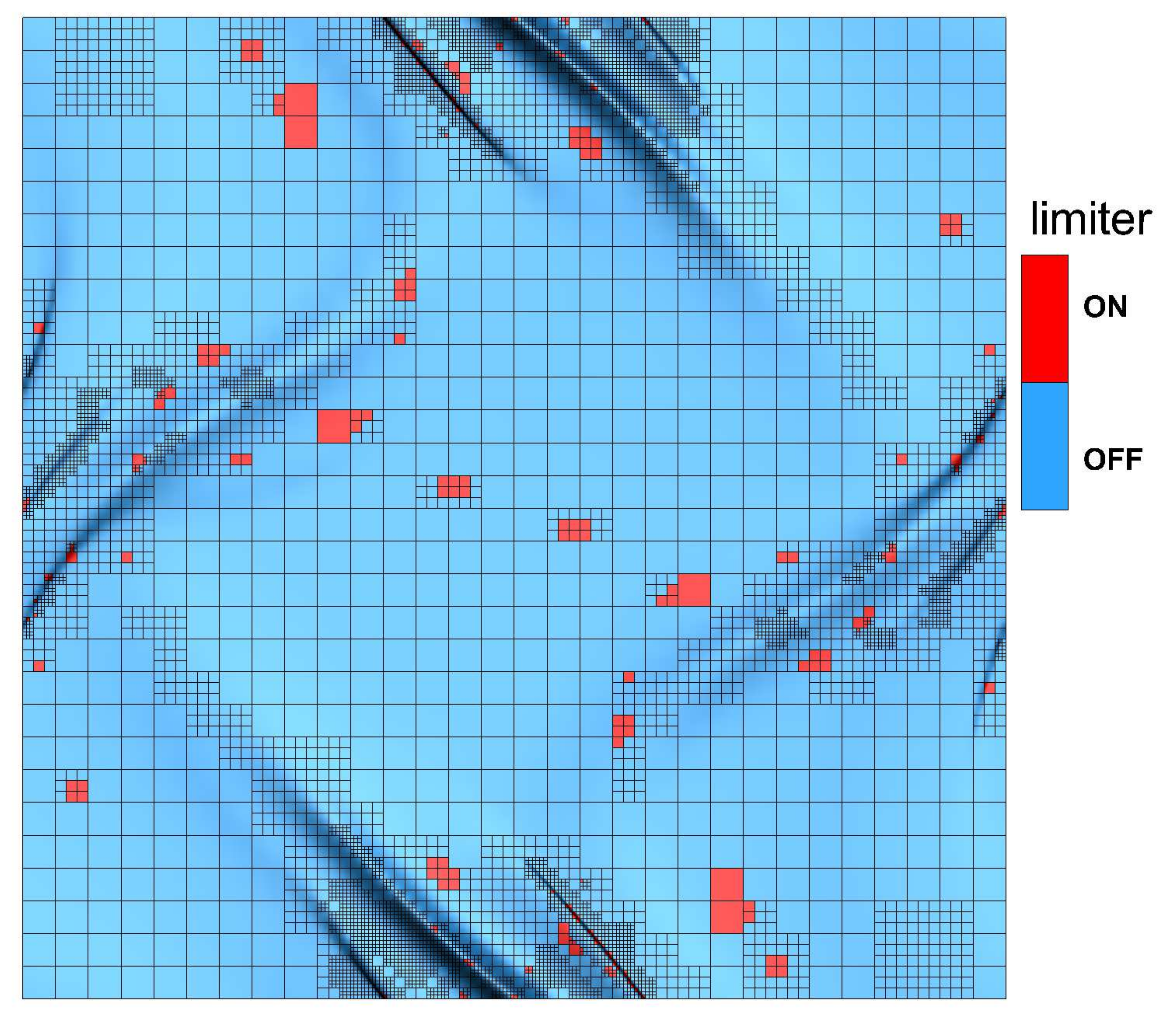}  \\  
		\includegraphics[width=0.4\textwidth]{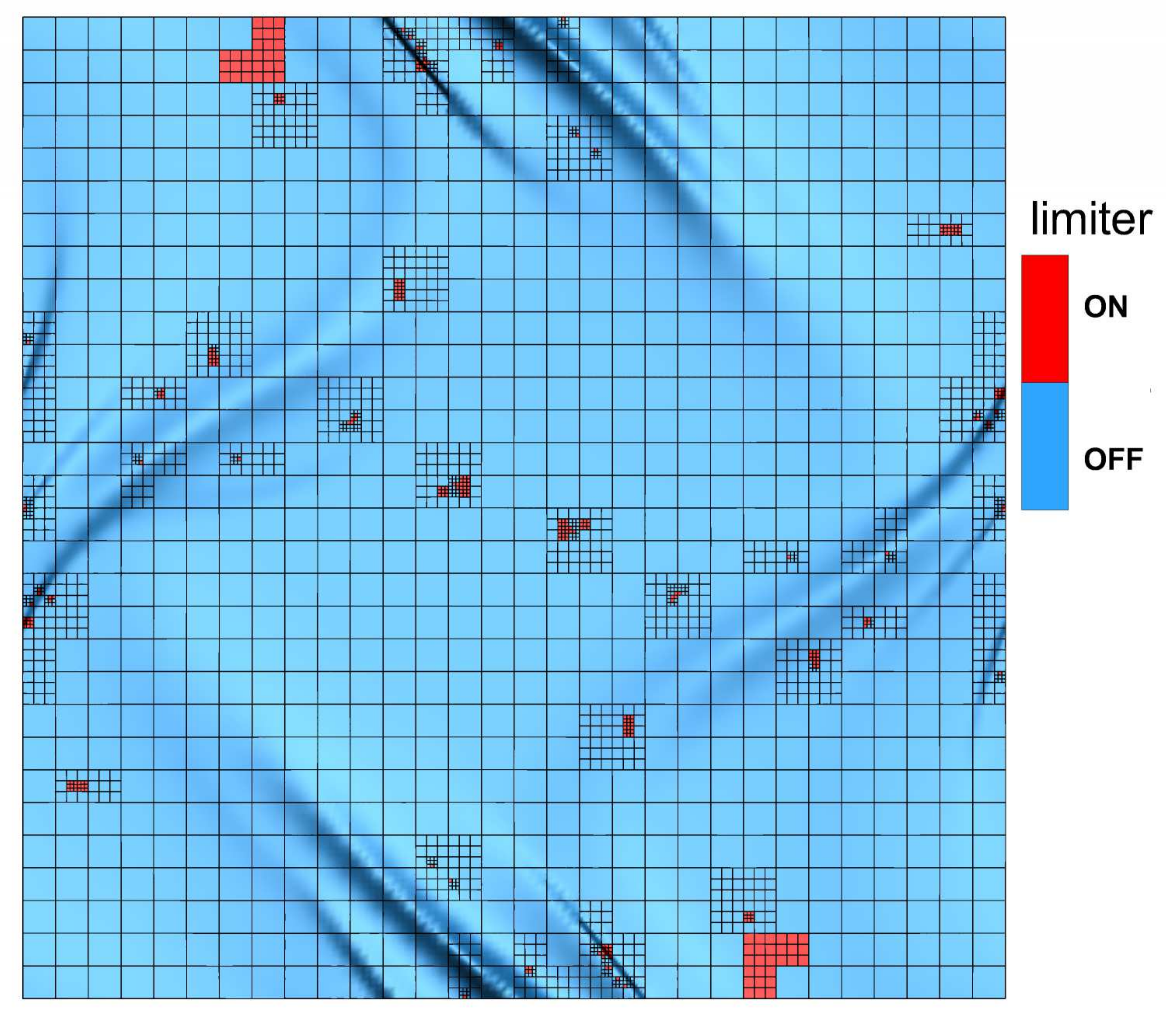}  
		\includegraphics[width=0.4\textwidth]{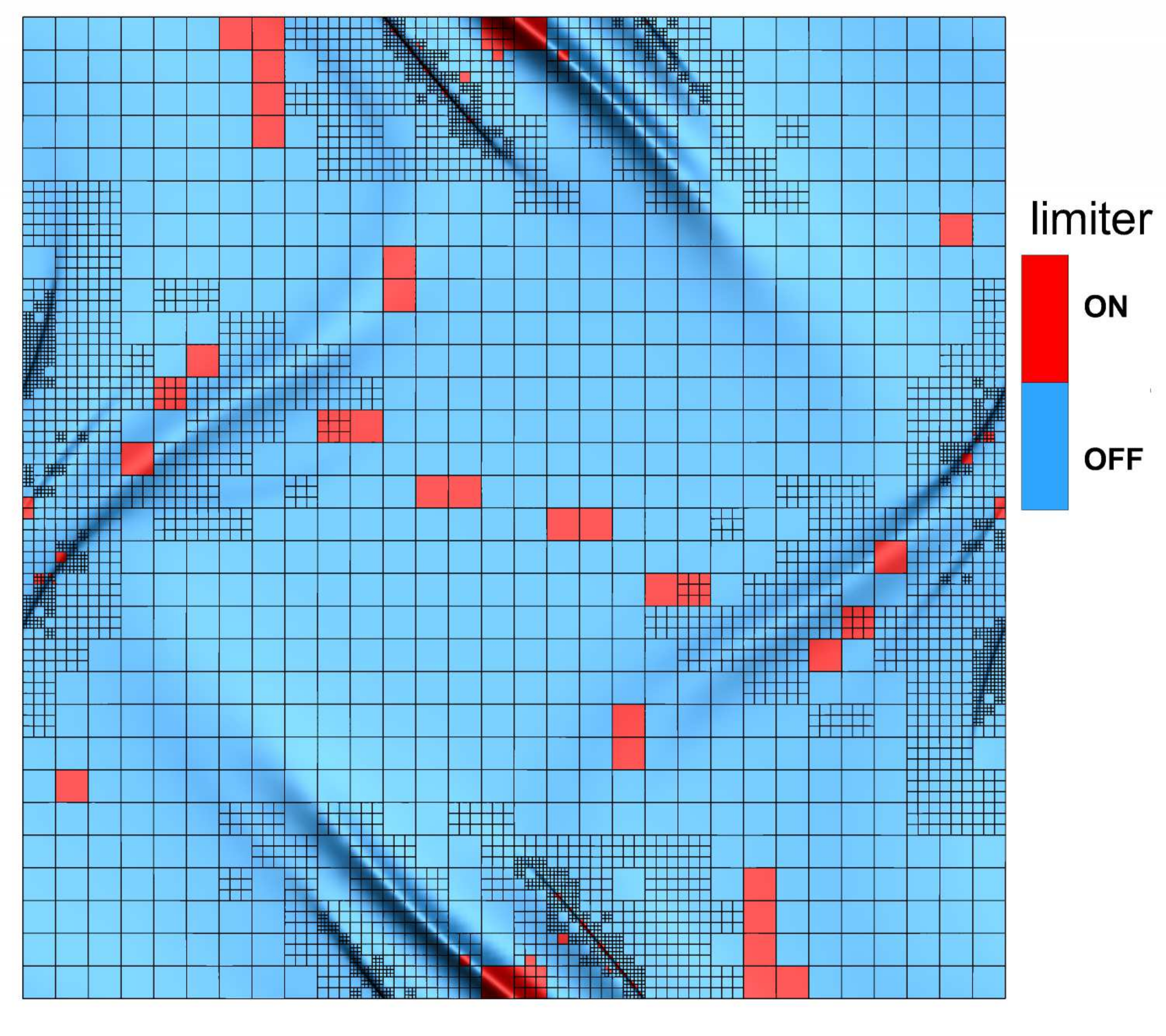}  
		\caption{Results for the GRMHD Orszag-Tang vortex problem in flat space-time (SRMHD) at  
			$t=2$ obtained with ADER-DG-$\mathbb{P}_5$ schemes, supplemented with \emph{a posteriori} 
			subcell finite volume limiter and using different refinement estimator functions $\chi$. (i) first  order-derivative estimator $\chi_1$ (top left); (ii) second-order derivative estimator $\chi_2$ (top right);  (iii) a new limiter-based estimator $\chi_3$  (row 2, left) and  (iv) a new multi-resolution estimator $\chi_4$ based on the difference between the discrete solution on two adjacent refinement levels (row 2 right). 
		}
		\label{fig:OrszagTang}
	\end{center}
\end{figure*}

\begin{figure*}
	\begin{center}
		\includegraphics[width=0.3\textwidth]{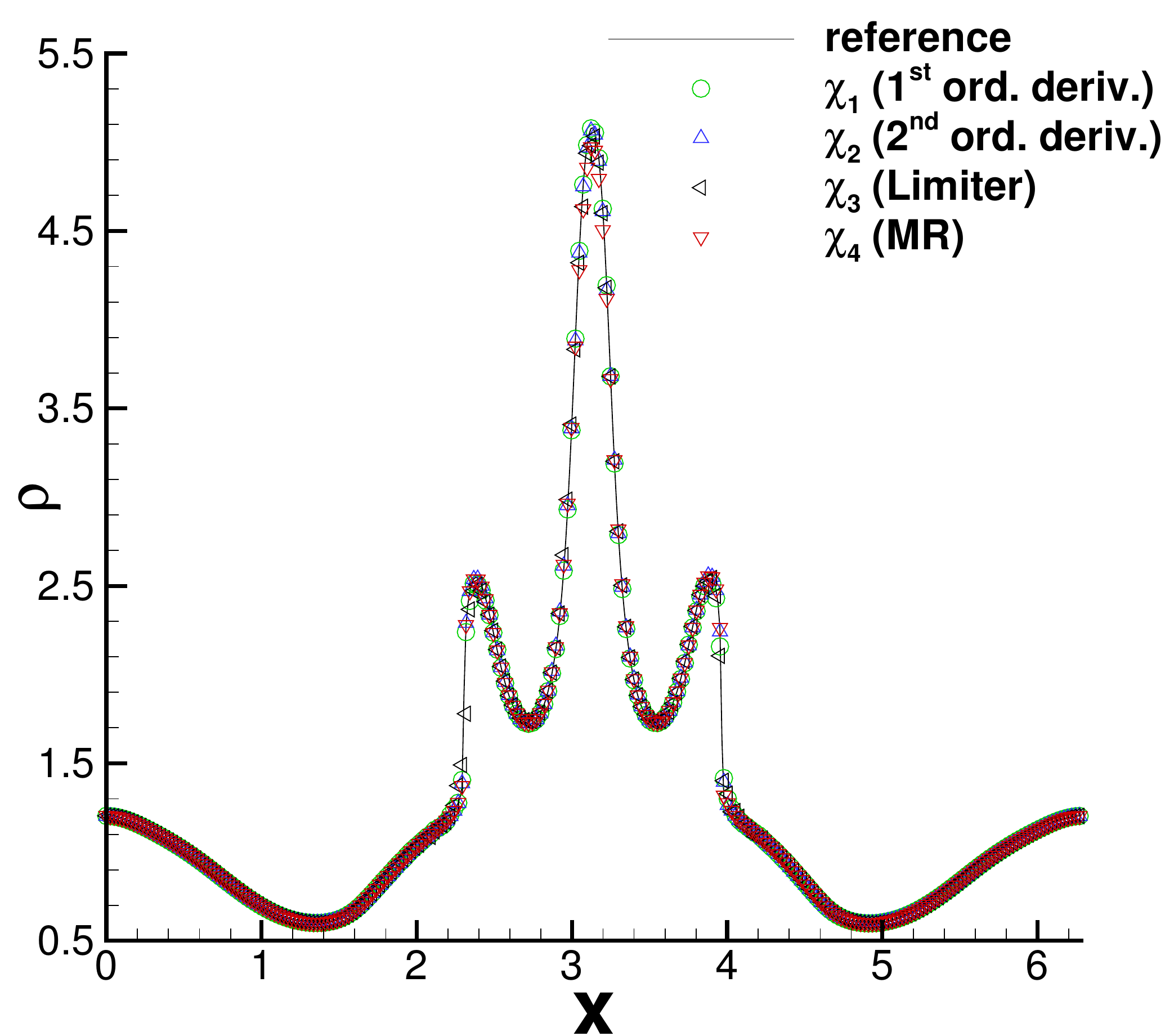} 
		\includegraphics[width=0.3\textwidth]{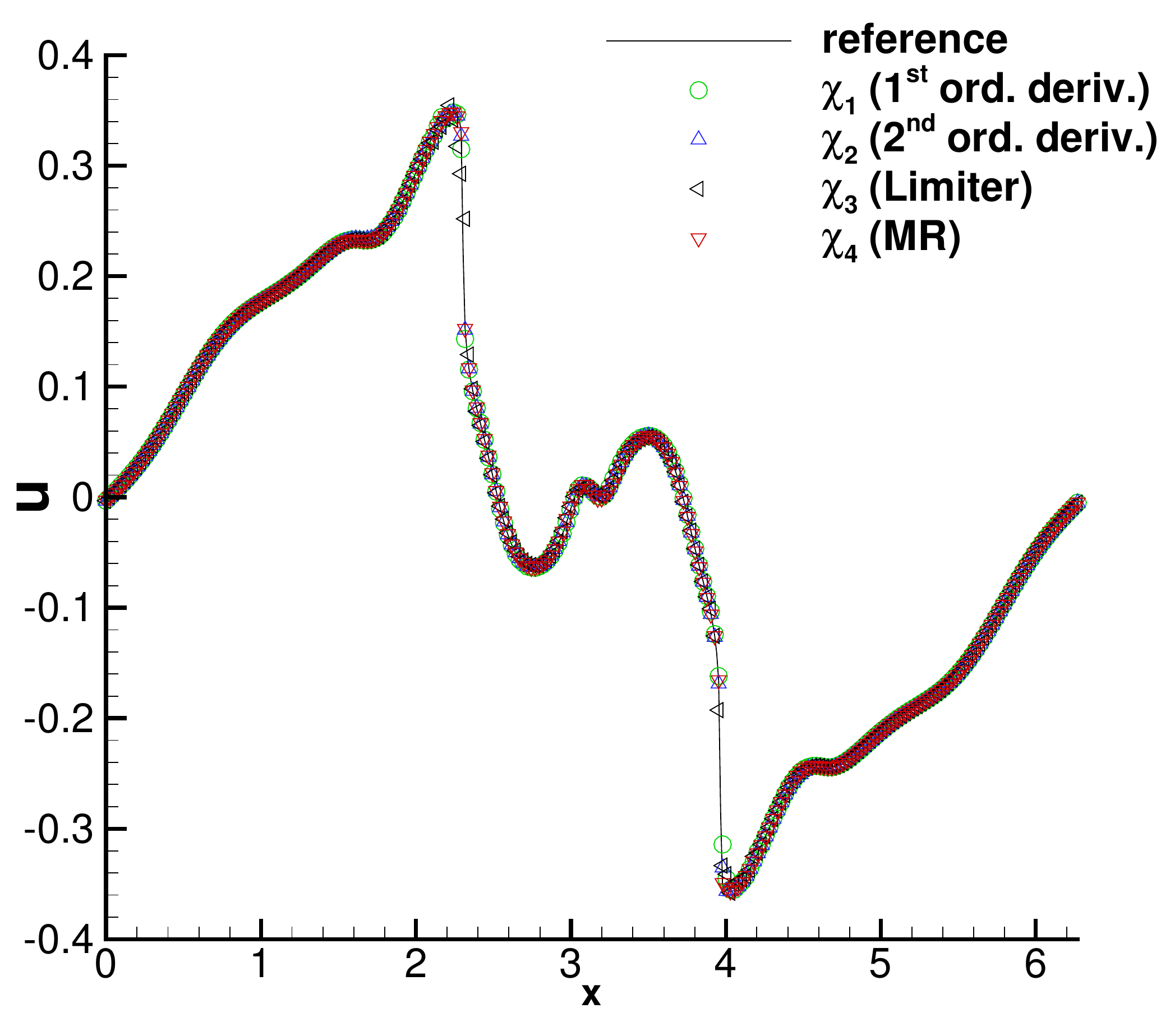}  
		\includegraphics[width=0.3\textwidth]{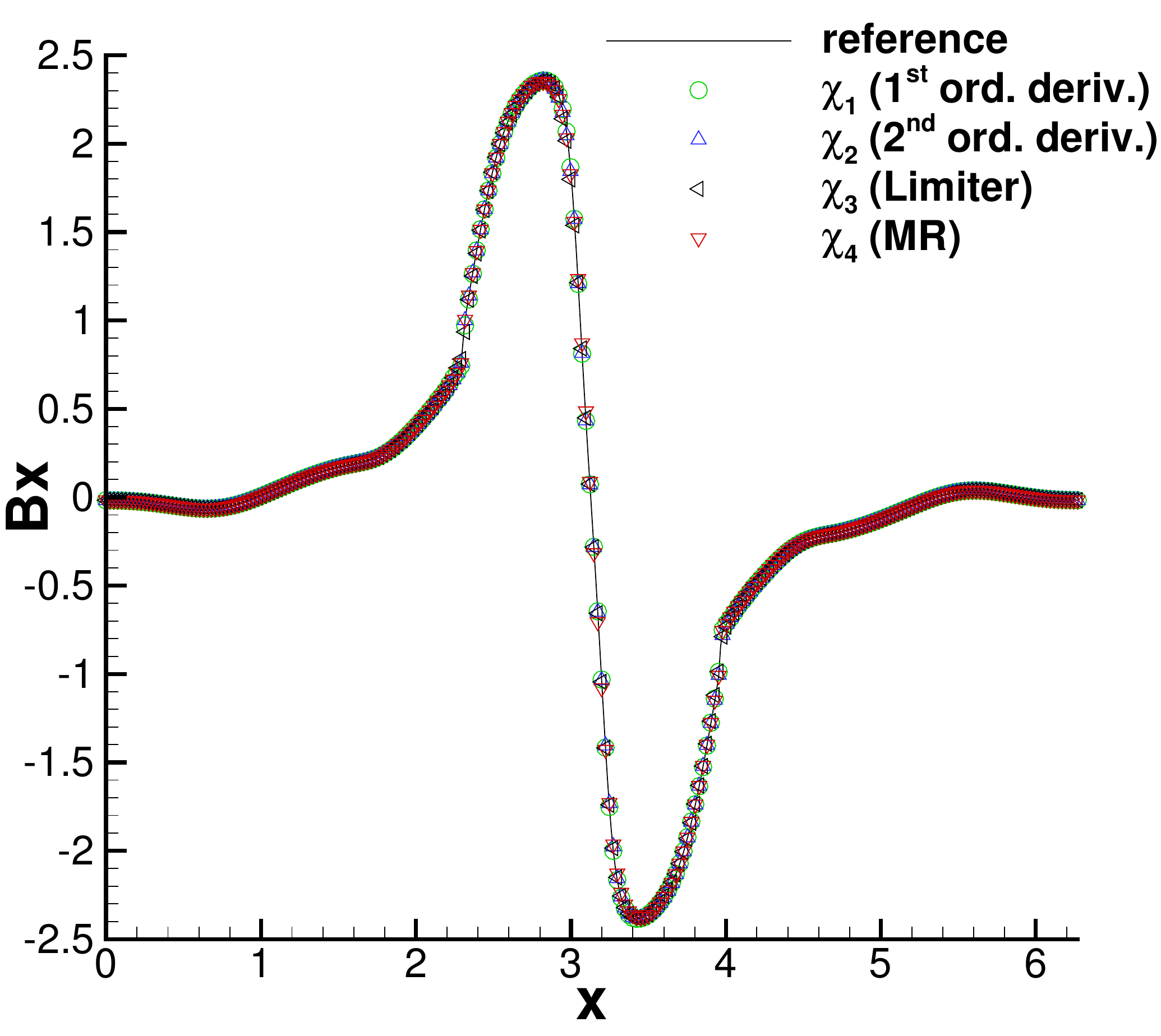} 
		\caption{Results for the GRMHD Orszag-Tang vortex problem in flat space-time (SRMHD) at  
			$t=2$ obtained with ADER-DG-$\mathbb{P}_5$ schemes supplemented with \emph{a posteriori} 
			subcell finite volume limiter and using different refinement estimator functions $\chi$. 
			A set of 1D cuts taken at $y=10^{-2}$ are shown. From left to right: the rest-mass density, the  velocity $u$ and the magnetic field component $B_x$. One can note an excellent agreement between the reference solution and the ones obtained on different AMR grids. }
		\label{fig:OrszagTang.cut}
	\end{center}
\end{figure*}

As second test problem we present the results obtained for the Orszag-Tang vortex system in flat Minkowski spacetime,
where the GRMHD equations reduce to the special relativitic MHD equations. The initial condition is given by
\begin{align*}
	&\left( \rho, u, v, w, p, B_x ,B_y,B_z \right) = \left( 1 , -
	\frac{3}{4\sqrt{2}}\sin y\,, \frac{3}{4\sqrt{2}}\sin x \,, 0, 1, - \sin
	y\,, \sin 2x \,, 0 \right)\,,
\end{align*}
and we set the adiabatic index to $\Gamma=4/3$. The computational domain is $\Omega = [0,2\pi]^2$
and is discretized with a dynamically adaptive AMR grid. 
For this test we chose the $P_5$ version of the ADER-DG scheme with FV subcell limiter and the 
rest mass density as indicator function for AMR, i.e. $\varphi(\Q)=\rho$. 
Fig. \ref{fig:OrszagTang.cut} shows 1D cuts through the numerical solution at time $t=2$ and at 
$y=0.01$, while Fig. \ref{fig:OrszagTang} shows the numerical results for the AMR-grid 
with limiter-status map (blue cells are unlimited, while limited cells are highlighted
in red), together with Schlieren images for the rest-mass density at time $t=2$. The same simulation 
has been repeated with different refinement estimator functions $\chi$ that tell the AMR algorithm 
where and when to refine and to coarsen the mesh: (i) a simple first order derivative estimator $\chi_1$
based on discrete gradients of the indicator function $\varphi(\Q)$, (ii)  the classical 
second order derivative estimator $\chi_2$ based on \cite{Loehner1987}, (iii) a novel estimator 
$\chi_3$ based on the action of the \textit{a posteriori} subcell finite volume limiter, i.e. the mesh is
refined where the limiter is active (iv)  a multi-resolution estimator $\chi_4$ based on the difference
in $L_\infty$ norm of the discrete solution on two different refinement levels $\ell$ and $\ell-1$.    
The reference solution is obtained on a uniform fine grid corresponding to the finest refinement level, i.e. 
a uniform composed of $270\times 270$ elements. The results shown in Fig. \ref{fig:OrszagTang} clearly 
show that the  numerical results obtained by means of different refinement estimator functions are comparable 
with each other and thus the proposed AMR algorithm is robust with respect to the particular
choice of the mesh.  

\textcolor{black}{As last test case we simulate a stationary neutron star in three space dimensions using the 
	Cowling approximation, i.e. assuming a fixed \textit{static} background spacetime. 
The initial data for the matter and the spacetime are both compatible with the Einstein field equations and are 
given by the solution of the Tolman–Oppenheimer–Volkoff (TOV) equations, which constitute a nonlinear 
ODE system in the radial coordinate that  can be numerically solved up to any precision at the aid of a fourth order Runge-Kutta 
scheme using a very fine grid. We setup a stable nonrotating TOV star without magnetic field and with central rest mass density 
$\rho(\mathbf{0},0)=1.28 \cdot 10^{-3}$ and adiabatic exponent $\Gamma=2$ in a computational domain $\Omega=[-10,+10]^3$ 
discretized with a fourth order ADER-DG scheme ($N=3$) using $32^3$ elements, which corresponds to $128^3$ spatial 
degrees of freedom. The pressure in the atmosphere outside the compact object is set to $p_{atm}=10^{-13}$. We run the 
simulation until a final time of $t=1000$ and measure the $L_\infty$ error norms of the rest mass density and the pressure 
against the exact solution, which is given by the initial condition. The error at $t=1000$ for the rest mass 
density is  $L_\infty(\rho)=1.553778 \cdot 10^{-5}$ while the error for the pressure is   
$L_\infty(p)=1.605334 \cdot 10^{-7}$. The simulation was carried out with the vectorized version of the code on 512 CPU  
cores of the SuperMUC phase 2 system (based on AVX2) and required only 3010 s of wallclock time. The same simulation with the established  
finite difference GR code WhiskyTHC \cite{WhiskyTHC} required 8991 s of wall clock time on the same machine with the same spatial mesh 
resolution and the same number of CPU cores. The time series of the relative error of the central rest mass density in the origin of the domain 
is plotted in the left panel of Fig. \ref{fig:TOVseries}. 
At the final time $t=1000$, the relative error of the central rest mass density is still below 0.1\%. In the right panel of Fig. \ref{fig:TOVseries} 
we show the contour surfaces of the pressure at the final time $t=1000$. 
In Fig. \ref{fig:TOVcut} we show a 1D cut along the $x$ axis, comparing the numerical solution at time $t=1000$ with the exact one. 
We note that the numerical scheme  
is very accurate, but it is \textit{not well-balanced} for the GRMHD equations, i.e. the method \textit{cannot} preserve the 
stationary equilibrium solution of the TOV equations \textit{exactly} at the discrete level. Therefore, further work along 
the lines of research reported recently in \cite{GaburroDumbserEuler} for the Euler equations with Newtonian gravity are 
needed, extending the framework of well-balanced methods \cite{Bermudez1994,Castro2006,Pares2006} also to general relativity. 
Finally, in Fig. \ref{fig:TOV3Dxy} we compare the exact and the numerical solution at  time $t=1000$ in the $x-y$ plane. 
\begin{figure}
	\begin{tabular}{cc}
		\includegraphics[width=0.45\textwidth]{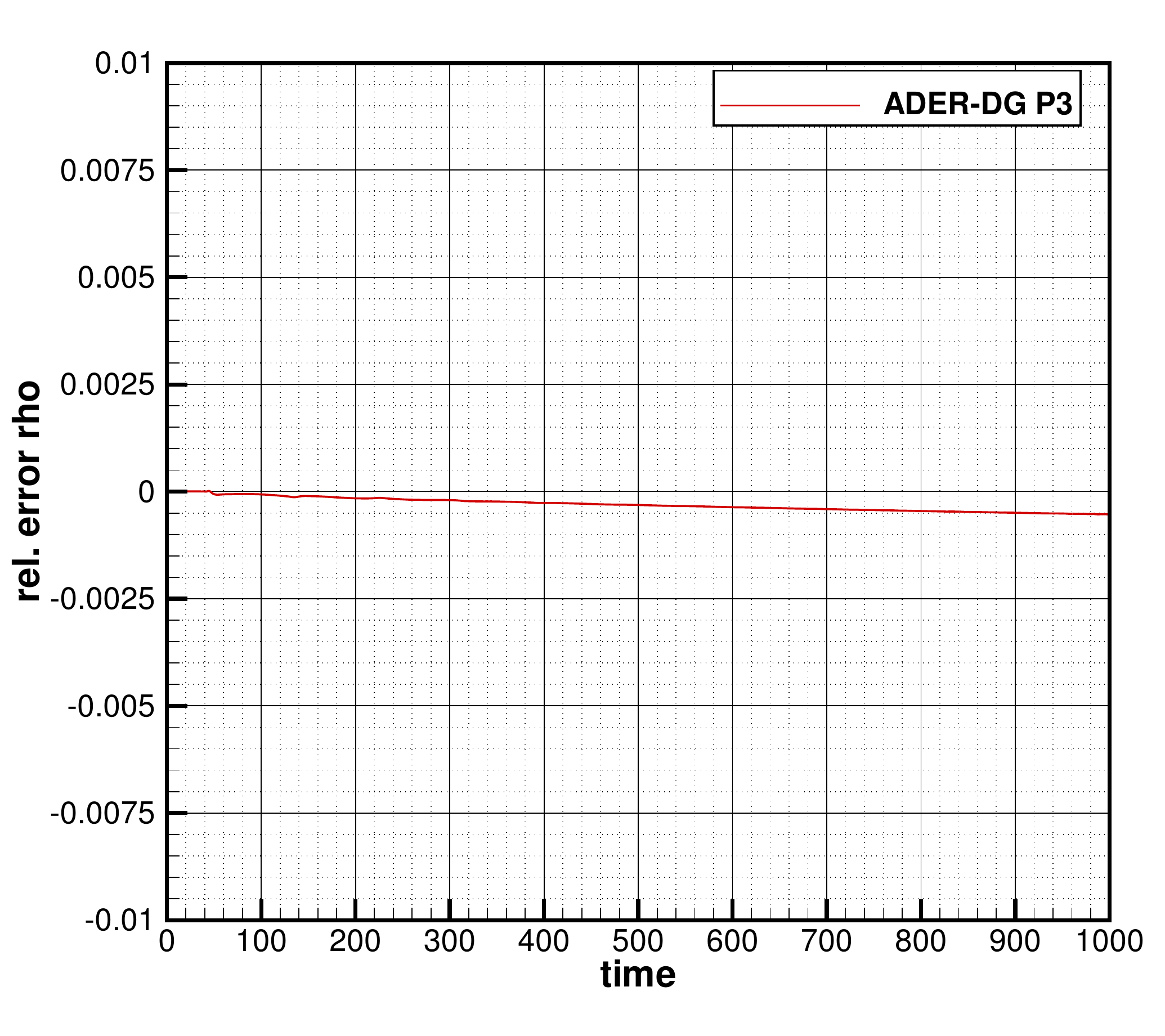} & 
		\includegraphics[width=0.45\textwidth]{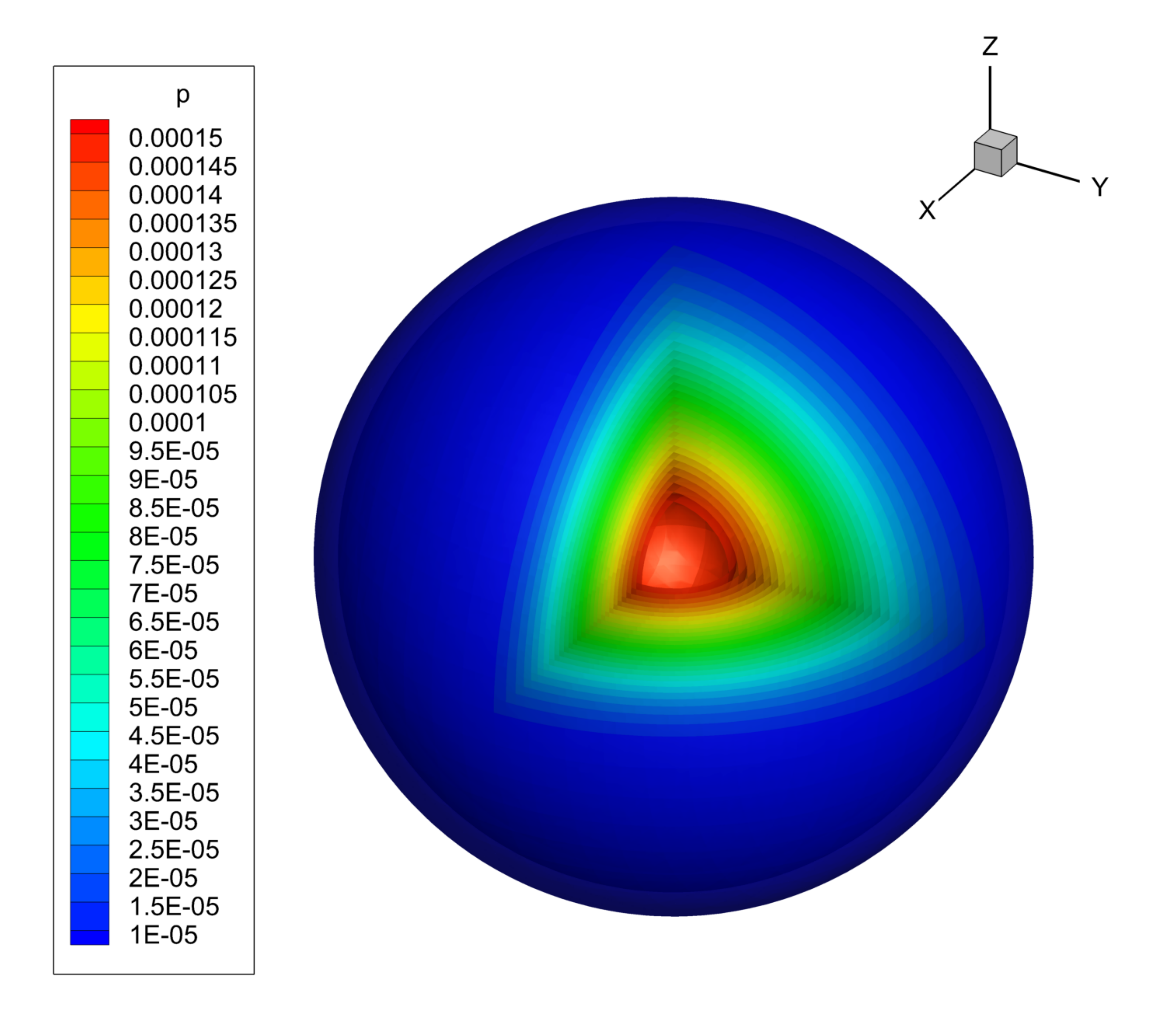}  
		\caption{Computational results for a stable 3D neutron star. 
			     Time series of the relative error of the central rest mass density $\left(\rho(\mathbf{0},t)-\rho(\mathbf{0},0) \right) / \rho(\mathbf{0},0)$ (left) and 3D view of of the pressure contour surfaces at time $t=1000$ (right). } 
		\label{fig:TOVseries}
	\end{tabular}
\end{figure}
\begin{figure}[!t]
	\begin{tabular}{cc}
		\includegraphics[width=0.45\textwidth]{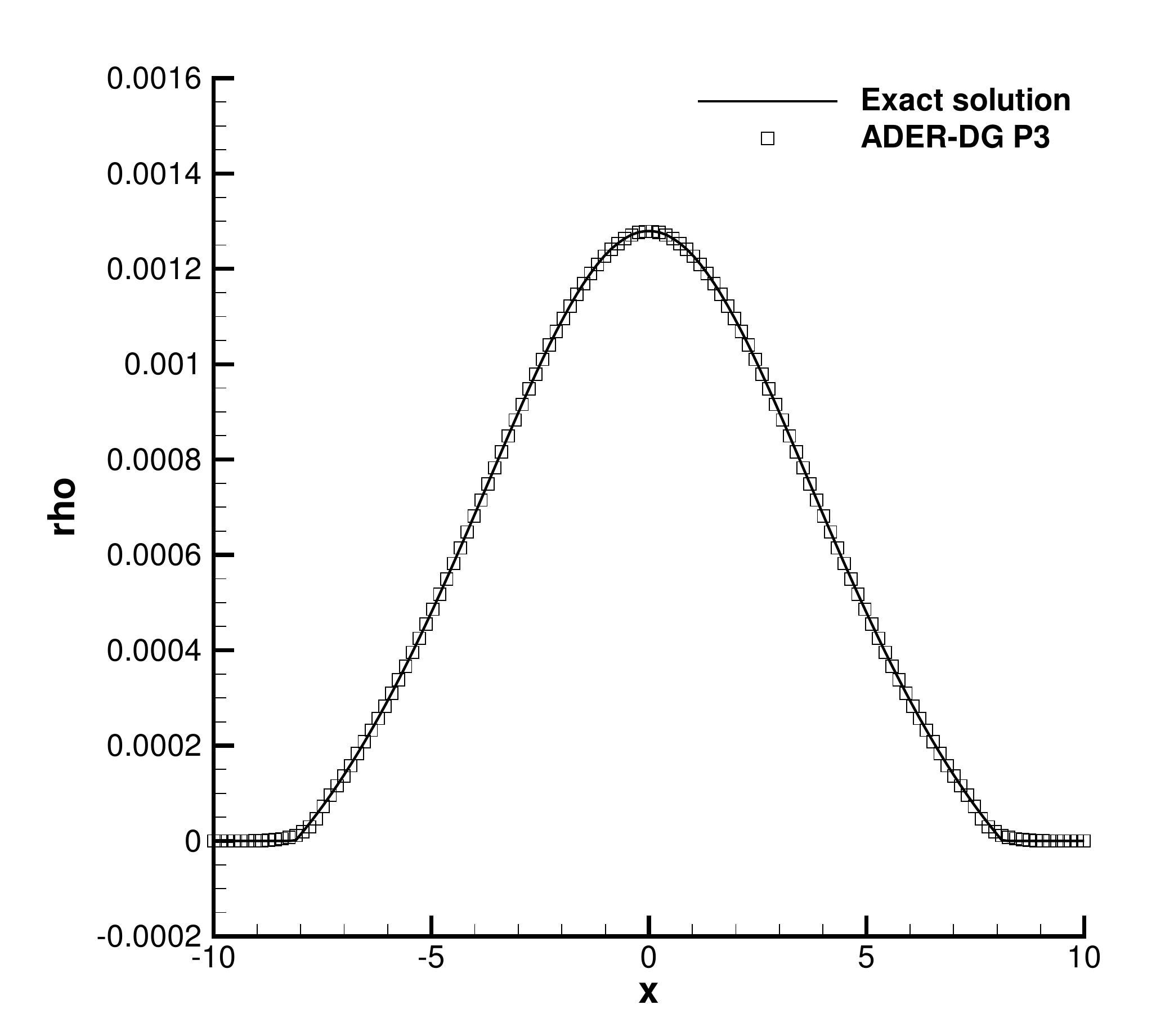} & 
		\includegraphics[width=0.45\textwidth]{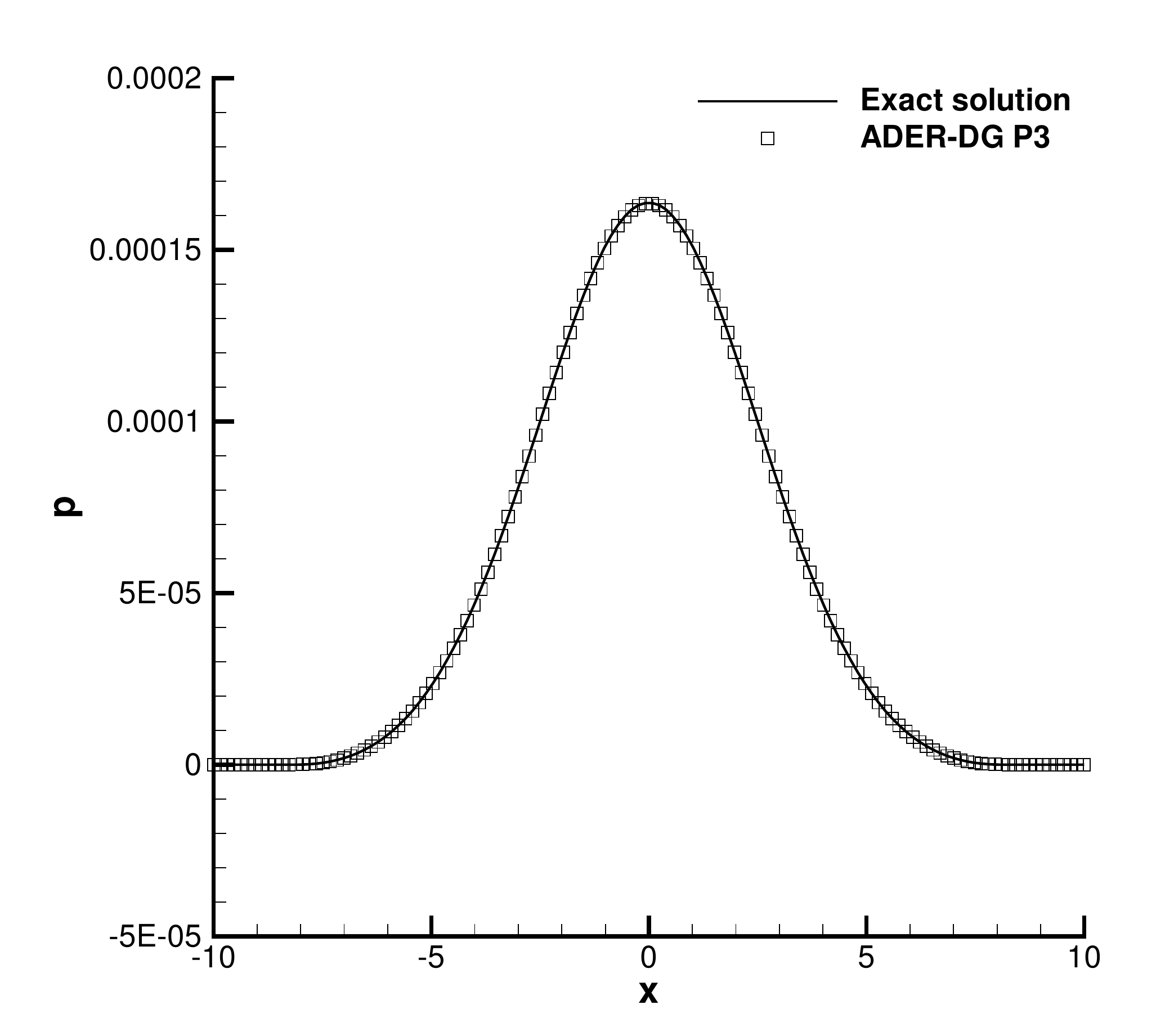}  
		\caption{Computational results for a stable 3D neutron star. Comparison of the numerical solution with the exact one at time $t=1000$ 
			on a 1D cut along the $x$-axis for the rest mass density (left) and the pressure (right).} 
		\label{fig:TOVcut}
	\end{tabular}
\end{figure}
\begin{figure}[!t]
	\begin{tabular}{cc}
		\includegraphics[width=0.45\textwidth]{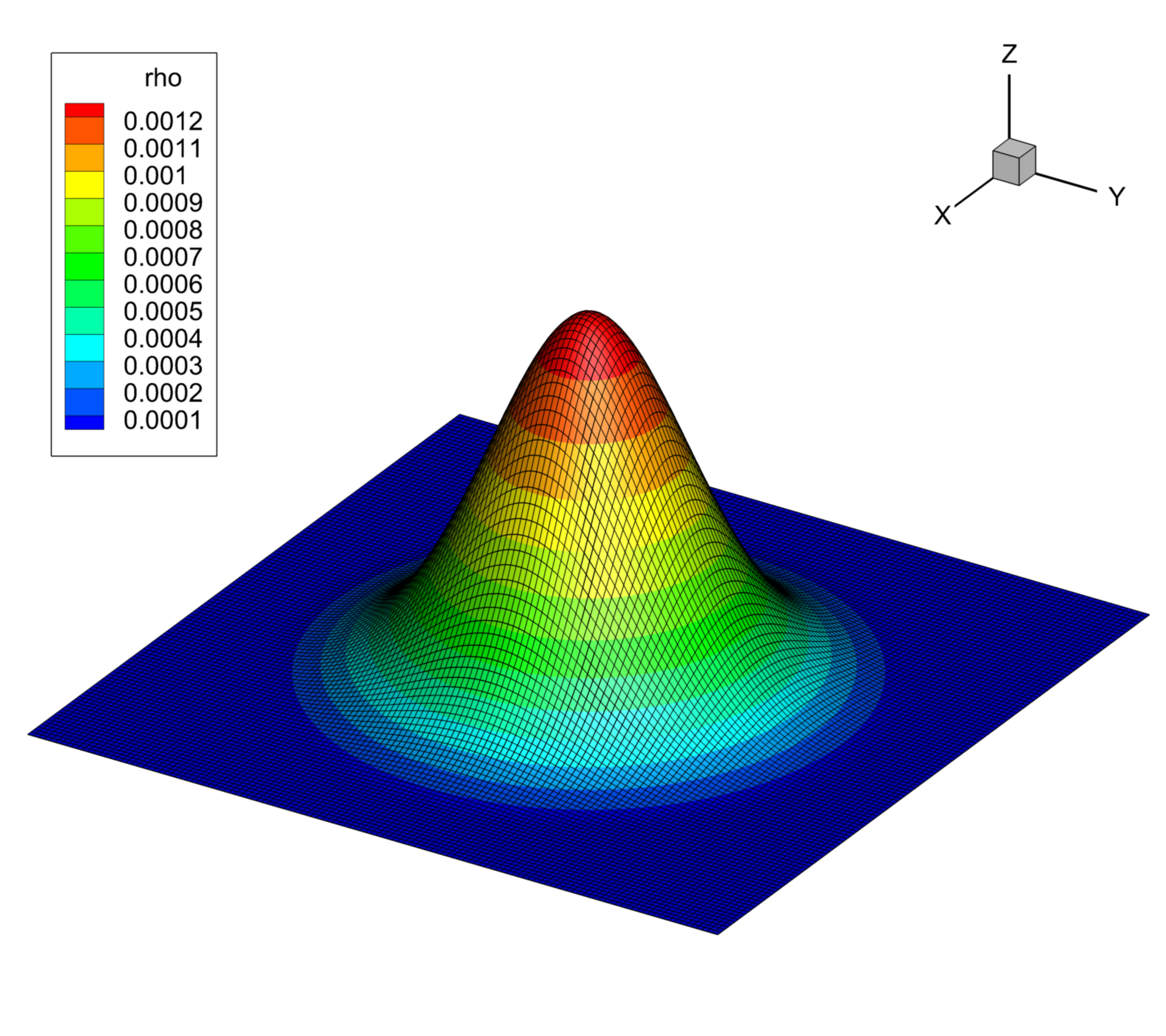} & 
		\includegraphics[width=0.45\textwidth]{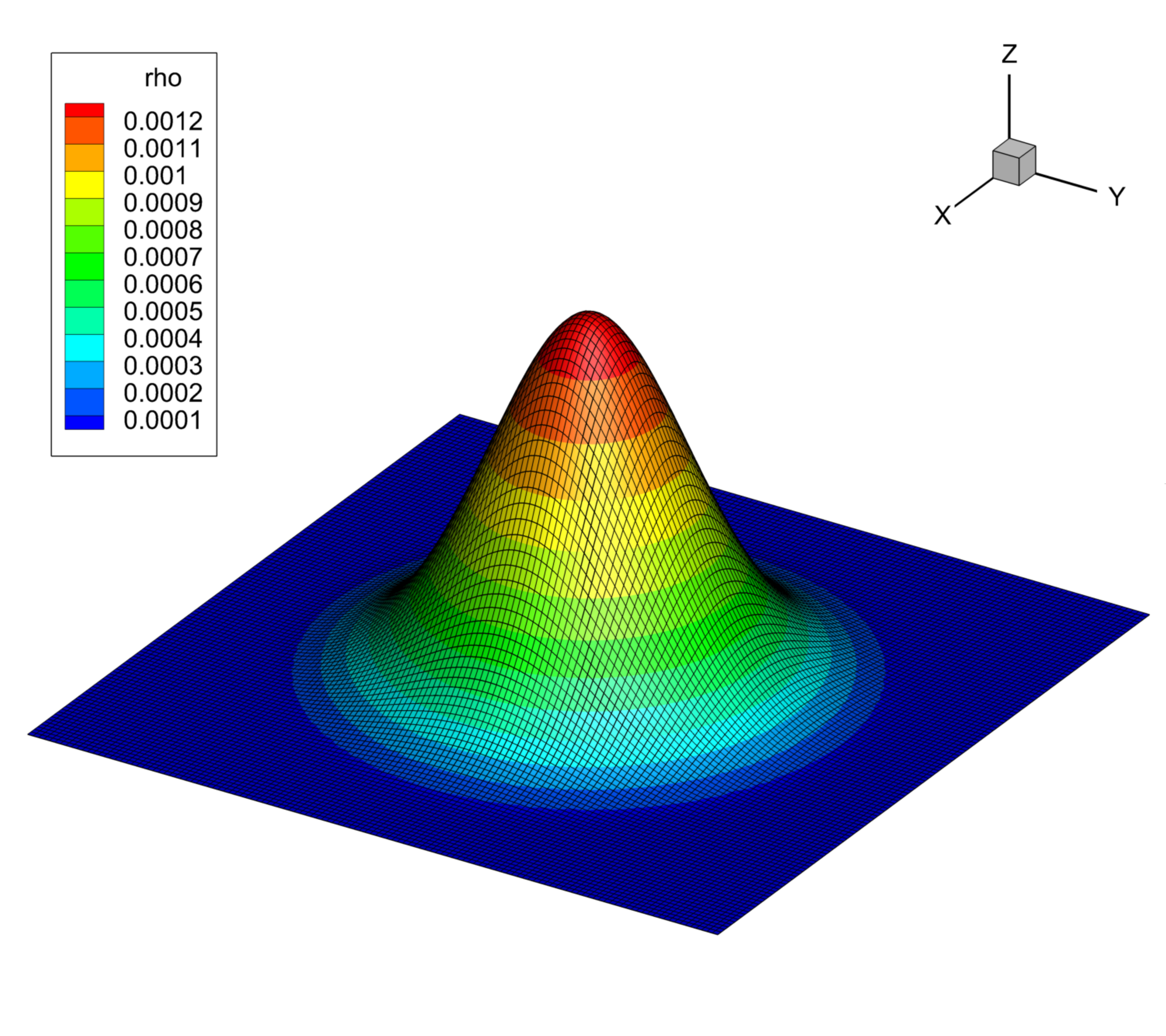}  
		\caption{Computational results for a stable 3D neutron star. Cut through the $x-y$ plane with pressure on the $z$ axis   
			and rest mass density contour colors. Exact solution (left) and numerical solution at time $t=1000$ (right).}  
		\label{fig:TOV3Dxy}
	\end{tabular}
\end{figure}
} 
	
\subsection{A strongly hyperbolic first order reduction of the CCZ4 formulation of the Einstein field equations (FO-CCZ4) } 	
	
The last PDE system under consideration here are the Einstein field equations that describe the evolution of dynamic spacetimes. Here we consider the so-called CCZ4 formulation \cite{Alic:2011a}, which is based
on the Z4 formalism that takes into account the involutions (stationary differential constraints) 
inherent in the Einstein equations via an augmented system similar to the generalized Lagrangian multiplier (GLM) 
approach of Dedner et al. \cite{Dedneretal} that takes care of the stationary divergence-free constraint of the magnetic
field in the MHD equations. In compact covariant notation the undamped Z4 Einstein equations in vacuum, which can be
derived from the Einstein-Hilbert action integral associated with the Z4 Lagrangian $\mathcal{L}=g^{\mu \nu} \left( R_{\mu \nu} + 2 \nabla_\mu Z\nu \right)$, read 
\begin{equation} 
   R_{\mu \nu} + \nabla_{(\mu} Z_{\nu)} =0,    
\end{equation}  	
where $g^{\mu \nu}$ is the 4-metric of the spacetime, $R_{\mu \nu}$ is the 4-Ricci tensor and the 4-vector $Z_\nu$ accounts for the stationary constraints of
the Einstein equations, as already mentioned before. 
After introducing the usual 3+1 ADM split of the 4-metric as 
\begin{equation}
 ds^2 = -\alpha^2 dt^2 + \gamma_{ij} \left( dx^i + \beta^i dt \right) \left( dx^j + \beta^j dt \right), 
\end{equation}
the equations can be cast into a time-dependent system of 25 partial differential equations that involve 
first order derivatives in time and both first and second order derivatives in space,  
see \cite{Alic:2011a}.  Nevertheless, the system is \textit{not} dissipative, but a rather unusual 
formulation of a wave equation, see \cite{Gundlach04a}. In the expression above, $\alpha$ denotes the so-called lapse, $\beta^i$ is the spatial shift vector and $\gamma_{ij}$ is the spatial metric. In the original form presented in 
\cite{Alic:2011a}, the PDE system does \textit{not} fit into the formalism given by Eqn. \eqref{eqn.pde}. 
After the introduction of 33 auxiliary 
variables, which are the spatial gradients of some of the 25 primary evolution quantities, it is possible
to derive a first order reduction of the system that contains a total of 58 evolution quantities. 
However, a naive procedure of converting the original 
second order evolution system into a first order system leads only to a \textit{weakly hyperbolic} formulation,
which is not suitable for numerical simulations since the initial value problem is not well posed in this case. 
Only after adding suitable first and second order ordering constraints, which arise from the
definition of the auxiliary variables, it is possible to obtain a provably strongly hyperbolic
and thus well-posed evolution system, denoted by FO-CCZ4 in the following. For all details of the derivation,
the strong hyperbolicity proof and numerical results achieved with high order ADER-DG schemes, the reader is referred to \cite{ADERCCZ4}. In order to give an idea about the complexity of the Einstein field equations, it should be mentioned that one single evaluation of the FO-CCZ4 system requires about 20,000 floating point operations! In order to obtain still 
a computationally efficient implementation, the entire PDE system has been carefully vectorized using blocks
of the size \texttt{VECTORLENGTH}, so that in the end a level of 99,9 \% of vectorization of the 
 code have been reached. Using a fourth order ADER-DG scheme ($N=3$) the time per degree of freedom update (TDU) 
 metric per core on a modern workstation with Intel i9-7900X CPU that supports the novel AVX 512 instructions 
 is TDU=4.7 $\mu$s.      

\section{Strong MPI scaling study for the FO-CCZ4 system}

A major focus of this paper is the efficient implementation of ADER-DG schemes for high performance computing (HPC) on massively parallel distributed memory supercomputers. For this purpose, we have very recently carried  out a systematic study of the strong MPI scaling efficiency of our new high order fully-discrete one-step ADER-DG schemes on the Hazel Hen supercomputer of the HLRS center in Stuttgart, Germany, using from 720 up to 180,000 CPU cores. We have furthermore carried out a systematic comparison with conventional Runge-Kutta DG schemes using the SuperMUC phase 1 system of the LRZ center in Munich, Germany. 

As already discussed before, the particular feature of ADER-DG schemes compared to traditional Runge-Kutta DG  schemes (RKDG) is that they are intrinsically \textit{communication-avoiding} and \textit{cache-blocking}, 
which makes them particularly well suited for modern massively parallel distributed memory supercomputers. 
As governing PDE system for the strong scaling test the novel first-order reduction of the CCZ4 formulation 
of the 3+1 Einstein field equations has been been adopted \cite{ADERCCZ4}. We recall that FO-CCZ4 is a very  large nonlinear hyperbolic PDE system that contains 58 evolution quantities. 

\begin{figure}
	\centering
	\begin{tabular}{cc}
	\includegraphics[width=0.45\textwidth]{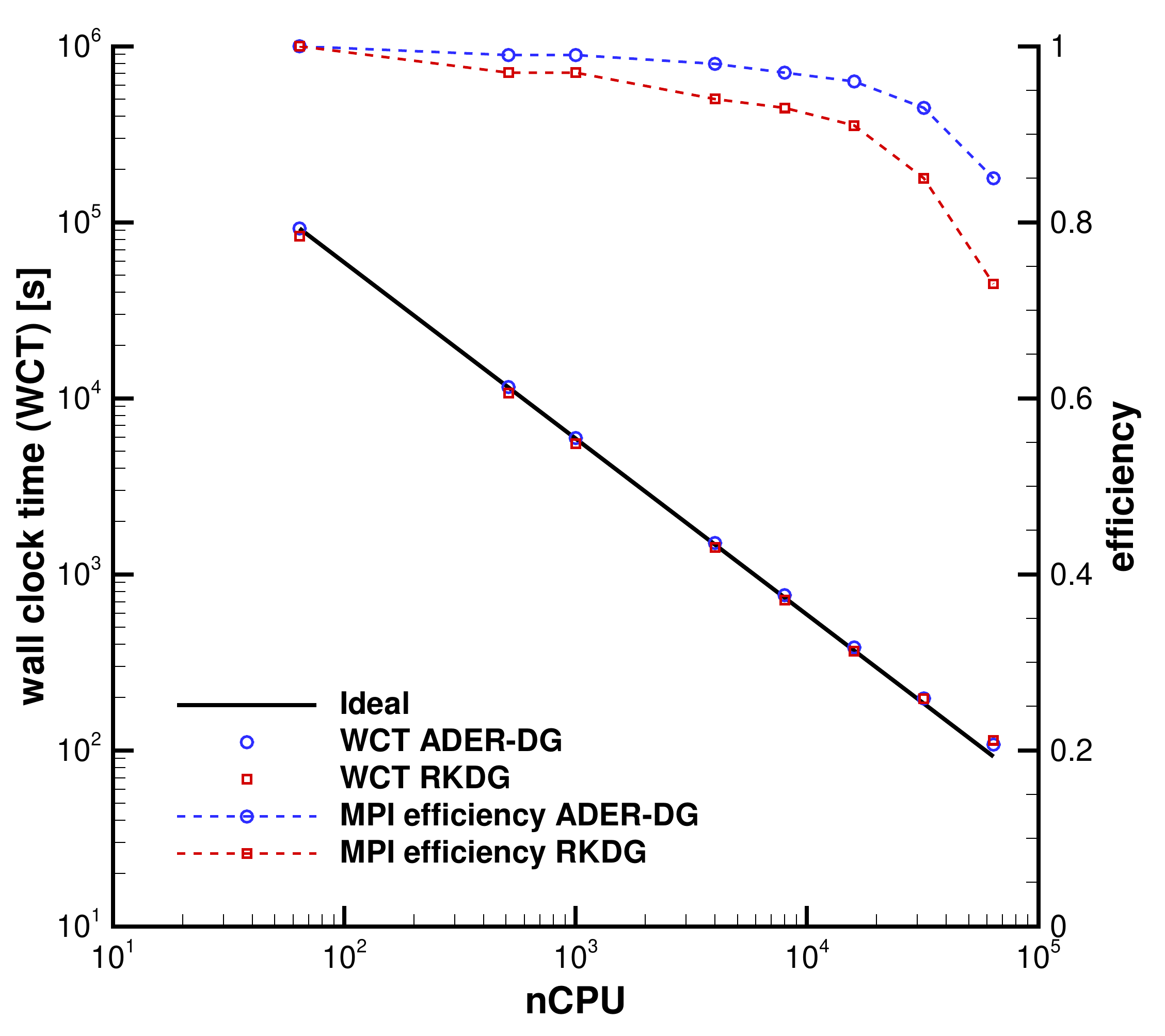} &  
	\includegraphics[width=0.46\textwidth]{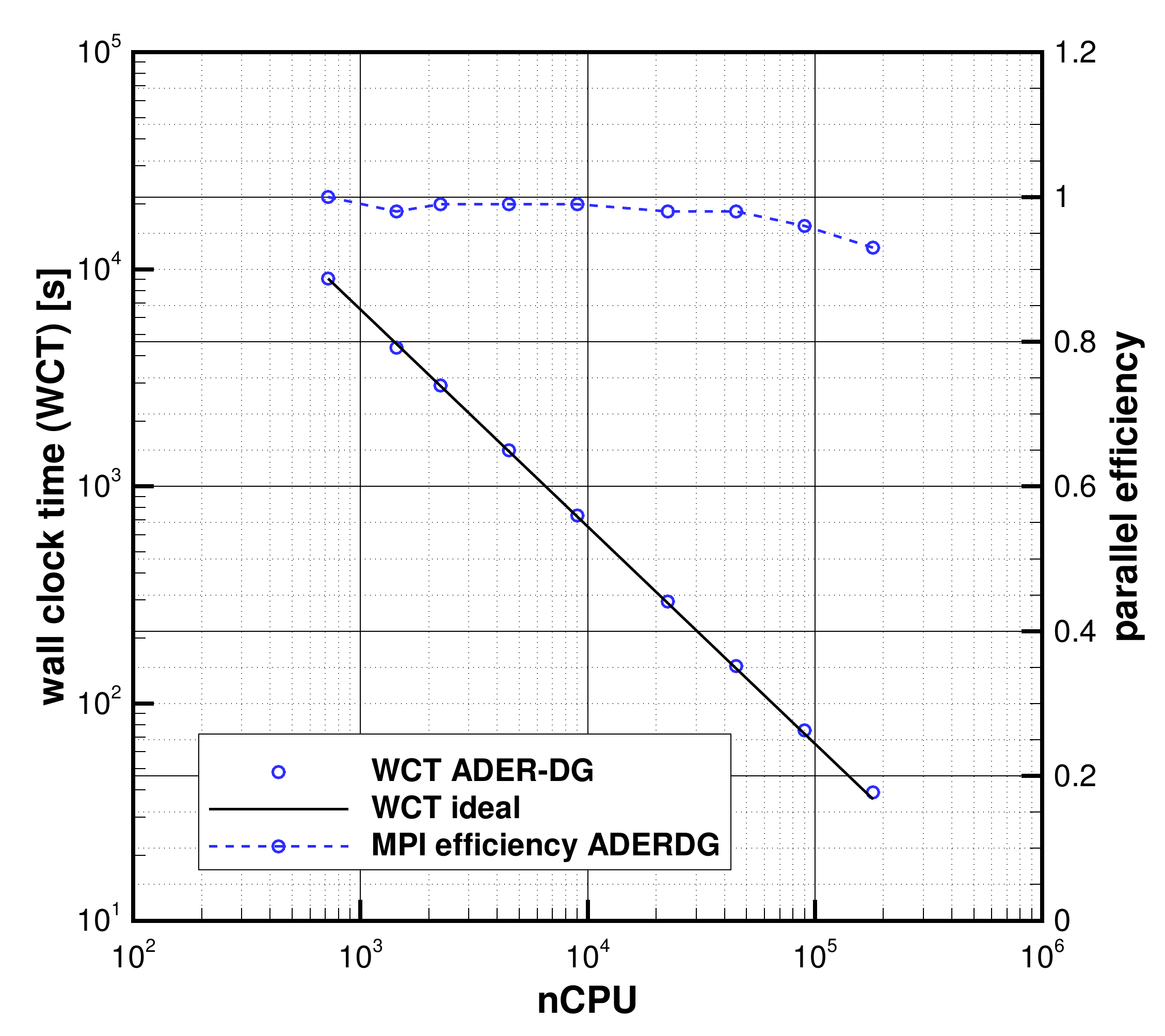}
	\end{tabular} 
	\caption{Strong MPI scaling study of ADER-DG schemes for the novel FO-CCZ4 formulation of the Einstein field 
		     equations recently proposed in \cite{ADERCCZ4}.  		  
		     Left: comparison of ADER-DG schemes with conventional Runge-Kutta DG schemes from 64 to 64,000 CPU 
		     cores on the SuperMUC phase 1 system of the LRZ supercomputing center (Garching, Germany). 
		     Right: strong scaling study from 720 to 180,000 CPU cores, including a full machine run on the Hazel Hen supercomputer of HLRS (Stuttgart, Germany) with ADER-DG schemes (right). Even on the full machine we observe still more than 90 \% of parallel efficiency. }
	\label{fig.scaling}
\end{figure}

The first strong scaling study on the SuperMUC phase 1 system uses 64 to 64,000 CPU cores. The test problem
 was the gauge wave problem \cite{ADERCCZ4} setup on the 3D domain $\Omega=[-0.5,0.5]^3$. For the test we have 
 compared a fourth order ADER-DG  scheme ($N=3$) with a 
 fourth order accurate RKDG scheme on a uniform Cartesian grid composed of $120^3$ elements. It has to 
 be stressed, that when using 64,000 CPU cores for this setup each CPU has to update only $3^3 = 27$ elements.  
 The wall clock time as a function of the used number of CPU cores (nCPU) and the obtained parallel
 efficiency with respect to an ideal linear scaling are reported in the left panel of Fig. \ref{fig.scaling}. 
 We find that ADER-DG schemes provide a better parallel efficiency than RKDG schemes, as expected. 

The second strong scaling study has been performed on the Hazel-Hen supercomputer, using 720 to 180,000 CPU cores. Again we have used a fourth order accurate ADER-DG scheme ($N=3$), this time using a uniform grid of $200\times 180\times 180$ elements, solving again the 3D gauge wave benchmark problem detailed in \cite{ADERCCZ4}.  
The measured wall-clock-times (WCT) as a function of the employed number of CPU cores, as well as the corresponding parallel scaling-efficiency are shown in  Fig. \ref{fig.scaling}. 
The results depicted in Fig. \ref{fig.scaling} clearly show that our new  ADER-DG schemes \textbf{scale very well}  up to 90,000 CPU cores with a parallel  efficiency \textbf{greater than}  \textbf{95\%}, and up to 180,000 cores with a  parallel efficiency that is still greater than \textbf{93\%}.    
Furthermore, the code was instrumented with manual FLOP counters in order to measure the floating point performance quantitatively. The full machine run on \textbf{180,000 CPU cores} of Hazel Hen took place on 7th of May 2018. During the run, each core has provided an average performance of 8.2 GFLOPS, leading to a total of \textbf{1.476 PFLOPS} of sustained performance. To our knowledge, this was the largest test run ever 
carried out with high order ADER-DG schemes for nonlinear hyperbolic systems of partial differential 
equations. For large runs with sustained petascale performance of ADER-DG schemes for linear hyperbolic 
PDE systems on unstructured tetrahedral meshes, see \cite{SeisSol2}. 

\section{Conclusions}

In this paper we have presented an efficient implementation of high order ADER-DG schemes on modern massively
parallel supercomputers using the ExaHyPE engine. 
The key ingredients are the communication-avoiding and cache-blocking properties of 
ADER-DG, together with an efficient vectorization of the high level user functions that provide 
the evaluation of the physical fluxes $\F(\Q)$, of the non-conservative products $\mathbf{B}(\Q) \cdot \nabla \Q$ and of the algebraic source terms $\S(\Q)$. The engine is highly versatile and flexible and allows to
solve a very broad spectrum of different hyperbolic PDE systems in a very efficient and highly 
scalable manner. In order to support this claim, we have provided a rather large set of different numerical examples solved with ADER-DG schemes. To show the excellent parallel scalability of the ADER-DG method, we have
provided strong scaling results on 64 to 64,000 CPU cores including a detailed and quantitative comparison 
with RKDG schemes. We have furthermore shown strong scaling results of the vectorized ADER-DG implementation 
for the FO-CCZ4 formulation of the Einstein field equations using 720 to 180,000 CPU cores of the Hazel 
Hen supercomputer at the HLRS in Stuttgart, Germany, where a sustained performance of more than one
petaflop has been reached. 

Future research in ExaHyPE will concern an extension of the GPR model to full general relativity, able 
to describe nonlinear elastic and plastic solids as well as viscous and ideal fluids in one single 
governing PDE system. We furthermore plan an implementation of the FO-CCZ4 system \cite{ADERCCZ4} 
directly based on AVX intrinsics, in order to further improve 
the performance of the scheme and to reduce computational time. The final aim of our developments are
the simulation of ongoing nonlinear dynamic rupture processes during earthquakes, as well 
as the inspiral and merger of binary neutron star systems and the associated generation of 
gravitational waves. Although both problems seem to be totally different and unrelated, it is indeed 
possible to write the mathematical formulation of both applications under the same form of a 
hyperbolic system given by \eqref{eqn.pde} and thus to solve both problems within the same computer 
software.  

\vspace{6pt} 

\authorcontributions{``Conceptualization and Methodology: M.D., M.B. and T.W.; Software: M.D., T.W., F.F. and M.T.; Validation: F.F. and M.T.; Strong MPI Scaling Study: F.F; Writing original draft and editing: M.D.; Supervision and Project Administration: M.B., T.W. and M.D.; Funding Acquisition: M.B., T.W. and M.D. ''. }

\funding{This research was funded by the European Union's Horizon 2020 Research and Innovation Programme under the 
	project \textit{ExaHyPE}, grant no. 671698 (call FETHPC-1-2014). The first three authors also acknowledge  funding from the Italian Ministry of Education, University and Research (MIUR) in the frame of the Departments of Excellence Initiative 2018--2022 attributed to DICAM of the University of Trento. MD also acknowledges support from the University of Trento in the frame of the Strategic Initiative \textit{Modeling and Simulation}. }

\acknowledgments{
The authors acknowledge the support of the HLRS supercomputing center in Stuttgart, Germany, for 
awarding access to the HazelHen supercomputer, as well as of the Leibniz Rechenzentrum (LRZ) in Munich, 
Germany for awarding access to the SuperMUC supercomputer. 
In particular, the authors are very grateful for the technical support provided by 
Dr. Bj\"orn Dick (HLRS) and Dr. Nicolay Hammer (LRZ). 
   
M.D., F.F. and M.T. are all members of the Italian INdAM Research group GNCS.	
}

\conflictsofinterest{The authors declare no conflict of interest.} 

%

\reftitle{References}

\bibliography{biblio}

\begin{thebibliography}{-------}
\providecommand{\natexlab}[1]{#1}

\bibitem[Riemann(1859)]{riemann1}
Riemann, B.
\newblock {\"Uber die Fortpflanzung ebener Luftwellen von endlicher
  Schwingungsweite}.
\newblock {\em G\"ottinger Nachrichten} {\bf 1859}, {\em 19}.

\bibitem[Riemann(1860)]{riemann2}
Riemann, B.
\newblock {\"Uber die Fortpflanzung ebener Luftwellen von endlicher
  Schwingungsweite}.
\newblock {\em Abhandlungen der K\"oniglichen Gesellschaft der Wissenschaften
  zu G\"ottingen} {\bf 1860}, {\em 8},~43--65.

\bibitem[Noether(1918)]{noether}
Noether, E.
\newblock {Invariante Variationsprobleme}.
\newblock {\em Nachrichten von der Gesellschaft der Wissenschaften zu
  G\"ottingen, mathematisch-physikalische Klasse} {\bf 1918}, pp. 235--257.

\bibitem[Courant \em{et~al.}(1928)Courant, Friedrichs, and Lewy]{CFL}
Courant, R.; Friedrichs, K.; Lewy, H.
\newblock {\"Uber die partiellen Differenzengleichungen der mathematischen
  Physik}.
\newblock {\em Mathematische Annalen} {\bf 1928}, {\em 100},~32--74.

\bibitem[Courant \em{et~al.}(1952)Courant, Isaacson, and Rees]{CIR}
Courant, R.; Isaacson, E.; Rees, M.
\newblock On the solution of nonlinear hyperbolic differential equations by
  finite differences.
\newblock {\em Comm. Pure Appl. Math.} {\bf 1952}, {\em 5},~243--255.

\bibitem[Courant and Hilbert(1962)]{Courant62}
Courant, R.; Hilbert, D.
\newblock {\em Methods of Mathematical Physics}; John Wiley and Sons, Inc: New
  York,  1962.

\bibitem[Courant and Friedrichs(1976)]{Courant62b}
Courant, R.; Friedrichs, K.O.
\newblock {\em Supersonic flows and shock waves}; Springer: Berlin,  1976.

\bibitem[Friedrichs(1958)]{FriedSymm}
Friedrichs, K.
\newblock Symmetric positive linear differential equations.
\newblock {\em Communications on Pure and Applied Mathematics} {\bf 1958}, {\em
  11},~333--418.

\bibitem[Godunov(1961)]{God1961}
Godunov, S.
\newblock An interesting class of quasilinear systems.
\newblock {\em Dokl. Akad. Nauk SSSR} {\bf 1961}, {\em 139(3)},~521--523.

\bibitem[Friedrichs and Lax(1971)]{FriedLax1971}
Friedrichs, K.; Lax, P.
\newblock Systems of conservation equations with a convex extension.
\newblock {\em Proceedings of the National Academy of Sciences} {\bf 1971},
  {\em 68},~1686--1688.

\bibitem[Godunov and Romenski(1995)]{Godunov:1995a}
Godunov, S.; Romenski, E.
\newblock {Thermodynamics, conservation laws, and symmetric forms of
  differential equations in mechanics of continuous media}.
\newblock  {Computational Fluid Dynamics Review 95}. {John Wiley, NY},  1995,
  pp. 19--31.

\bibitem[Godunov and Romenski(2003)]{Godunov:2003a}
Godunov, S.; Romenski, E.
\newblock {\em {Elements of Continuum Mechanics and Conservation Laws}}; Kluwer
  Academic/ Plenum Publishers,  2003.

\bibitem[Godunov(1972)]{God1972MHD}
Godunov, S.
\newblock Symmetric form of the magnetohydrodynamic equation.
\newblock {\em Numerical Methods for Mechanics of Continuum Medium} {\bf 1972},
  {\em 3},~26--34.

\bibitem[Godunov and Romenski(1972)]{GodunovRomenski72}
Godunov, S.; Romenski, E.
\newblock Nonstationary equations of the nonlinear theory of elasticity in
  {Euler} coordinates.
\newblock {\em Journal of Applied Mechanics and Technical Physics} {\bf 1972},
  {\em 13},~868--885.

\bibitem[Romenski(1998)]{Rom1998}
Romenski, E.
\newblock Hyperbolic systems of thermodynamically compatible conservation laws
  in continuum mechanics.
\newblock {\em Mathematical and Computer Modelling} {\bf 1998}, {\em
  28(10)},~115--130.

\bibitem[Peshkov and Romenski(2016)]{PeshRom2014}
Peshkov, I.; Romenski, E.
\newblock A hyperbolic model for viscous {Newtonian} flows.
\newblock {\em Continuum Mechanics and Thermodynamics} {\bf 2016}, {\em
  28},~85--104.

\bibitem[Dumbser \em{et~al.}(2016)Dumbser, Peshkov, Romenski, and
  Zanotti]{HPRmodel}
Dumbser, M.; Peshkov, I.; Romenski, E.; Zanotti, O.
\newblock {High order ADER schemes for a unified first order hyperbolic
  formulation of continuum mechanics: Viscous heat-conducting fluids and
  elastic solids}.
\newblock {\em Journal of Computational Physics} {\bf 2016}, {\em
  314},~824--862.

\bibitem[Dumbser \em{et~al.}(2017)Dumbser, Peshkov, Romenski, and
  Zanotti]{HPRmodelMHD}
Dumbser, M.; Peshkov, I.; Romenski, E.; Zanotti, O.
\newblock {H}igh order {ADER} schemes for a unified first order hyperbolic
  formulation of {N}ewtonian continuum mechanics coupled with electro-dynamics.
\newblock {\em Journal of Computational Physics} {\bf 2017}, {\em
  348},~298--342.

\bibitem[Boscheri \em{et~al.}(2016)Boscheri, Dumbser, and
  Loub\`ere]{LagrangeHPR}
Boscheri, W.; Dumbser, M.; Loub\`ere, R.
\newblock {Cell centered direct Arbitrary-Lagrangian-Eulerian ADER-WENO finite
  volume schemes for nonlinear hyperelasticity}.
\newblock {\em Computers and Fluids} {\bf 2016}, {\em 134-135},~111--129.

\bibitem[Neumann and Richtmyer(1950)]{vonneumann50}
Neumann, J.V.; Richtmyer, R.D.
\newblock A method for the numerical calculation of hydrodynamic shocks.
\newblock {\em Journal of Applied Physics} {\bf 1950}, {\em 21},~232--237.

\bibitem[{S.K.}Godunov(1959)]{godunov:linear}
{S.K.}Godunov.
\newblock A finite difference Method for the Computation of discontinuous
  solutions of the equations of fluid dynamics.
\newblock {\em Mathematics of the USSR, Sbornik} {\bf 1959}, {\em
  47},~357--393.

\bibitem[Kolgan(1972)]{Kolgan:1972}
Kolgan, V.P.
\newblock Application of the minimum-derivative principle in the construction
  of finite-difference schemes for numerical analysis of discontinuous
  solutions in gas dynamics.
\newblock {\em Transactions of the Central Aerohydrodynamics Institute} {\bf
  1972}, {\em 3},~68--77.
\newblock in Russian.

\bibitem[{van Leer}(1974)]{leer2}
{van Leer}, B.
\newblock Towards the Ultimate Conservative Difference Scheme {II}:
  Monotonicity and conservation combined in a second order scheme.
\newblock {\em Journal of Computational Physics} {\bf 1974}, {\em
  14},~361--370.

\bibitem[{van Leer}(1979)]{leer5}
{van Leer}, B.
\newblock Towards the Ultimate Conservative Difference Scheme {V}: A second
  Order sequel to {Godunov}'s Method.
\newblock {\em Journal of Computational Physics} {\bf 1979}, {\em
  32},~101--136.

\bibitem[Harten and Osher(1987)]{HartenENO}
Harten, A.; Osher, S.
\newblock Uniformly high–order accurate nonoscillatory schemes I.
\newblock {\em SIAM J. Num. Anal.} {\bf 1987}, {\em 24},~279--309.

\bibitem[Jiang and Shu(1996)]{shu_efficient_weno}
Jiang, G.; Shu, C.
\newblock Efficient Implementation of Weighted {ENO} Schemes.
\newblock {\em Journal of Computational Physics} {\bf 1996}, {\em
  126},~202--228.

\bibitem[Cockburn and Shu(1989)]{cbs1}
Cockburn, B.; Shu, C.W.
\newblock {TVB} {Runge}-{Kutta} local projection discontinuous {Galerkin}
  finite element method for conservation laws {II}: general framework.
\newblock {\em Mathematics of Computation} {\bf 1989}, {\em 52},~411--435.

\bibitem[Cockburn \em{et~al.}(1989)Cockburn, Lin, and Shu]{cbs2}
Cockburn, B.; Lin, S.Y.; Shu, C.
\newblock {TVB} {Runge}-{Kutta} local projection discontinuous {Galerkin}
  finite element method for conservation laws {III}: one dimensional systems.
\newblock {\em Journal of Computational Physics} {\bf 1989}, {\em 84},~90--113.

\bibitem[Cockburn \em{et~al.}(1990)Cockburn, Hou, and Shu]{cbs3}
Cockburn, B.; Hou, S.; Shu, C.W.
\newblock The {Runge}-{Kutta} local projection discontinuous {Galerkin} finite
  element method for conservation laws {IV}: the multidimensional case.
\newblock {\em Mathematics of Computation} {\bf 1990}, {\em 54},~545--581.

\bibitem[Cockburn and Shu(1998{\natexlab{a}})]{cbs4}
Cockburn, B.; Shu, C.W.
\newblock The {Runge}-{Kutta} discontinuous {Galerkin} method for conservation
  laws {V}: multidimensional systems.
\newblock {\em Journal of Computational Physics} {\bf 1998}, {\em
  141},~199--224.

\bibitem[Cockburn and Shu(1998{\natexlab{b}})]{CBS-convection-diffusion}
Cockburn, B.; Shu, C.W.
\newblock The Local Discontinuous {Galerkin} Method for Time-Dependent
  Convection Diffusion Systems.
\newblock {\em SIAM Journal on Numerical Analysis} {\bf 1998}, {\em
  35},~2440--2463.

\bibitem[Cockburn and Shu(2001)]{CBS-convection-dominated}
Cockburn, B.; Shu, C.W.
\newblock {Runge}-{Kutta} Discontinuous {Galerkin} Methods for
  Convection-Dominated Problems.
\newblock {\em Journal of Scientific Computing} {\bf 2001}, {\em 16},~173--261.

\bibitem[Shu(2016)]{ShuDGWENOReview}
Shu, C.
\newblock {High order WENO and DG methods for time-dependent
  convection-dominated PDEs: A brief survey of several recent developments}.
\newblock {\em Journal of Computational Phyiscs} {\bf 2016}, {\em
  316},~598--613.

\bibitem[Dumbser \em{et~al.}(2008)Dumbser, Balsara, Toro, and
  Munz]{Dumbser2008}
Dumbser, M.; Balsara, D.; Toro, E.; Munz, C.
\newblock A Unified Framework for the Construction of One-Step Finite--Volume
  and discontinuous {Galerkin} schemes.
\newblock {\em Journal of Computational Physics} {\bf 2008}, {\em
  227},~8209--8253.

\bibitem[Dumbser \em{et~al.}(2009)Dumbser, Castro, {Par\'es}, and Toro]{ADERNC}
Dumbser, M.; Castro, M.; {Par\'es}, C.; Toro, E.
\newblock {ADER} Schemes on Unstructured Meshes for Non--Conservative
  Hyperbolic Systems: Applications to Geophysical Flows.
\newblock {\em Computers and Fluids} {\bf 2009}, {\em 38},~1731--–1748.

\bibitem[Dumbser \em{et~al.}(2014)Dumbser, Zanotti, Loub{\`{e}}re, and
  Diot]{DGLimiter1}
Dumbser, M.; Zanotti, O.; Loub{\`{e}}re, R.; Diot, S.
\newblock {A posteriori subcell limiting of the discontinuous Galerkin finite
  element method for hyperbolic conservation laws}.
\newblock {\em Journal of Computational Physics} {\bf 2014}, {\em 278},~47--75.

\bibitem[Zanotti \em{et~al.}(2015)Zanotti, Fambri, Dumbser, and
  Hidalgo]{DGLimiter2}
Zanotti, O.; Fambri, F.; Dumbser, M.; Hidalgo, A.
\newblock {Space--time adaptive ADER discontinuous Galerkin finite element
  schemes with a posteriori sub--cell finite volume limiting}.
\newblock {\em Computers and Fluids} {\bf 2015}, {\em 118},~204--224.

\bibitem[Dumbser and Loub{\`{e}}re(2016)]{DGLimiter3}
Dumbser, M.; Loub{\`{e}}re, R.
\newblock {A simple robust and accurate a posteriori sub--cell finite volume
  limiter for the discontinuous Galerkin method on unstructured meshes}.
\newblock {\em Journal of Computational Physics} {\bf 2016}, {\em
  319},~163--199.

\bibitem[Titarev and Toro(2002)]{toro3}
Titarev, V.; Toro, E.
\newblock {ADER}: Arbitrary High Order {Godunov} Approach.
\newblock {\em Journal of Scientific Computing} {\bf 2002}, {\em 17},~609--618.

\bibitem[Toro and Titarev(2002)]{toro4}
Toro, E.; Titarev, V.
\newblock Solution of the generalized {Riemann} problem for advection-reaction
  equations.
\newblock {\em Proc. Roy. Soc. London} {\bf 2002}, {\em 458},~271--281.

\bibitem[Titarev and Toro(2005)]{titarevtoro}
Titarev, V.; Toro, E.
\newblock {ADER} schemes for three-dimensional nonlinear hyperbolic systems.
\newblock {\em Journal of Computational Physics} {\bf 2005}, {\em
  204},~715--736.

\bibitem[Toro and Titarev(2006)]{Toro:2006a}
Toro, E.F.; Titarev, V.A.
\newblock {Derivative Riemann solvers for systems of conservation laws and ADER
  methods}.
\newblock {\em Journal of Computational Physics} {\bf 2006}, {\em
  212},~150--165.

\bibitem[Bungartz \em{et~al.}(2010)Bungartz, Mehl, Neckel, and
  Weinzierl]{Peano1}
Bungartz, H.; Mehl, M.; Neckel, T.; Weinzierl, T.
\newblock {The PDE framework Peano applied to fluid dynamics: An efficient
  implementation of a parallel multiscale fluid dynamics solver on octree-like
  adaptive Cartesian grids}.
\newblock {\em Computational Mechanics} {\bf 2010}, {\em 46},~103--114.

\bibitem[Weinzierl and Mehl(2011)]{Peano2}
Weinzierl, T.; Mehl, M.
\newblock {Peano-A traversal and storage scheme for octree-like adaptive
  Cartesian multiscale grids}.
\newblock {\em SIAM Journal on Scientific Computing} {\bf 2011}, {\em
  33},~2732--2760.

\bibitem[Khokhlov(1998)]{Khokhlov1998}
Khokhlov, A.
\newblock Fully Threaded Tree Algorithms for Adaptive Refinement Fluid Dynamics
  Simulations.
\newblock {\em Journal of Computational Physics} {\bf 1998}, {\em
  143},~519--543.

\bibitem[Dumbser \em{et~al.}(2013)Dumbser, Zanotti, Hidalgo, and
  Balsara]{AMR3DCL}
Dumbser, M.; Zanotti, O.; Hidalgo, A.; Balsara, D.
\newblock {ADER-WENO Finite Volume Schemes with Space-Time Adaptive Mesh
  Refinement}.
\newblock {\em Journal of Computational Physics} {\bf 2013}, {\em
  248},~257--286.

\bibitem[Dumbser \em{et~al.}(2014)Dumbser, Hidalgo, and Zanotti]{AMR3DNC}
Dumbser, M.; Hidalgo, A.; Zanotti, O.
\newblock {High Order Space-Time Adaptive ADER-WENO Finite Volume Schemes for
  Non-Conservative Hyperbolic Systems}.
\newblock {\em Computer Methods in Applied Mechanics and Engineering} {\bf
  2014}, {\em 268},~359--387.

\bibitem[Zanotti \em{et~al.}(2015)Zanotti, Fambri, Dumbser, and
  Hidalgo]{Zanotti2015c}
Zanotti, O.; Fambri, F.; Dumbser, M.; Hidalgo, A.
\newblock Space-time adaptive {ADER} discontinuous {{G}alerkin} finite element
  schemes with a posteriori sub-cell finite volume limiting.
\newblock {\em Computers and Fluids} {\bf 2015}, {\em 118},~204 -- 224.

\bibitem[Zanotti \em{et~al.}({2015})Zanotti, Fambri, and Dumbser]{Zanotti2015d}
Zanotti, O.; Fambri, F.; Dumbser, M.
\newblock Solving the relativistic magnetohydrodynamics equations with {ADER}
  discontinuous {G}alerkin methods, a posteriori subcell limiting and adaptive
  mesh refinement.
\newblock {\em Mon. Not. R. Astron. Soc.} {\bf {2015}}, {\em 452},~3010--3029.

\bibitem[Fambri \em{et~al.}(2017)Fambri, Dumbser, and Zanotti]{ADERDGVisc}
Fambri, F.; Dumbser, M.; Zanotti, O.
\newblock Space-time adaptive {ADER}-{DG} schemes for dissipative flows:
  {C}ompressible {N}avier-{S}tokes and resistive {MHD} equations.
\newblock {\em Computer Physics Communications} {\bf 2017}, {\em
  220},~297--318.

\bibitem[Fambri \em{et~al.}(2018)Fambri, Dumbser, K\"oppel, Rezzolla, and
  Zanotti]{ADERGRMHD}
Fambri, F.; Dumbser, M.; K\"oppel, S.; Rezzolla, L.; Zanotti, O.
\newblock {ADER discontinuous Galerkin schemes for general-relativistic ideal
  magnetohydrodynamics}.
\newblock {\em Monthly Notices of the Royal Astronomical Society} {\bf 2018},
  {\em 477},~4543--4564.

\bibitem[{Berger} and {Oliger}(1984)]{Berger-Oliger1984}
{Berger}, M.J.; {Oliger}, J.
\newblock {Adaptive Mesh Refinement for Hyperbolic Partial Differential
  Equations}.
\newblock {\em Journal of Computational Physics} {\bf 1984}, {\em
  53},~484--512.

\bibitem[{Berger} and {Jameson}(1985)]{Berger-Jameson1985}
{Berger}, M.J.; {Jameson}, A.
\newblock {Automatic adaptive grid refinement for the Euler equations}.
\newblock {\em AIAA Journal} {\bf 1985}, {\em 23},~561--568.

\bibitem[{Berger} and {Colella}(1989)]{Berger-Colella1989}
{Berger}, M.J.; {Colella}, P.
\newblock {Local adaptive mesh refinement for shock hydrodynamics}.
\newblock {\em Journal of Computational Physics} {\bf 1989}, {\em 82},~64--84.

\bibitem[software()]{LeVequeCLAWPACK}
software, C.
\newblock Available at http://depts.washington.edu/clawpack/.

\bibitem[{Berger} and {LeVeque}(1998)]{BergerLeveque1998}
{Berger}, M.J.; {LeVeque}, R.
\newblock {Adaptive mesh refinement using wave-propagation algorithms for
  hyperbolic systems}.
\newblock {\em SIAM Journal on Numerical Analysis} {\bf 1998}, {\em
  35},~2298--2316.

\bibitem[Bell \em{et~al.}(1994)Bell, Berger, Saltzman, and Welcome]{Bell1994}
Bell, J.; Berger, M.; Saltzman, J.; Welcome, M.
\newblock Three-dimensional adaptive mesh refinement for hyperbolic
  conservation laws.
\newblock {\em SIAM Journal on Scientific Computing} {\bf 1994}, {\em
  15},~127--138.

\bibitem[{Dezeeuw} and {Powell}(1993)]{Dezeeuw1993}
{Dezeeuw}, D.; {Powell}, K.G.
\newblock {An Adaptively Refined Cartesian Mesh Solver for the Euler
  Equations}.
\newblock {\em Journal of Computational Physics} {\bf 1993}, {\em 104},~56--68.

\bibitem[Balsara(2001)]{BalsaraAMR}
Balsara, D.
\newblock Divergence-free adaptive mesh refinement for magnetohydrodynamics.
\newblock {\em Journal of Computational Physics} {\bf 2001}, {\em
  174},~614--648.

\bibitem[{Teyssier}(2002)]{Teyssier2002}
{Teyssier}, R.
\newblock {Cosmological hydrodynamics with adaptive mesh refinement. A new high
  resolution code called RAMSES}.
\newblock {\em Astronomy \& Astrophysics} {\bf 2002}, {\em 385},~337--364.

\bibitem[{Keppens} \em{et~al.}(2003){Keppens}, {Nool}, {T{\'o}th}, and
  {Goedbloed}]{Keppens2003}
{Keppens}, R.; {Nool}, M.; {T{\'o}th}, G.; {Goedbloed}, J.P.
\newblock {Adaptive Mesh Refinement for conservative systems: multi-dimensional
  efficiency evaluation}.
\newblock {\em Computer Physics Communications} {\bf 2003}, {\em
  153},~317--339.

\bibitem[Ziegler(2008)]{Ziegler2008}
Ziegler, U.
\newblock The NIRVANA code: Parallel computational MHD with adaptive mesh
  refinement.
\newblock {\em Computer Physics Communications} {\bf 2008}, {\em
  179},~227--244.

\bibitem[{Mignone} \em{et~al.}(2012){Mignone}, {Zanni}, {Tzeferacos}, {van
  Straalen}, {Colella}, and {Bodo}]{Mignone2012}
{Mignone}, A.; {Zanni}, C.; {Tzeferacos}, P.; {van Straalen}, B.; {Colella},
  P.; {Bodo}, G.
\newblock {The PLUTO Code for Adaptive Mesh Computations in Astrophysical Fluid
  Dynamics}.
\newblock {\em The Astrophysical Journal Supplement Series} {\bf 2012}, {\em
  198},~7.

\bibitem[Cunningham \em{et~al.}(2009)Cunningham, Frank, Varni\`ere, Mitran, and
  Jones]{Cunningham2009}
Cunningham, A.; Frank, A.; Varni\`ere, P.; Mitran, S.; Jones, T.W.
\newblock Simulating Magnetohydrodynamical Flow with Constrained Transport and
  Adaptive Mesh Refinement: Algorithms and Tests of the AstroBEAR Code.
\newblock {\em The Astrophysical Journal Supplement Series} {\bf 2009}, {\em
  182},~519.

\bibitem[Keppens \em{et~al.}(2012)Keppens, Meliani, van Marle, Delmont, Vlasis,
  and van~der Holst]{Keppens2012}
Keppens, R.; Meliani, Z.; van Marle, A.; Delmont, P.; Vlasis, A.; van~der
  Holst, B.
\newblock Parallel, grid-adaptive approaches for relativistic hydro and
  magnetohydrodynamics.
\newblock {\em Journal of Computational Physics} {\bf 2012}, {\em
  231},~718--744.

\bibitem[et~al.(2017)]{BHAC}
et~al., O.P.
\newblock The black hole accretion code.
\newblock {\em Computational Astrophysics and Cosmology} {\bf 2017}, {\em
  4:1},~1--42.

\bibitem[et~al.(2014)]{AMRSurvey}
et~al., A.D.
\newblock A survey of high level frameworks in block-structured adaptive mesh
  refinement packages.
\newblock {\em Journal of Parallel and Distributed Computing} {\bf 2014}, {\em
  74},~3217--3227.

\bibitem[Baeza and Mulet(2006)]{Mulet1}
Baeza, A.; Mulet, P.
\newblock Adaptive mesh refinement techniques for high--order shock capturing
  schemes for multi--dimensional hydrodynamic simulations.
\newblock {\em International Journal for Numerical Methods in Fluids} {\bf
  2006}, {\em 52},~455--471.

\bibitem[{Colella} \em{et~al.}(2009){Colella}, {Dorr}, {Hittinger}, {Martin},
  and {McCorquodale}]{Colella2009}
{Colella}, P.; {Dorr}, M.; {Hittinger}, J.; {Martin}, D.F.; {McCorquodale}, P.
\newblock {High-order finite-volume adaptive methods on locally rectangular
  grids}.
\newblock {\em Journal of Physics Conference Series} {\bf 2009}, {\em
  180},~012010.

\bibitem[B\"urger \em{et~al.}(2013)B\"urger, Mulet, and Villada]{Burger2012}
B\"urger, R.; Mulet, P.; Villada, L.
\newblock {Spectral WENO schemes with Adaptive Mesh Refinement for models of
  polydisperse sedimentation}.
\newblock {\em ZAMM -- Journal of Applied Mathematics and Mechanics /
  Zeitschrift f?ewandte Mathematik und Mechanik} {\bf 2013}, {\em
  93},~373--386.

\bibitem[Ivan and Groth(2014)]{Ivan2014}
Ivan, L.; Groth, C.
\newblock {High-order solution-adaptive central essentially non-oscillatory
  (CENO) method for viscous flows}.
\newblock {\em Journal of Computational Physics} {\bf 2014}, {\em
  257},~830--862.

\bibitem[Buchm\"uller \em{et~al.}(2016)Buchm\"uller, Dreher, and
  Helzel]{Buchmuller2015}
Buchm\"uller, P.; Dreher, J.; Helzel, C.
\newblock {Finite volume WENO methods for hyperbolic conservation laws on
  Cartesian grids with adaptive mesh refinement}.
\newblock {\em Applied Mathematics and Computation} {\bf 2016}, {\em
  272},~460--478.

\bibitem[Semplice \em{et~al.}(2016)Semplice, Coco, and
  Russo]{SCR:CWENOquadtree}
Semplice, M.; Coco, A.; Russo, G.
\newblock Adaptive Mesh Refinement for Hyperbolic Systems based on Third-Order
  Compact {WENO} Reconstruction.
\newblock {\em Journal of Scientific Computing} {\bf 2016}, {\em 66},~692--724.

\bibitem[Shen \em{et~al.}(2011)Shen, Qiu, and Christlieb]{FDWENOAMR}
Shen, C.; Qiu, J.; Christlieb, A.
\newblock {Adaptive mesh refinement based on high order finite difference WENO
  scheme for multi-scale simulations}.
\newblock {\em Journal of Computational Physics} {\bf 2011}, {\em
  230},~3780--3802.

\bibitem[Stroud(1971)]{stroud}
Stroud, A.
\newblock {\em Approximate Calculation of Multiple Integrals}; Prentice-Hall
  Inc.: Englewood Cliffs, New Jersey,  1971.

\bibitem[Castro \em{et~al.}(2006)Castro, Gallardo, and Par\'es]{Castro2006}
Castro, M.; Gallardo, J.; Par\'es, C.
\newblock High-order finite volume schemes based on reconstruction of states
  for solving hyperbolic systems with nonconservative products. Applications to
  shallow-water systems.
\newblock {\em Mathematics of Computation} {\bf 2006}, {\em 75},~1103--1134.

\bibitem[Par\'es(2006)]{Pares2006}
Par\'es, C.
\newblock Numerical methods for nonconservative hyperbolic systems: a
  theoretical framework.
\newblock {\em SIAM Journal on Numerical Analysis} {\bf 2006}, {\em
  44},~300--321.

\bibitem[Rhebergen \em{et~al.}(2008)Rhebergen, Bokhove, and van~der
  Vegt]{Rhebergen2008}
Rhebergen, S.; Bokhove, O.; van~der Vegt, J.
\newblock Discontinuous {Galerkin} finite element methods for hyperbolic
  nonconservative partial differential equations.
\newblock {\em Journal of Computational Physics} {\bf 2008}, {\em
  227},~1887--1922.

\bibitem[Dumbser \em{et~al.}(2010)Dumbser, Hidalgo, Castro, Par\'es, and
  Toro]{USFORCE2}
Dumbser, M.; Hidalgo, A.; Castro, M.; Par\'es, C.; Toro, E.
\newblock {FORCE} Schemes on Unstructured Meshes {II}: Non--Conservative
  Hyperbolic Systems.
\newblock {\em Computer Methods in Applied Mechanics and Engineering} {\bf
  2010}, {\em 199},~625--647.

\bibitem[M{\"{u}}ller \em{et~al.}(2013)M{\"{u}}ller, Par{\'{e}}s, and
  Toro]{MuellerToro1}
M{\"{u}}ller, L.O.; Par{\'{e}}s, C.; Toro, E.F.
\newblock {Well-balanced high-order numerical schemes for one-dimensional blood
  flow in vessels with varying mechanical properties}.
\newblock {\em Journal of Computational Physics} {\bf 2013}, {\em 242},~53--85.

\bibitem[M\"uller and Toro(2013)]{MuellerToro2}
M\"uller, L.; Toro, E.
\newblock Well-balanced high-order solver for blood flow in networks of vessels
  with variable properties.
\newblock {\em International Journal for Numerical Methods in Biomedical
  Engineering} {\bf 2013}, {\em 29},~1388--1411.

\bibitem[Gaburro \em{et~al.}(2017)Gaburro, Dumbser, and
  Castro]{GaburroDumbserSWE}
Gaburro, E.; Dumbser, M.; Castro, M.
\newblock {Direct Arbitrary-Lagrangian-Eulerian finite volume schemes on moving
  nonconforming unstructured meshes}.
\newblock {\em Computers and Fluids} {\bf 2017}, {\em 159},~254--275.

\bibitem[Gaburro \em{et~al.}(2018)Gaburro, Castro, and
  Dumbser]{GaburroDumbserEuler}
Gaburro, E.; Castro, M.; Dumbser, M.
\newblock {Well balanced Arbitrary-Lagrangian-Eulerian finite volume schemes on
  moving nonconforming meshes for the Euler equations of gasdynamics with
  gravity}.
\newblock {\em Monthly Notices of the Royal Astronomical Society} {\bf 2018},
  {\em 477},~2251--2275.

\bibitem[Toro \em{et~al.}(2009)Toro, Hidalgo, and Dumbser]{USFORCE}
Toro, E.; Hidalgo, A.; Dumbser, M.
\newblock {FORCE} Schemes on Unstructured Meshes {I}: Conservative Hyperbolic
  Systems.
\newblock {\em Journal of Computational Physics} {\bf 2009}, {\em
  228},~3368--3389.

\bibitem[Dumbser and Toro(2011)]{OsherNC}
Dumbser, M.; Toro, E.F.
\newblock A Simple Extension of the {Osher} {Riemann} Solver to
  Non-Conservative Hyperbolic Systems.
\newblock {\em Journal of Scientific Computing} {\bf 2011}, {\em 48},~70--88.

\bibitem[Castro \em{et~al.}(2010)Castro, Pardo, Par\'es, and
  Toro]{CastroPardoPares}
Castro, M.; Pardo, A.; Par\'es, C.; Toro, E.
\newblock On some fast well-balanced first order solvers for nonconservative
  systems.
\newblock {\em Mathematics of Computation} {\bf 2010}, {\em 79},~1427--1472.

\bibitem[Dumbser and Toro(2011)]{OsherUniversal}
Dumbser, M.; Toro, E.F.
\newblock On Universal {Osher}--Type Schemes for General Nonlinear Hyperbolic
  Conservation Laws.
\newblock {\em Communications in Computational Physics} {\bf 2011}, {\em
  10},~635--671.

\bibitem[{Einfeldt} \em{et~al.}(1991){Einfeldt}, {Roe}, {Munz}, and
  {Sjogreen}]{Einfeldt1991}
{Einfeldt}, B.; {Roe}, P.L.; {Munz}, C.D.; {Sjogreen}, B.
\newblock {On Godunov-type methods near low densities}.
\newblock {\em Journal of Computational Physics} {\bf 1991}, {\em
  92},~273--295.

\bibitem[Dumbser and Balsara(2016)]{HLLEMNC}
Dumbser, M.; Balsara, D.
\newblock A New, Efficient Formulation of the {HLLEM Riemann} Solver for
  General Conservative and Non-Conservative Hyperbolic Systems.
\newblock {\em Journal of Computational Physics} {\bf 2016}, {\em
  304},~275--319.

\bibitem[Dumbser \em{et~al.}(2008)Dumbser, Enaux, and Toro]{DumbserEnauxToro}
Dumbser, M.; Enaux, C.; Toro, E.
\newblock Finite Volume Schemes of Very High Order of Accuracy for Stiff
  Hyperbolic Balance Laws.
\newblock {\em Journal of Computational Physics} {\bf 2008}, {\em
  227},~3971--4001.

\bibitem[Dumbser and Zanotti(2009)]{DumbserZanotti}
Dumbser, M.; Zanotti, O.
\newblock Very High Order {PNPM} Schemes on Unstructured Meshes for the
  Resistive Relativistic {MHD} Equations.
\newblock {\em Journal of Computational Physics} {\bf 2009}, {\em
  228},~6991--7006.

\bibitem[Harten \em{et~al.}(1987)Harten, Engquist, Osher, and
  Chakravarthy]{eno}
Harten, A.; Engquist, B.; Osher, S.; Chakravarthy, S.
\newblock Uniformly high order essentially non-oscillatory schemes, {III}.
\newblock {\em Journal of Computational Physics} {\bf 1987}, {\em
  71},~231--303.

\bibitem[Dumbser \em{et~al.}(2007)Dumbser, K\"aser, Titarev, and
  Toro]{DumbserKaeser07}
Dumbser, M.; K\"aser, M.; Titarev, V.; Toro, E.
\newblock Quadrature-Free Non-Oscillatory Finite Volume Schemes on Unstructured
  Meshes for Nonlinear Hyperbolic Systems.
\newblock {\em Journal of Computational Physics} {\bf 2007}, {\em
  226},~204--243.

\bibitem[Taube \em{et~al.}(2007)Taube, Dumbser, Balsara, and Munz]{taube_jsc}
Taube, A.; Dumbser, M.; Balsara, D.; Munz, C.
\newblock Arbitrary High Order Discontinuous {Galerkin} Schemes for the
  Magnetohydrodynamic Equations.
\newblock {\em Journal of Scientific Computing} {\bf 2007}, {\em 30},~441--464.

\bibitem[Hidalgo and Dumbser(2011)]{HidalgoDumbser}
Hidalgo, A.; Dumbser, M.
\newblock {ADER} Schemes for Nonlinear Systems of Stiff
  Advection-Diffusion-Reaction Equations.
\newblock {\em Journal of Scientific Computing} {\bf 2011}, {\em 48},~173--189.

\bibitem[Jackson(2017)]{Jackson}
Jackson, H.
\newblock On the eigenvalues of the ADER-WENO Galerkin predictor.
\newblock {\em Journal of Computational Physics} {\bf 2017}, {\em
  333},~409--413.

\bibitem[Zanotti and Dumbser(2016)]{ADERPrim}
Zanotti, O.; Dumbser, M.
\newblock Efficient conservative ADER schemes based on WENO reconstruction and
  space-time predictor in primitive variables.
\newblock {\em Computational Astrophysics and Cosmology} {\bf 2016}, {\em
  3},~1.

\bibitem[Owren and Zennaro(1992)]{OwrenZennaro}
Owren, B.; Zennaro, M.
\newblock {Derivation of efficient, continuous, explicit Runge--Kutta methods}.
\newblock {\em SIAM J. Sci. and Stat. Comput.} {\bf 1992}, {\em
  13},~1488--1501.

\bibitem[Gassner \em{et~al.}(2011)Gassner, Dumbser, Hindenlang, and
  Munz]{Gassner2011a}
Gassner, G.; Dumbser, M.; Hindenlang, F.; Munz, C.
\newblock {Explicit one--step time discretizations for discontinuous {G}alerkin
  and finite volume schemes based on local predictors}.
\newblock {\em Journal of Computational Physics} {\bf 2011}, {\em
  230},~4232--4247.

\bibitem[Heinecke \em{et~al.}(2015)Heinecke, Pabst, and Henry]{libxsmm}
Heinecke, A.; Pabst, H.; Henry, G.
\newblock {LIBXSMM: A High Performance Library for Small Matrix
  Multiplications}.
\newblock Technical report, SC'15: The International Conference for High
  Performance Computing, Networking, Storage and Analysis, Austin (Texas),
  2015.
\newblock https://github.com/hfp/libxsmm.

\bibitem[Breuer \em{et~al.}(2014{\natexlab{a}})Breuer, Heinecke, Bader, and
  Pelties]{SeisSol1}
Breuer, A.; Heinecke, A.; Bader, M.; Pelties, C.
\newblock {Accelerating SeisSol by generating vectorized code for sparse matrix
  operators}.
\newblock {\em Advances in Parallel Computing} {\bf 2014}, {\em 25},~347--356.

\bibitem[Breuer \em{et~al.}(2014{\natexlab{b}})Breuer, Heinecke, Rettenberger,
  Bader, Gabriel, and Pelties]{SeisSol2}
Breuer, A.; Heinecke, A.; Rettenberger, S.; Bader, M.; Gabriel, A.; Pelties, C.
\newblock {Sustained petascale performance of seismic simulations with SeisSol
  on SuperMUC}.
\newblock {\em Lecture Notes in Computer Science (LNCS)} {\bf 2014}, {\em
  8488},~1--18.

\bibitem[Charrier and Weinzierl(2018)]{StopTalking}
Charrier, D.; Weinzierl, T.
\newblock {Stop talking to me -- a communication-avoiding ADER-DG realisation}.
\newblock {\em SIAM Journal on Scientific Computing} {\bf 2018}.
\newblock submitted to. https://arxiv.org/abs/1801.08682.

\bibitem[Boscheri and Dumbser(2017)]{ALEDG}
Boscheri, W.; Dumbser, M.
\newblock {Arbitrary--Lagrangian--Eulerian Discontinuous Galerkin schemes with
  a posteriori subcell finite volume limiting on moving unstructured meshes}.
\newblock {\em Journal of Computational Physics} {\bf 2017}, {\em
  346},~449--479.

\bibitem[Clain \em{et~al.}(2011)Clain, Diot, and Loub{\`e}re]{CDL1}
Clain, S.; Diot, S.; Loub{\`e}re, R.
\newblock A high-order finite volume method for systems of conservation
  laws—Multi-dimensional Optimal Order Detection ({MOOD}).
\newblock {\em Journal of Computational Physics} {\bf 2011}, {\em 230},~4028 --
  4050.

\bibitem[Diot \em{et~al.}(2012)Diot, Clain, and Loub{\`e}re]{CDL2}
Diot, S.; Clain, S.; Loub{\`e}re, R.
\newblock Improved detection criteria for the Multi-dimensional Optimal Order
  Detection ({MOOD}) on unstructured meshes with very high-order polynomials.
\newblock {\em Computers and Fluids} {\bf 2012}, {\em 64},~43 -- 63.

\bibitem[Diot \em{et~al.}(2013)Diot, Loub{\`e}re, and Clain]{CDL3}
Diot, S.; Loub{\`e}re, R.; Clain, S.
\newblock The {MOOD} method in the three-dimensional case: Very-High-Order
  Finite Volume Method for Hyperbolic Systems.
\newblock {\em International Journal of Numerical Methods in Fluids} {\bf
  2013}, {\em 73},~362--392.

\bibitem[Loub{\`e}re \em{et~al.}(2014)Loub{\`e}re, Dumbser, and
  Diot]{ADER_MOOD_14}
Loub{\`e}re, R.; Dumbser, M.; Diot, S.
\newblock A New Family of High Order Unstructured MOOD and ADER Finite Volume
  Schemes for Multidimensional Systems of Hyperbolic Conservation Laws.
\newblock {\em Communication in Computational Physics} {\bf 2014}, {\em
  16},~718--763.

\bibitem[Levy \em{et~al.}(1999)Levy, Puppo, and Russo]{LPR:99}
Levy, D.; Puppo, G.; Russo, G.
\newblock Central {WENO} schemes for hyperbolic systems of conservation laws.
\newblock {\em M2AN Math. Model. Numer. Anal.} {\bf 1999}, {\em 33},~547--571.

\bibitem[Levy \em{et~al.}(2000)Levy, Puppo, and Russo]{LPR:2001}
Levy, D.; Puppo, G.; Russo, G.
\newblock Compact central {WENO} schemes for multidimensional conservation
  laws.
\newblock {\em SIAM J. Sci. Comput.} {\bf 2000}, {\em 22},~656--672.

\bibitem[Dumbser \em{et~al.}(2017)Dumbser, Boscheri, Semplice, and
  Russo]{ADER_CWENO}
Dumbser, M.; Boscheri, W.; Semplice, M.; Russo, G.
\newblock Central weighted {ENO} schemes for hyperbolic conservation laws on
  fixed and moving unstructured meshes.
\newblock {\em {SIAM} J. Sci. Comput.} {\bf 2017}, {\em 39},~A2564--A2591.

\bibitem[Hu and Shu(1999)]{HuShuTri}
Hu, C.; Shu, C.
\newblock Weighted essentially non-oscillatory schemes on triangular meshes.
\newblock {\em Journal of Computational Physics} {\bf 1999}, {\em
  150},~97--127.

\bibitem[Sedov(1959)]{Sedov}
Sedov, L.
\newblock {\em Similarity and Dimensional Methods in Mechanics}; Academic
  Press: New York,  1959.

\bibitem[Kamm and Timmes(2007)]{SedovExact}
Kamm, J.; Timmes, F.
\newblock On efficient generation of numerically robust Sedov solutions.
\newblock {\em Technical Report LA-UR-07-2849,LANL} {\bf 2007}.

\bibitem[Tavelli \em{et~al.}(2018)Tavelli, Dumbser, Charrier, Rannabauer,
  Weinzierl, and Bader]{AMRDIM}
Tavelli, M.; Dumbser, M.; Charrier, D.; Rannabauer, L.; Weinzierl, T.; Bader,
  M.
\newblock {A simple diffuse interface approach on adaptive Cartesian grids for
  the linear elastic wave equations with complex topography}.
\newblock {\em Journal of Computational Physics} {\bf 2018}.
\newblock submitted to. https://arxiv.org/abs/1804.09491.

\bibitem[Dumbser and K\"aser(2006)]{gij2}
Dumbser, M.; K\"aser, M.
\newblock {An arbitrary high--order discontinuous Galerkin method for elastic
  waves on unstructured meshes -- II. The three--dimensional isotropic case}.
\newblock {\em Geophysical Journal International} {\bf 2006}, {\em
  167},~319--336.

\bibitem[Godunov \em{et~al.}(1961)Godunov, Zabrodin, and
  Prokopov]{Godunov:1961a}
Godunov, S.K.; Zabrodin, A.V.; Prokopov, G.P.
\newblock {A Difference Scheme for Two-Dimensional Unsteady Aerodynamics}.
\newblock {\em J. Comp. Math. and Math. Phys. USSR} {\bf 1961}, {\em
  2},~1020--1050.

\bibitem[Ant{\'o}n \em{et~al.}(2006)Ant{\'o}n, Zanotti, Miralles, Mart{\'\i},
  Ib{\'a}{\~n}ez, Font, and Pons]{Anton06}
Ant{\'o}n, L.; Zanotti, O.; Miralles, J.A.; Mart{\'\i}, J.M.; Ib{\'a}{\~n}ez,
  J.M.; Font, J.A.; Pons, J.A.
\newblock Numerical 3+1 general relativistic magnetohydrodynamics: a local
  characteristic approach.
\newblock {\em Astrophys. J.} {\bf 2006}, {\em 637},~296.

\bibitem[Zanna \em{et~al.}(2007)Zanna, Zanotti, Bucciantini, and
  Londrillo]{ZannaZanotti}
Zanna, L.D.; Zanotti, O.; Bucciantini, N.; Londrillo, P.
\newblock {ECHO}: an {Eulerian} Conservative High Order scheme for general
  relativistic magnetohydrodynamics and magnetodynamics.
\newblock {\em Astronomy and Astrophysics} {\bf 2007}, {\em 473},~11--30.

\bibitem[L{\"o}hner(1987)]{Loehner1987}
L{\"o}hner, R.
\newblock An adaptive finite element scheme for transient problems in {CFD}.
\newblock {\em Computer Methods in Applied Mechanics and Engineering} {\bf
  1987}, {\em 61},~323 -- 338.

\bibitem[Radice \em{et~al.}(2014)Radice, Rezzolla, and Galeazzi]{WhiskyTHC}
Radice, D.; Rezzolla, L.; Galeazzi, F.
\newblock High--order fully general--relativistic hydrodynamics: new approaches
  and tests.
\newblock {\em Classical and Quantum Gravity} {\bf 2014}, {\em 31},~075012.

\bibitem[Berm\'udez and V\'azquez(1994)]{Bermudez1994}
Berm\'udez, A.; V\'azquez, M.
\newblock Upwind methods for hyperbolic conservation laws with source terms.
\newblock {\em Computers and Fluids} {\bf 1994}, {\em 23},~1049--1071.

\bibitem[{Alic} \em{et~al.}(2012){Alic}, {Bona-Casas}, {Bona}, {Rezzolla}, and
  {Palenzuela}]{Alic:2011a}
{Alic}, D.; {Bona-Casas}, C.; {Bona}, C.; {Rezzolla}, L.; {Palenzuela}, C.
\newblock {Conformal and covariant formulation of the Z4 system with
  constraint-violation damping}.
\newblock {\em Phys. Rev. D} {\bf 2012}, {\em 85}.

\bibitem[Dedner \em{et~al.}(2002)Dedner, Kemm, Kr\"oner, Munz, Schnitzer, and
  Wesenberg]{Dedneretal}
Dedner, A.; Kemm, F.; Kr\"oner, D.; Munz, C.D.; Schnitzer, T.; Wesenberg, M.
\newblock Hyperbolic Divergence Cleaning for the {MHD} Equations.
\newblock {\em Journal of Computational Physics} {\bf 2002}, {\em
  175},~645--673.

\bibitem[Gundlach and Martin-Garcia(2004)]{Gundlach04a}
Gundlach, C.; Martin-Garcia, J.
\newblock Symmetric hyperbolic form of systems of second-order evolution
  equations subject to constraints.
\newblock {\em Phys. Rev. D} {\bf 2004}, {\em 70},~044031,
  \href{http://xxx.lanl.gov/abs/gr-qc/0402079}{{\normalfont [gr-qc/0402079]}}.

\bibitem[Dumbser \em{et~al.}(2018)Dumbser, Guercilena, K\"oppel, Rezzolla, and
  Zanotti]{ADERCCZ4}
Dumbser, M.; Guercilena, F.; K\"oppel, S.; Rezzolla, L.; Zanotti, O.
\newblock {Conformal and covariant Z4 formulation of the Einstein equations:
  strongly hyperbolic first--order reduction and solution with discontinuous
  Galerkin schemes}.
\newblock {\em Physical Review D} {\bf 2018}, {\em 97},~084053.

\end{thebibliography}




\end{document}